\RequirePackage{fix-cm}
\documentclass{svjour3}                     
\smartqed  
\usepackage[margin=1in]{geometry}
\usepackage[ansinew]{inputenc}
\usepackage[english]{babel}
\usepackage{fancyhdr}
\usepackage[active]{srcltx}
\usepackage{amsmath,amssymb}
\usepackage{xcolor}
\usepackage{array}
\usepackage{textcomp}
\usepackage{tikz}
\usepackage[ruled,vlined]{algorithm2e}
\usepackage{cuted}
\usepackage{placeins}
\usepackage[normalem]{ulem}
\usepackage{hyperref}

\def\makeheadbox{{%
\hbox to0pt{\vbox{\baselineskip=10dd\hrule\hbox
to\hsize{\vrule\kern3pt\vbox{\kern3pt
\hbox{\bfseries Nonlinear Dynamics}
\hbox{This is a pre-peer-review, pre-copyedit version of this article.}
\kern3pt}\hfil\kern3pt\vrule}\hrule}%
\hss}}}
\newcommand{\bfm}[1]{\boldsymbol{#1}}
\newcommand{\bsn}[1]{\boldsymbol{#1}}
\newcommand{\auu}{\bfm{a}}
\newcommand{\buu}{\bfm{b}}
\newcommand{\fuu}{\bfm{f}}
\newcommand{\guu}{\bfm{g}}
\newcommand{\huu}{\bfm{h}}
\newcommand{\kuu}{\bfm{k}}
\newcommand{\uuu}{\bfm{u}}
\newcommand{\Au}{\mathbf{A}}
\newcommand{\au}{\mathbf{a}}
\newcommand{\Bu}{\mathbf{B}}
\newcommand{\bu}{\mathbf{b}}

\newcommand{\cu}{\mathbf{c}}
\newcommand{\Du}{\mathbf{D}}
\newcommand{\Eu}{\mathbf{E}}
\newcommand{\eu}{\mathbf{e}}
\newcommand{\Fu}{\mathbf{F}}
\newcommand{\Gu}{\mathbf{G}}
\newcommand{\Hu}{\mathbf{H}}
\newcommand{\Iu}{\mathbf{I}}
\newcommand{\Mu}{\mathbf{M}}
\newcommand{\Ku}{\mathbf{K}}
\newcommand{\Uu}{\mathbf{U}}
\newcommand{\Vu}{\mathbf{V}}
\newcommand{\zerou}{\mathbf{0}}
\newcommand{\gu}{\mathbf{g}}
\newcommand{\hu}{\mathbf{h}}
\newcommand{\Lu}{\mathbf{L}}
\newcommand{\Pu}{\mathbf{P}}
\newcommand{\pu}{\mathbf{p}}
\newcommand{\Qu}{\mathbf{Q}}
\newcommand{\Ru}{\mathbf{R}}
\newcommand{\ru}{\mathbf{r}}
\newcommand{\Su}{\mathbf{S}}
\newcommand{\su}{\mathbf{s}}
\newcommand{\Tu}{\mathbf{T}}
\newcommand{\vu}{\mathbf{v}}

\newcommand{\uu}{\mathbf{u}}
\newcommand{\xu}{\mathbf{x}}
\newcommand{\yu}{\mathbf{y}}
\newcommand{\zu}{\mathbf{z}}
\newcommand{\Ar}{\mathrm{A}}
\newcommand{\ar}{\mathrm{a}}
\newcommand{\Br}{\mathrm{B}}
\newcommand{\br}{\mathrm{b}}
\newcommand{\gr}{\mathrm{g}}
\newcommand{\hr}{\mathrm{h}}
\newcommand{\Ir}{\mathrm{I}}
\newcommand{\Nr}{\mathrm{N}}
\newcommand{\nr}{\mathrm{n}}
\newcommand{\pr}{\mathrm{p}}
\newcommand{\Qr}{\mathrm{Q}}
\newcommand{\rr}{\mathrm{r}}
\newcommand{\sr}{\mathrm{s}}
\newcommand{\ur}{\mathrm{u}}
\newcommand{\vr}{\mathrm{v}}
\newcommand{\xr}{\mathrm{x}}
\newcommand{\zr}{\mathrm{z}}
\newcommand{\alphau}{\bsn{\alpha}}
\newcommand{\betau}{\bsn{\beta}}
\newcommand{\gammau}{\bsn{\gamma}}
\newcommand{\Deltau}{\bsn{\Delta}}
\newcommand{\Lambdau}{\bsn{\Lambda}}
\newcommand{\Gammau}{\bsn{\Gamma}}
\newcommand{\sigmauu}{\bsn{\sigma}}
\newcommand{\Omegau}{\bsn{\Omega}}
\newcommand{\phiu}{\bsn{\phi}}
\newcommand{\Phiu}{\bsn{\Phi}}
\newcommand{\Psiu}{\bsn{\Psi}}
\newcommand{\vepsuu}{\bsn{\varepsilon}}
\newcommand{\Upsu}{\bsn{\Upsilon}}
\newcommand{\Xiu}{\bsn{\Xi}}
\newcommand{\img}{\mathfrak{i}}

\newcommand{\UTu}{\tilde{\mathbf{U}}}
\newcommand{\www}{\tilde{\bfm{u}}}
\newcommand{\cC}{\mathcal{C}}
\newcommand{\cA}{\mathcal{A}}
\newcommand{\cZ}{\mathcal{Z}}

\newcommand{\gray}[1]{{\color{black!50!white} {#1}}}
\newcommand{\ale}[1]{{\color{blue!50!red} {#1}}}
\begin{document}

\title{Model Order Reduction based on Direct Normal Form: Application to Large Finite Element MEMS Structures Featuring Internal Resonance}

\titlerunning{Model Order Reduction based on Direct Normal Form}        

\author{Andrea Opreni$^{1}$         \and
        Alessandra Vizzaccaro$^{2}$ \and
        Attilio Frangi$^{1}$        \and
        Cyril Touz\'e$^{3}$		
}


\institute{$^{1}$ Department of Civil and Environmental Engineering\\
				  Politecnico di Milano\\
				  Piazza Leonardo da Vinci, 32, 20133 Milano MI\\
				  Tel.: +39.349.377.6864\\
				  \email{andrea.opreni@polimi.it}	\\
				  \\	
		   $^{2}$ Department of Engineering Mathematics\\
		   		  University of Bristol\\
		   		  Bristol BS8 1UB
				  \\
				  \\   
		   $^{3}$ Institute of Mechanical Sciences and Industrial Applications (IMSIA)\\ ENSTA Paris - CNRS - EDF - CEA - Institut Polytechnique de Paris\\
		          828 boulevard des mar\'echaux 91762 Palaiseau cedex
}


\maketitle

\begin{abstract}
Dimensionality reduction in mechanical vibratory systems poses challenges for distributed structures including geometric nonlinearities, mainly because of the lack of invariance of the linear subspaces. A reduction method based on direct normal form computation for large finite element (FE) models is here detailed. The main advantage resides in operating directly from the physical space, hence avoiding the computation of the complete eigenfunctions spectrum. Explicit solutions are given, thus enabling a fully non-intrusive version of the reduction method. The reduced dynamics is obtained from the normal form of the geometrically nonlinear mechanical problem, free of non-resonant monomials, and truncated to the selected master coordinates, thus making a direct link with the parametrisation of invariant manifolds. The method is fully expressed with a complex-valued formalism by detailing the homological equations in a systematic manner, and the link with real-valued expressions is established. A special emphasis is put on the treatment of second-order internal resonances and the specific case of a 1:2 resonance is made  explicit. Finally, applications to large-scale  models of Micro-Electro-Mechanical structures featuring 1:2 and 1:3 resonances are reported, along with considerations on computational efficiency.
\keywords{invariant manifold parametrisation \and normal form \and nonlinear normal modes \and harmonic balance \and model order reduction \and non-intrusive method}
\end{abstract}

\section{Introduction}
\label{sec:intro}
The  rate of development of computational power has enabled the simulation of complex systems. However, computational costs for the analysis of mechanical components are still not competitive enough on an industrial scale, where first guesses on the performance of mechanical components must be obtained rapidly. This limitation led to the development of techniques aimed at reducing the computational cost of numerical models while retaining sufficient reliability, {\em i.e.} Model Order Reduction (MOR) techniques \cite{Steindl01,mignolet13,BaurBenner14}. Restricting ourselves to the case of vibratory systems including large-amplitude displacements and thus exciting geometric nonlinearities, the first idea used since decades has been to project the nonlinear equations of motion onto a selected subset of the linear modes basis, see {\em e.g.}~\cite{Nayfeh79,Amabmodal} to cite only two examples. Unfortunately, the loss of invariance of linear eigenspaces, expressed through important nonlinear coupling terms that may happen with high-frequency modes~\cite{Vizza3d}, makes this approach not efficient in terms of accuracy~\cite{ShawPierre91,ShawPierre93,touze03-NNM,mignolet13,TouzeCISM}. The Proper Orthogonal Decomposition (POD) method offers a gain  since being able to modify slightly the orientation of the subspaces to better fit the curvatures of the nonlinear data, but it is still restricted to the use of linear orthogonal subspaces~\cite{Kerschen2005POD,KryslPOD,TOUZE:JFS:2007,SampaioPOD}.\\

On the other hand, nonlinear reduction methods start by defining a nonlinear relationship between the original coordinates and those of the reduced dynamics, hence providing a more accurate treatment of the nonlinear trajectories and thus faster convergence with fewer master modes. First attempts can be traced back to the work by Rosenberg~\cite{Rosenberg62} who introduced the term {\em Nonlinear Normal Mode} (NNM)~\cite{KerschenNNM09,PeetersNNM09}, which has been further developed in the works by Shaw and Pierre who first recognised that the concept of invariant manifold is key to compute accurate ROMs for nonlinear vibratory systems~\cite{ShawPierre91,ShawPierre93,PesheckMultiNNM,ShawCISM}. The normal form theory together with a truncation to a small subset of master modes have been proposed in~\cite{touzeLMA,touze03-NNM,TOUZE:JSV:2006,TouzeCISM}. On the mathematical point of view, invariant manifolds are investigated for a long time, see {\em e.g.}~\cite{Roberts89,LlaveIM97}, and the methods have been written in a unified and abstract general formalism thanks to the parametrisation methods~\cite{Cabre1,Cabre2,Cabre3}. The book by Haro {\em et al.}~\cite{Haro} gives a clear presentation of these  developments with a presentation oriented toward applications. In particular, it shows that the computations led by previous authors in vibration theory followed the guidelines of the two main parametrisation styles of invariant manifolds. While Shaw and Pierre used the {\em graph style}, the {\em normal form style} was used in~\cite{touzeLMA,touze03-NNM,TOUZE:JSV:2006}, with a real formulation in order to better fit the classical oscillatory framework.\\

A remarkable advancement within the correct and formal definition and properties of  nonlinear normal modes has been given by Haller {\em et al.} in a series of papers~\cite{Haller2016,PONSIOEN2018,VERASZTO,PONSIOEN2020}. In the conservative framework, existence and uniqueness of the searched invariant structures are given by the Lyapunov center theorem~\cite{Lyapunov1907}, stating that under non-resonance conditions, a two-dimensional manifold densely filled with periodic orbits exists for each couple of imaginary eigenvalues. As remarked in~\cite{Haller2016}, the picture is completely different for dissipative systems since the whole phase space is foliated by invariant manifolds tangent at the origin to linear subspace (see also \cite{NeildNF01,CIRILLO2016}). In order to give a correct, accurate and unequivocal definition of the searched reduction subspaces, Haller {\em et al.} introduced the notion of spectral submanifolds (SSM) as the smoothest nonlinear continuation of a spectral subspace of the linearised system. They proved that existence and uniqueness of SSM are related to non-resonance conditions and a  spectral quotient computed from the real part of the spectrum of the linearised system. In the conservative case, the manifolds are named as Lyapunov Subcenter Manifold (LSM)~\cite{Kelley69,VERASZTO}.
Retrospectively, the earlier computations of NNMs turned out to be computations of LSM when damping was not taken into account. On the computational point of view, Haller {\em et al.} used normal form style~\cite{Haro,PONSIOEN2018}, and proposed an automated computational framework allowing to go to any order of asymptotic development while taking into account damping and various nonlinear terms.\\

Applications of methods inherited from dynamical systems theory to finite element (FE) problems has remained scarce until recently. A major motivation was that all the methods described in the last paragraph have as common starting point the mechanical system written in modal space. Unfortunately, FE methods, where meshes with millions of degrees-of-freedom (dofs) are routinely used, make this step out of reach. The direct computation of NNMs from large Finite Element Models  
\cite{KerschenNNM09,PeetersNNM09,RENSON2016} has been addressed with harmonic balance approaches or shooting procedures. Thanks to their purely computational nature, these techniques are very general and can be applied to a wide range of problems, but they require a certain amount of computational effort to be applied to realistic structures and fail at serving the aim of generating agile reduced order models (ROMs).
On the other hand, for FE structures including geometric nonlinearities, a number of methods have been proposed in the last 20 years in order to formulate ROMs that can be easily computed, in the best case in a non-intrusive manner, {\em i.e.} without the need to enter new calculations at the elementary level in the code, such that any research-oriented or commercial FE code could be used as starting point. The stiffness evaluation procedure (STEP) uses the linear modes as projection basis~\cite{muravyov,Perez2014}, but can be simply used as a technique to compute non-intrusively linear and nonlinear characteristics. The implicit condensation and expansion (ICE) method relies on a set of applied forces to keep trace of non-resonant couplings among modes~\cite{Hollkamp2008,kuether2015,FRANGI2019,NicolaidouIceKE}. Finally modal derivatives (MD) have been introduced with the aim of taking the amplitude dependence of modes into account, thus pursuing the same goal as the NNMs~\cite{IDELSOHN1985,Weeger2016,Jain2017,Rutzmoser}. 
ICE and MD assume the manifold to be velocity independent, assumption which is the more fulfilled, the larger the slow/fast separation between the slave and master coordinates (\cite{HallerSF,Vizzaccaro:NNMvsMD,YichangICE}).
In~\cite{Vizzaccaro:NNMvsMD}, it is estimated that a ratio between eigenfrequencies of the slave modes and those of the master modes of 
at least 3 ensures the correctness of MDs in the prediction of the hardening/softening behaviour (\cite{YichangVib,YichangNODYCON,YichangICE}). 
The inclusion of velocity dependence in the MD approach have been proposed in~\cite{RixenMDNF19} to overcome this limitation, and leads to similar formulations than those reported in~\cite{artDNF2020}.
\\

Consequently there is still a real need of direct applications of invariant-based reduction methods, using the general theorems from dynamical systems theory, and that could be applied non-intrusively on the dofs of FEM problems in physical space. Pursuing previous developments using normal form theory~\cite{touze03-NNM,TOUZE:JSV:2006}, explicit and non-intrusive real formulas have been proposed in~\cite{artDNF2020}. All the coefficients of second and third-order nonlinear mappings have been derived by rewriting the normal form approach proposed in~\cite{touze03-NNM,TOUZE:JSV:2006} directly in the physical space. A real formalism was used so that all along the calculation, the results are expressed in the form of oscillator-like equations. 
Damping and forcing were taken into account and illustrative examples on blades and beams were reported in~\cite{artDNF2020,YichangVib}. At the same time, the computation based on SSM led by the group of Haller also proposed a direct computation that has been written in the open code SSM tool 2.0~\cite{SSMtool20}.\\


The present article aims at completely rewriting the direct normal form approach proposed in~\cite{artDNF2020} in order to explain more clearly its settings and possible further developments. Secondly, the special case of second-order internal resonance, that had not been treated in~\cite{touze03-NNM,TOUZE:JSV:2006,artDNF2020}, is here detailed.
Finally, applications to large-scale Micro-Electro-Mechanical systems (MEMS) are derived in order to show the potentiality of the method to deal with millions of dofs yet providing fast and accurate ROMs. MEMS structures are generally actuated at resonance, they are subjected to geometric nonlinearities due to large transformations and have very small damping values as they operate in near-vacuum packages, hence  showing highly nonlinear dynamical features that are rarely observed at the macro scale~\cite{czaplewski2018a,czaplewski2018b,czaplewski2020,frangi2020,ruzziconi2021a,ruzziconi2021b}. Furthermore, nonlinear dynamic properties of MEMS can be tailored to yield performance that would not be accessible through operation in the linear regime \cite{alcheikh2021,hajjaj2020,shoshani2020,li2017}. A remarkable example of successful application of nonlinear mode interaction for the development of highly efficient mechanical filters is for instance reported in  \cite{hafiz2020}. This implies that the proposed method would play a major role in this field since it would ensure a fast and efficient estimation of the frequency response functions of structures within times that are compatible with industrial design requirements.\\

The paper is organised as follows. Section~\ref{sec:theory} details the reduction method based on direct normal form (DNF). A first-order (state-space) formulation along with a complex-valued formalism is introduced in Sections~\ref{sec:param} and~\ref{subsec:dnf0}, with the general idea of rewriting the DNF in a more complete and symmetric formalism, a necessary step for further developments to include more easily new forces  and/or going to higher-orders. Section~\ref{subsec:homolog} details the second and third-order homological equations, written from the physical space, such that the method can be easily understood and adopted by a wider group of researchers. These steps are important novelties as compared to the calculations presented in~\cite{artDNF2020} where only the translations of the generic method presented in~\cite{touze03-NNM,TOUZE:JSV:2006} for direct application in physical space was given, in a real formalism. Once the complete normal form is computed, the reduction method is explained in Section~\ref{sec:ROM}. More insights in the computations are given for the simple case of single master coordinate in Section~\ref{subsec:rom1master}. The case of second-order internal resonance is then tackled in Section~\ref{subsec:IRmulti}. 
Finally, Section~\ref{sec:real} explains how one can pass from the complex formalism detailed here to the real  nonlinear mappings provided in~\cite{artDNF2020}, hence bridging earlier works.

Section~\ref{sec:results} shows the numerical results obtained for three different MEMS structures. First, 
a MEMS micromirror undergoing large rotations is considered. In this case a single master mode is needed to correctly retrieve the hardening behaviour up to very large angles of rotations. This example is important since other reduction methods had been tested before, all of them giving incorrect results, as illustrated with application of the ICE method. Then a MEMS beam resonator that features a 1:3 internal resonance \cite{czaplewski2019}, and a MEMS arch resonator showing 1:2 internal resonance are selected in order to show the ability of the method to deal with systems that feature internally resonant modes. All numerical simulations in the present work are compared with full-order harmonic balance finite element simulations of the systems \cite{opreni2021,DetrouxHBStability,Blahos2020}.\\

\section{Equations of motion and direct normal form approach}\label{sec:theory}

The aim of this Section is to derive a general formalism to write direct expressions for the computation of the normal transform on the physical dofs. In particular, homological equations adapted to the framework of mechanical systems including geometric nonlinearities are given.

\subsection{Equations of motion and linear decomposition}
\label{sec:param}

Let us consider the equation of motion associated with an undamped mechanical system subjected to geometric nonlinearities:
\begin{equation}
\label{eq:damped_dyn}
\Mu\ddot{\Uu}+\Ku\Uu +\Gu (\Uu ,\Uu)+\Hu (\Uu ,\Uu ,\Uu) = 0,
\end{equation}
where $\Mu$ denotes the mass matrix, $\Uu$ the vector of nodal displacements of dimension $\Nr$, $\dot{(\cdot)}$ the time derivative operator, $\Ku$ the stiffness matrix, $\Gu (\Uu ,\Uu)$ and 
$\Hu (\Uu ,\Uu ,\Uu)$ the quadratic and cubic force terms, respectively. This starting point is common to any finite element discretisation of the linear momentum equation for mechanical systems. In particular in the framework of three-dimensional linear elasticity with large transformations, Eq.~\eqref{eq:damped_dyn} is exact. For the sake of completeness, Appendix~\ref{sec:notation} gives some details on the adopted notations, and Appendix~\ref{sec:nl_expr} recalls how the quadratic and cubic terms are obtained from the decomposition of the internal force vector.

The solution of a generalised eigenvalue problem formulated on the linear part of Eq.~\eqref{eq:damped_dyn} yields a finite set of real-valued mass-normalised eigenmodes $\phiu_s$ which are collected column-wise in the eigenvector matrix $\Phiu$. The following equalities hold:
\begin{subequations}
\begin{align}
	\Phiu^{\mathrm{T}}\Mu\Phiu =&\, \Iu, \\
	\Phiu^{\mathrm{T}}\Ku\Phiu =&\, \Omegau^{2},
\end{align}
\end{subequations}
where $\Omegau$ is a diagonal matrix that stores the eigenfrequencies $\omega_s$. Equation~\eqref{eq:damped_dyn} can be mapped to modal coordinates $\uu$ through the linear transformation $\Uu=\Phiu\uu$, leading to:
\begin{equation}\label{eq:modal_dyn}
	\ddot{\uu} + \Omegau^{2}\uu + \gu(\uu,\uu) + \hu(\uu,\uu,\uu) = 0,
\end{equation}
with:
\begin{subequations}
	\begin{align}
		\gu(\uu,\uu) =&\, \Phiu^{\mathrm{T}}\Gu (\Uu ,\Uu), \\ \hu(\uu,\uu,\uu) =&\, \Phiu^{\mathrm{T}}\Hu (\Uu ,\Uu,\Uu).
	\end{align}
\end{subequations}
Modal decomposition of Eq.~\eqref{eq:modal_dyn} is computationally infeasible for large systems and must be avoided. In order to derive a general formalism for computing the normal transform of Eq.~\eqref{eq:damped_dyn},
let us first express Eq.~\eqref{eq:damped_dyn} in state-space formalism, {\em i.e.} as a first-order dynamical system by introducing the velocity $\Vu$ as the time derivative of $\Uu$:
\begin{align}
\label{eq:damped_dyn_first}
    \left[\begin{array}{cc}
		\Mu & \zerou \\
		\zerou & -\Iu
	\end{array}\right] 
	\left\{\begin{array}{l}
		\dot{\Vu} \\
		\dot{\Uu}
	\end{array}\right\}
	+
	\left[\begin{array}{cc}
 		\zerou &  \Ku \\
		\Iu &  \zerou
	\end{array}\right] 
	\left\{\begin{array}{l}
		\Vu \\
		\Uu
	\end{array}\right\}
	+ 
	\left\{\begin{array}{c}
	\Gu (\Uu ,\Uu) + \Hu (\Uu ,\Uu ,\Uu) \\
	0
	\end{array}\right\} = \zerou.
\end{align}
Note that the previous developments on the DNF approach proposed in~\cite{artDNF2020} aimed at avoiding this first-order formulation in order to keep real expressions and oscillator-like equations, with second-order derivatives in time, throughout the calculations. This could be realised at the cost of an important loss in the symmetries of the problem. Moreover, generalisations of the second-order formalism to other linear and nonlinear forces, {\em e.g.} linear viscous damping, electrostatic and piezoelectric forces, become more difficult to handle compared to a first-order approach. 
Consequently, we propose in this article a complete rewriting of the DNF approach, using state-space formulation as starting point, but keeping in all calculations the link to the usual mechanical representation as given in Eq.~\eqref{eq:damped_dyn}, so that further generalisations will be easier to derive.

For the sake of conciseness, the diagonalisation of the linear part of Eq.~\eqref{eq:damped_dyn_first} is reported in Appendix~\ref{sec:diag} and leads to:
\begin{equation}\label{eq:gen_coor}
	\dot{\pu} = \Lambdau\pu + \Gammau(\pu,\pu) + \Deltau(\pu,\pu,\pu),
\end{equation}
with $\pu$ the generalised coordinates, $\Lambdau$ a diagonal matrix, and $\Gammau(\pu,\pu)$, $\Deltau(\pu,\pu,\pu)$ nonlinear operators in generalised coordinates. In index notation, it reads:
\begin{equation}\label{eq:gen_coor_index}
	\dot{\pr}_s =  \lambda_{s}\pr_s + \sum_{k,l=1}^{2\Nr} \Gamma_{skl}\pr_k\pr_l + \sum_{k,l,m=1}^{2\Nr} \Delta_{sklm} \pr_k\pr_l\pr_m ,\quad \forall\; s=1,...,2\Nr,
\end{equation}
with $\lambda_s$ the diagonal entries of $\Lambdau$, and $\Gamma_{skl}$, $\Delta_{sklm}$ the scalar coupling coefficients. Note that due to the state-space formulation,  $k,l,m$ are summed from 1 to $2\Nr$ in Eq.~\eqref{eq:gen_coor_index}, and the system is now of size $2\Nr$. 
Further details on the derivation of Eq.~\eqref{eq:gen_coor} are reported in Appendix~\ref{sec:diag}.
In particular, the structure of the diagonal $\Lambdau$ is such that:
\begin{equation}
\label{eq:lambdamain}
	\Lambdau = 
	\left[\begin{array}{cc}
		\img\Omegau & \zerou \\
		\zerou & -\img\Omegau
	\end{array}\right],
\end{equation}
and the $\lambda_s$ coefficients read:
\begin{subequations}\label{eq:fo_coeff}
	\begin{align}
		\forall\, s =&\, 1,...,\Nr,\nonumber \\
		\lambda_{s} = &\, +\img\omega_{s},  \\
    	\lambda_{s+N} =&\, -\img\omega_{s},
	\end{align}
\end{subequations}	
with $\img$ the imaginary unit. Moreover displacements $\Uu$, velocities $\Vu$, modal displacements $\uu = [\ur_1, ..., \ur_N]^{\mathrm{T}}$, and modal velocities $\vu= [\vr_1, ..., \vr_N]^{\mathrm{T}}$ are related to the generalised coordinates via the following relationships:
\begin{subequations}
	\begin{align}
		\Uu = &\, \sum_{s=1}^{\Nr} \phiu_s \ur_s\\
	    \Vu = &\, \sum_{s=1}^{\Nr} \phiu_{s} \vr_s \\
	    \ur_s = &\, \pr_s+\pr_{s+\Nr}, \\
	    \vr_s = &\, \lambda_{s}\pr_s+\lambda_{s+\Nr}\pr_{s+\Nr}.
	\end{align}
\end{subequations}	

As highlighted by Eq.~\eqref{eq:gen_coor}, the diagonalisation of the linear part of Eq.~\eqref{eq:damped_dyn_first} through a linear change of coordinates does not guarantee a decoupling of the nonlinear terms. 
This implies that a given couple $\pr_s$ and $\pr_{s+\Nr}$ does not define an invariant subspace of the  system. The next steps of the developments aim at deriving a nonlinear mapping that could express the dynamics in an invariant-based span of the phase space, such that once the nonlinear mapping is computed with correct truncation, reduced-order models could be easily obtained by keeping only a few master coordinates. This will be realised thanks to the normal form computation, directly applied to Eq.~\eqref{eq:damped_dyn_first}.

\subsection{Direct Normal Form Setting}\label{subsec:dnf0}
The normal form approach has been first introduced by Poincar{\'e}, leading to well-known theorems~\cite{Poincare,Dulac1912}. In the context of vibratory systems, it has been introduced for model order reduction purpose in~\cite{touze03-NNM,TOUZE:JSV:2006,TouzeCISM}, with the additional idea of truncating to a few subset of master coordinates once the full nonlinear mapping is computed, hence retrieving the parametrisation method of invariant manifold with normal form style \cite{Haro,Haller2016}, as will be further discussed in Sec. \ref{sec:ROM}.
In this contribution, we follow the guidelines of the normal transform by first deriving the complete mapping for all coordinates, and then truncating to obtain a ROM. Consequently a nonlinear relationship between the original $(\Uu ,\Vu)$ variables (displacement and velocity $\Nr$-dimensional vectors) and the {\em normal} $\zu$ coordinate, a $2\Nr$-dimensional vector describing the dynamics in an invariant-based span, is introduced as:
\begin{subequations}\label{eq:mapping_nl}
	\begin{align}
		\Uu = \Psiu(\zu), \\
		\Vu = \Upsu(\zu),
	\end{align}
\end{subequations}
with $\Psiu(\zu)$, $\Upsu(\zu)$ nonlinear polynomial mappings. 
Mappings are expanded in terms of their components and written in either tensorial form:
\begin{subequations}
	\begin{align}
		\Psiu(\zu) = &\, \Psiu^{(1)}\zu + \Psiu^{(2)}(\zu,\zu) + \Psiu^{(3)}(\zu,\zu,\zu) + O(\|\zu\|^{4}), \\
		\Upsu(\zu) = &\, \Upsu^{(1)}\zu + \Upsu^{(2)}(\zu,\zu) + \Upsu^{(3)}(\zu,\zu,\zu) + O(\|\zu\|^{4}),
	\end{align}
\end{subequations}
where the in-parenthesis upperscript refers to the order of the polynomial, or in indicial form:
\begin{subequations}
	\begin{align}
		\Psiu(\zu) = &\, \sum_{s=1}^{2\Nr} \Psiu^{(1)}_s\zr_s + \sum_{k,l=1}^{2\Nr} \Psiu^{(2)}_{kl}\zr_k\zr_l + \sum_{k,l,m=1}^{2\Nr}\Psiu^{(3)}_{klm}\zr_k\zr_l\zr_m  + O(\|\zu\|^{4}), \\
		\Upsu(\zu) = &\, \sum_{s=1}^{2\Nr} \Upsu^{(1)}_s\zr_s + \sum_{k,l=1}^{2\Nr} \Upsu^{(2)}_{kl}\zr_k\zr_l + \sum_{k,l,m=1}^{2\Nr}\Upsu^{(3)}_{klm}\zr_k\zr_l\zr_m + O(\|\zu\|^{4}),
	\end{align}
\end{subequations}
where a bold capital letter with one or more indices denotes a vector. Since modal coordinates define invariant subspaces in the absence of nonlinearities, the first-order maps should correspond to the usual linear decomposition. This is also in line with the general idea of finding a continuation of the linear mode subspace where the higher-order terms will account for the amplitude-dependence of modal quantities, given by the curvatures of the invariant manifolds~\cite{Boivin95,TouzeCISM,BuzaHaller}. Note that this is a common idea with the quadratic manifold approach including MD as proposed in~\cite{Jain2017,Rutzmoser}, the only difference being that in the present derivation, velocities are fully taken into account.
The linear terms of the mapping are thus simply expressed as:
\begin{subequations}
	\begin{align}
		\label{eq:fo_map}
		& \forall\,s= 1,...,\Nr, \nonumber \\
		& \Psiu_{s}^{(1)} = \phiu_s, \\
		& \Psiu_{s+\Nr}^{(1)} = \phiu_s, \\
		& \Upsu_{s}^{(1)} = \lambda_{s}\phiu_s, \\
		& \Upsu_{s+\Nr}^{(1)} = \lambda_{s+\Nr}\phiu_s.
	\end{align}
\end{subequations}
This choice makes first-order terms identical to the linear transformation that maps $\Uu$ and $\Vu$ to the generalised coordinates $\pu$ (see Appendix~\ref{sec:diag}). In order to derive complete expressions for the homological equations~\cite{Murdock,Jezequel91,LamarqueUP}, one needs to first express the  dynamics of the normal variable $\zu$, {\em i.e.}\ the normal form of the initial problem, as:
\begin{equation}\label{eq:red_dyn}
	\dot{\zu} = \fuu(\zu),
\end{equation}
where $\fuu(\zu)$ is expressed as a polynomial function of $\zu$ with unknown coefficients:
\begin{equation}
	\fuu(\zu) = \fuu^{(1)}\zu + \fuu^{(2)}(\zu,\zu) + \fuu^{(3)}(\zu,\zu,\zu) + O(\|\zu\|^{4}),
\end{equation}
or in index form:
\begin{equation}
	f_s(\zu) =  f^{(1)}_{s}\zr_s + \sum_{k,l=1}^{2\Nr}  f^{(2)}_{skl}\zr_k\zr_l + \sum_{k,l,m=1}^{2\Nr} f^{(3)}_{sklm}\zr_k\zr_l\zr_m + O(\|\zu\|^{4}), \quad \forall\;s=\,1,...,2\Nr.
\end{equation}
This first-order form for the dynamics of the normal coordinates is in full accordance with~\cite{Haro,Haller2016,PONSIOEN2018}. Since $\Psiu^{(1)}_s=\phiu_s$, then $\fuu^{(1)}$ is taken equal to $\Lambdau$,
or equivalently $f^{(1)}_{s}=\lambda_s$. This is the logical consequence of using an identity-tangent change of coordinates that leaves the linear diagonalised part of the dynamics unchanged.
The homological equations are obtained by discarding the dependence on time and equating to zero the collections of terms with same powers in the asymptotic developments, order by order. To that purpose, the time derivatives of the mappings are obtained from Eqs.~\eqref{eq:mapping_nl} and \eqref{eq:red_dyn} as:
\begin{subequations}\label{eq:xdot}
	\begin{align}
		\dot{\Uu} =&\, \Psiu^{(1)}\fuu^{(1)}\zu +  \nonumber \\
			   &\,
			   \Psiu^{(1)}\fuu^{(2)}(\zu,\zu)  + \Psiu^{(2)}(\fuu^{(1)}\zu,\zu) + \Psiu^{(2)}(\zu,\fuu^{(1)}\zu) + \nonumber \\
			   &\, \Psiu^{(1)}\fuu^{(3)}(\zu,\zu,\zu)  + \Psiu^{(2)}(\fuu^{(2)}(\zu,\zu),\zu) + \Psiu^{(2)}(\zu,\fuu^{(2)}(\zu,\zu)) +	\nonumber \\  
			   &\,
			   \Psiu^{(3)}(\fuu^{(1)}\zu,\zu,\zu) + \Psiu^{(3)}(\zu,\fuu^{(1)}\zu,\zu) + \Psiu^{(3)}(\zu,\zu,\fuu^{(1)}\zu) +	O(\|\zu \|^{4}),     \\
		\dot{\Vu} =&\, \Upsu^{(1)}\fuu^{(1)}\zu +  \nonumber \\
			   &\,
			   \Upsu^{(1)}\fuu^{(2)}(\zu,\zu)  + \Upsu^{(2)}(\fuu^{(1)}\zu,\zu) + \Upsu^{(2)}(\zu,\fuu^{(1)}\zu) + \nonumber \\
			   &\, \Upsu^{(1)}\fuu^{(3)}(\zu,\zu,\zu)  + \Upsu^{(2)}(\fuu^{(2)}(\zu,\zu),\zu) + \Upsu^{(2)}(\zu,\fuu^{(2)}(\zu,\zu)) +	\nonumber \\  
			   &\,
			   \Upsu^{(3)}(\fuu^{(1)}\zu,\zu,\zu) + \Upsu^{(3)}(\zu,\fuu^{(1)}\zu,\zu) + \Upsu^{(3)}(\zu,\zu,\fuu^{(1)}\zu) +	O(\|\zu \|^{4}),
	\end{align}
\end{subequations}
where some properties defined in Appendix~\ref{sec:notation} have been used. Alternatively, using indicial notation:
\begin{subequations}
	\begin{align}
		\dot{\Uu} = \,& \sum_{k=1}^{2\Nr}\lambda_{k}\Psiu^{(1)}_{k}\zr_k + \nonumber \\
					&\, \sum_{k,l=1}^{2\Nr}\left[ \sum_{s=1}^{2\Nr}f^{(2)}_{skl}\Psiu^{(1)}_{s}  + (\lambda_{k}+\lambda_{l})\Psiu^{(2)}_{kl} \right] \zr_k\zr_l + \nonumber \\
					&\, \sum_{k,l,m=1}^{2\Nr}\left[ \sum_{s=1}^{2\Nr}\left(f^{(3)}_{sklm}\Psiu^{(1)}_{s} + f^{(2)}_{skl}\Psiu^{(2)}_{sm} + f^{(2)}_{slm}\Psiu^{(2)}_{ks}\right) +  (\lambda_{k}+\lambda_{l}+\lambda_{m})\Psiu^{(3)}_{klm} \right] \zr_k\zr_l\zr_m + O(\|\zu \|^{4}), \\
		\dot{\Vu} = \,& \sum_{k=1}^{2\Nr}\lambda_{k}\Upsu^{(1)}_{k}\zr_k + \nonumber \\
					&\, \sum_{k,l=1}^{2\Nr}\left[ \sum_{s=1}^{2\Nr}f^{(2)}_{skl}\Upsu^{(1)}_{s}  + (\lambda_{k}+\lambda_{l})\Upsu^{(2)}_{kl} \right] \zr_k\zr_l + \nonumber \\
					&\, \sum_{k,l,m=1}^{2\Nr}\left[ \sum_{s=1}^{2\Nr}\left(f^{(3)}_{sklm}\Upsu^{(1)}_{s} + f^{(2)}_{skl}\Upsu^{(2)}_{sm} + f^{(2)}_{slm}\Upsu^{(2)}_{ks}\right) +  (\lambda_{k}+\lambda_{l}+\lambda_{m})\Upsu^{(3)}_{klm} \right] \zr_k\zr_l\zr_m + O(\|\zu \|^{4}).
	\end{align} 
\end{subequations}
The same operation can be applied for system nonlinearities in order to correctly distinguish each order appearing in the expansions:
\begin{subequations}
	\begin{align}\label{eq:nl_terms}
	\Gu(\Uu,\Uu) =&\, 
				  \Gu(\Psiu^{(1)}\zu,\Psiu^{(1)}\zu) + \nonumber \\
				  &\, 
				  \Gu(\Psiu^{(2)}(\zu,\zu),\Psiu^{(1)}\zu) + \nonumber \\
				  &\,
				  \Gu(\Psiu^{(1)}\zu,\Psiu^{(2)}(\zu,\zu)) + \nonumber \\
				  &\,
				  \Gu(\Psiu^{(3)}(\zu,\zu,\zu),\Psiu^{(1)}\zu) + \nonumber \\
				  &\,  
				  \Gu(\Psiu^{(1)}\zu,\Psiu^{(3)}(\zu,\zu,\zu)) + \nonumber \\
				  &\, 
				  \Gu(\Psiu^{(2)}(\zu,\zu),\Psiu^{(2)}(\zu,\zu))  +  O(\|\zu \|^{5}), \\
	\Hu(\Uu,\Uu,\Uu) = &\, 
					\Hu(\Psiu^{(1)}\zu,\Psiu^{(1)}\zu,\Psiu^{(1)}\zu) + \nonumber \\
					&\, 
					\Hu(\Psiu^{(2)}(\zu,\zu),\Psiu^{(1)}\zu,\Psiu^{(1)}\zu) + \nonumber \\
					&\,
					 \Hu(\Psiu^{(1)}\zu,\Psiu^{(2)}(\zu,\zu),\Psiu^{(1)}\zu)  + \nonumber \\
					&\,
					\Hu(\Psiu^{(1)}\zu,\Psiu^{(1)}\zu,\Psiu^{(2)}(\zu,\zu))  +  O(\|\zu \|^{5}) .
	\end{align}
\end{subequations}
As before the same equation in full indicial notation is also provided for clarity:
\begin{subequations}
	\begin{align}\label{eq:nl_terms_index}
	\Gu(\Uu,\Uu) =&\, \sum_{k,l=1}^{2\Nr}\Gu(\Psiu^{(1)}_k,\Psiu^{(1)}_l)\zr_k\zr_l + \nonumber \\
				  &\, \sum_{k,l,m=1}^{2\Nr}\left[\Gu(\Psiu^{(2)}_{kl},\Psiu^{(1)}_m) + \Gu(\Psiu^{(1)}_k,\Psiu^{(2)}_{lm})\right]\zr_k\zr_l\zr_m + \nonumber \\
				  &\, \sum_{k,l,m,n=1}^{2\Nr}\left[ \Gu(\Psiu^{(3)}_{klm},\Psiu^{(1)}_n) +  \Gu(\Psiu^{(1)}_{k},\Psiu^{(3)}_{lmn}) + \Gu(\Psiu^{(2)}_{kl},\Psiu^{(2)}_{mn}) \right] \zr_k\zr_l\zr_m\zr_n   + O(\|\zu \|^{5}) , \\
	\Hu(\Uu,\Uu,\Uu) =&\, \sum_{k,l,m=1}^{2\Nr}\Hu(\Psiu^{(1)}_k,\Psiu^{(1)}_l,\Psiu^{(1)}_m)\zr_k\zr_l\zr_m + \nonumber \\	 
				  &\,   \sum_{k,l,m,n=1}^{2\Nr}\left[ \Hu(\Psiu^{(2)}_{kl},\Psiu^{(1)}_{m},\Psiu^{(1)}_{n}) +  \Hu(\Psiu^{(1)}_{k},\Psiu^{(2)}_{lm},\Psiu^{(1)}_{n}) +   \Hu(\Psiu^{(1)}_{k},\Psiu^{(1)}_{l},\Psiu^{(2)}_{mn})\right]  \zr_k\zr_l\zr_m\zr_n + O(\|\zu \|^{5}) ,
	\end{align}
\end{subequations}

where expansions of nonlinear terms up to fourth-order is done to better highlight the nested structure of the derived equations. 
It is worth stressing that the coefficients of the nonlinear mapping $\Psiu(\zu)$, $\Upsu(\zu)$ 
are unknown at this stage,
together with the expression of the normal form $\fuu(\zu)$. 
The aim is to find the simplest possible vector field $\fuu$ upon nonlinear transforms and this is achieved step-by-step by writing for each order the homological equations that collect the same powers of the normal variable $\zu$.

\subsection{Homological Equations}
\label{subsec:homolog}
Computation of $\Psiu(\zu)$, $\Upsu(\zu)$, and $\fuu(\zu)$ is performed by substituting Eqs.~\eqref{eq:mapping_nl} and \eqref{eq:xdot} in Eq.~\eqref{eq:damped_dyn_first} and by collecting terms of equal power in $\zu$ in hierarchical order, hence providing for each order the corresponding homological equation~\cite{Murdock,Haro}. The first-order homological equation reads:
\begin{equation}\label{eq:fo_homological}
	\left[\begin{array}{cc}
		\Mu & \zerou \\
		\zerou & -\Iu
	\end{array}\right] 
	\left\{\begin{array}{l}
		\Upsu^{(1)}\fuu^{(1)}\zu \\
		\Psiu^{(1)}\fuu^{(1)}\zu
	\end{array}\right\} + 
	\left[\begin{array}{cc}
		\zerou & \Ku \\
		\Iu & \zerou
	\end{array}\right]
	\left\{\begin{array}{l}
		\Upsu^{(1)}\zu \\
		\Psiu^{(1)}\zu
	\end{array}\right\} = 
	\zerou.
\end{equation}
The $2\Nr$ equations have $4\Nr$ unknowns, that are the $2\Nr$ coefficients of $\fuu^{(1)}$ and  $\Nr$ vector pairs $(\Upsu^{(1)},\Psiu^{(1)})$. Therefore, the solution is not unique and any triplet $(\Upsu^{(1)},\Psiu^{(1)},\fuu^{(1)}$) that satisfies Eq.~\eqref{eq:fo_homological} is a valid  choice. The second part of Eqs.~\eqref{eq:fo_homological} allows writing
\begin{equation}
\Upsu^{(1)}_{s} = \Psiu^{(1)}_{s} \fuu^{(1)}_{s},\quad\forall \; s=1, .., 2\Nr,
\end{equation}
showing that a direct relationship exists between the two parts of the mapping (displacement/velocity). Using the identity-tangent solution as justified before leads to the choice $\fuu^{(1)}=\Lambdau$ such that in indicial notation one can simply write $\Upsu^{(1)}_s = \lambda_{s}\Psiu^{(1)}_s$,   $\forall\, s = 1,...,2\Nr $. The first line of Eq.~\eqref{eq:fo_homological} then becomes:
\begin{equation}
		\left[ \lambda_{s}^{2}\Mu + \Ku \right] \Psiu^{(1)}_s =  \zerou ,\quad
		\forall \; s=1, .., 2\Nr,
\end{equation}
which makes appear the usual linear eigenproblem in vibration theory, thus justifying again the   adopted choices $\Psiu^{(1)} = \Phiu$ and $\fuu^{(1)}=\Lambdau$.\\

Collecting second-order terms yields the second-order homological equation written in a manner that explicitly encompasses the usual mechanical context as:
\begin{align}\label{eq:so_homological}
	&
	\left[\begin{array}{cc}
		\Mu & \zerou \\
		\zerou & -\Iu
	\end{array}\right] 
	\left\{\begin{array}{l}
		\Upsu^{(1)}\fuu^{(2)}(\zu,\zu)\\
		\Psiu^{(1)}\fuu^{(2)}(\zu,\zu)
	\end{array}\right\} +  \nonumber \\
	&
	\left[\begin{array}{cc}
		\Mu & \zerou \\
		\zerou & -\Iu
	\end{array}\right] 
	\left\{\begin{array}{l}
		\Upsu^{(2)}(\fuu^{(1)}\zu,\zu) +  \Upsu^{(2)}(\zu,\fuu^{(1)}\zu) \\
		\Psiu^{(2)} (\fuu^{(1)}\zu,\zu) + \Psiu^{(2)} (\zu,\fuu^{(1)}\zu)
	\end{array}\right\} +	
	\nonumber \\
	&
	\left[\begin{array}{cc}
		\zerou & \Ku \\
		\Iu & \zerou
	\end{array}\right]
	\left\{\begin{array}{l}
		\Upsu^{(2)} (\zu,\zu)\\
		\Psiu^{(2)} (\zu,\zu)
	\end{array}\right\} = 
	\left\{\begin{array}{c}
		-\Gu (\Psiu^{(1)}\zu,\Psiu^{(1)}\zu)\\
		\mathbf{0}
	\end{array}\right\}.
\end{align}
Three remarks on this second-order homological equation can be made. First, one can observe that 
$\Psiu^{(2)}$ and $\Upsu^{(2)}$ are expressed in terms of the lower order terms, a feature that is typical of asymptotic developments. Second, the quadratic part of the nonlinear internal force vector $\Gu$ now directly appears at the right-hand side of Eq.~\eqref{eq:so_homological}, underlining that the formalism can be extended in order to include different nonlinear forces, with more complex features like velocity-dependence or even more involved physics with new variables, {\em e.g.} electrostatic or piezo-electric couplings. Even though this is not addressed in the present contribution, we believe that the framework developed in this work is sufficiently general so that inclusions of new forces can be now simply incorporated by modifying the right-hand side. 
The last comment regards the dependence of the velocity mapping  $\Upsu$ on the displacement mapping  $\Psiu$. As in Eq.~\eqref{eq:fo_homological}, the second part of \eqref{eq:so_homological} can be extracted to show that a linear relation between $\Upsu^{(2)}$ and $\Psiu^{(n)}$, $\fuu^{(n)}$ with $n\leq 2$, exists. This implies that in the present context of conservative vibratory systems, the velocity mapping $\Upsu$ is not independent of the displacement mapping $\Psiu$, thus allowing us to keep only $\Psiu$ at all orders, and then derive $\Upsu$  thanks to the known relationships, only in case this mapping is needed. 
At second-order, the relationship between velocity and displacement mappings read:\\
\begin{equation}\label{eq:so_homologicalUpsilon}
\Upsu^{(2)}_{kl} = (\lambda_{k}+\lambda_{l})\Psiu^{(2)}_{kl} + \sum_{s=1}^{2\Nr}f_{skl}^{(2)}\Psiu^{(1)}_s , \qquad \forall\, k,l = 1,...,2\Nr.
\end{equation}
On the other hand the first lines of Eq.~\eqref{eq:so_homological} can be written in terms of the displacement mapping $\Psiu$ only, as:
\begin{align}\label{eq:so_homologicalPsi}
    \sum_{s=1}^{2\Nr}\left[(\lambda_{s}+\lambda_{k}+\lambda_{l})\Mu\right] f^{(2)}_{skl}\Psiu^{(1)}_s  +  \left[(\lambda_{k}+\lambda_{l})^2\Mu  + \Ku \right] \Psiu^{(2)}_{kl} = 
    -\Gu(\Psiu^{(1)}_k,\Psiu^{(1)}_l), \quad \forall\, k,l = 1,...,2\Nr.
\end{align}
%

As for Eq.~\eqref{eq:fo_homological}, the system is underdetermined, {\em i.e.} the number of unknowns 
$(\Psiu^{(2)},\Upsu^{(2)},\fuu^{(2)})$ is larger that the number of equations. 
The aim of the normal form  is to simplify as much as possible the resulting dynamical system, {\em i.e.} to choose the solution such that $\fuu$ has the smallest number of terms, and is, in the ideal case, simply linear. Retaining this choice to solve Eq.~\eqref{eq:so_homologicalPsi} leads to selecting $f^{(2)}_{skl}=0,\,\forall\, s,k,l=1,...,2\Nr$ whenever it is possible, and thus to express $\Psiu^{(2)}_{kl}$ as:
\begin{align}
\label{eq:so_map}
	\Psiu^{(2)}_{kl} =  -\left[(\lambda_{k}+\lambda_{l})^{2}\Mu + \Ku \right]^{-1}\Gu(\Psiu^{(1)}_{k},\Psiu^{(1)}_{l}), \quad \forall\, k,l = 1,...,2\Nr.
\end{align}
When Eq.~\eqref{eq:so_map} is solvable, the normal form of the system will be free of quadratic coupling terms.
The condition for which $\left[(\lambda_{k}+\lambda_{l})^2\Mu  + \Ku \right]$ is singular 
is related to the usual second-order resonances. 
As already remarked in previous developments with normal transforms~\cite{touze03-NNM}, since the spectrum is composed of purely imaginary conjugate pairs, these resonances occur if 
$\omega_j^2 = (\pm \omega_{k} \pm \omega_{l})^2$. 
For the rest of the development presented herein, we now assume that no second-order resonance between the eigenfrequencies are present. The treatment of second-order resonances is further investigated in Sec. \ref{subsec:IRmulti} to enlarge the scope of the present developments and explicitly show how to take them into account.

Since the quadratic terms can be cancelled under this assumption, and since mechanical systems are known to show nonlinear response at finite amplitudes ({\em e.g.} hardening/softening behaviour), it is important to derive at least the third-order dynamics. 
The third-order homological equation reads:
\begin{align}\label{eq:to_homological}
	&
	\left[\begin{array}{cc}
		\Mu & \zerou \\
		\zerou & -\Iu
	\end{array}\right] 
	\left\{\begin{array}{l}
		\Upsu^{(1)}\fuu^{(3)}(\zu,\zu,\zu) \\
		\Psiu^{(1)}\fuu^{(3)}(\zu,\zu,\zu) 
	\end{array}\right\} + \nonumber \\
	&
	\left[\begin{array}{cc}
		\Mu & \zerou \\
		\zerou & -\Iu
	\end{array}\right] 
	\left\{\begin{array}{l}
		\Upsu^{(2)}(\fuu^{(2)}(\zu,\zu),\zu) + \Upsu^{(2)}(\zu,\fuu^{(2)}(\zu,\zu))  \\
		\Psiu^{(2)} (\fuu^{(2)}(\zu,\zu),\zu) + \Psiu^{(2)} (\zu,\fuu^{(2)}(\zu,\zu))
	\end{array}\right\} + \nonumber \\
	&
	\left[\begin{array}{cc}
		\Mu & \zerou \\
		\zerou & -\Iu
	\end{array}\right] 
	\left\{\begin{array}{l}
		\Upsu^{(3)} (\fuu^{(1)}\zu,\zu,\zu) + \Upsu^{(3)}(\zu,\fuu^{(1)}\zu,\zu) + \Upsu^{(3)}(\zu,\zu,\fuu^{(1)}\zu) \\
		\Psiu^{(3)} (\fuu^{(1)}\zu,\zu,\zu) + \Psiu^{(3)}(\zu,\fuu^{(1)}\zu,\zu) + \Psiu^{(3)} (\zu,\zu,\fuu^{(1)}\zu)
	\end{array}\right\} + \nonumber \\
	&
	\left[\begin{array}{cc}
		\zerou & \Ku \\
		\Iu & \zerou
	\end{array}\right]
	\left\{\begin{array}{l}
		\Upsu^{(3)} (\zu,\zu,\zu)\\
		\Psiu^{(3)} (\zu,\zu,\zu)
	\end{array}\right\} =  \nonumber \\
	&
	\left\{\begin{array}{c}
		-\Gu (\Psiu^{(2)}(\zu,\zu),\Psiu^{(1)}\zu)-\Gu (\Psiu^{(1)}\zu,\Psiu^{(2)}(\zu,\zu))-\Hu (\Psiu^{(1)}\zu,\Psiu^{(1)}\zu,\Psiu^{(1)}\zu)\\
		\mathbf{0}
	\end{array}\right\},
\end{align}
which  is reported in an extended form to highlight its hierarchical structure.  Again, the second part of Eq.~\eqref{eq:to_homological} provides a direct relationship between the velocity mapping $\Upsu$ and its counterpart for displacements $\Psiu$. This result is in fact true for any homological equation in direct form of a given order:  every velocity mapping $\Upsu^{(n)}$ is a linear function of displacement mappings $\Psiu^{(m)}$ with $m\leq n$. This is in accordance with the results provided by the real-valued mapping introduced in~\cite{artDNF2020}, where equivalent findings were reported yet not systematically proven.\\

Eq.~\eqref{eq:to_homological} can be simplified under the assumption of no second-order resonances. In that case, one has $\fuu^{(2)}=0$, such that a compact and simplified third-order homological equation can be expressed as:
\begin{subequations} 
\label{eq:to_short_short}
	\begin{align}
	& \forall\,k,l,m = 1,...,2\Nr, \nonumber \\
	&
	\sum_{s=1}^{2\Nr}\left[ \left( \lambda_{s}+\lambda_{k}+\lambda_{l}+\lambda_{m} \right)\Mu \right] f^{(3)}_{sklm} \Psiu^{(1)}_{s} + 
	\left[ \left( \lambda_{k}+\lambda_{l}+\lambda_{m} \right)^{2}\Mu  + \Ku \right] \Psiu^{(3)}_{klm} = \nonumber \\
	&
	- \Gu(\Psiu^{(2)}_{kl},\Psiu^{(1)}_{m}) - \Gu(\Psiu^{(1)}_{k},\Psiu^{(2)}_{lm}) - \Hu(\Psiu^{(1)}_{k},\Psiu^{(1)}_{l},\Psiu^{(1)}_{m}), \label{eq:to_short_shorta}\\
	&
	\Upsu^{(3)}_{klm} = (\lambda_{k}+\lambda_{l}+\lambda_{m})\Psiu^{(3)}_{klm} + \sum_{s=1}^{2\Nr}\Psiu^{(1)}_sf^{(3)}_{sklm}.  \label{eq:to_short_shortb}
	\end{align}
\end{subequations}
One can notice in particular that Eq.~\eqref{eq:to_short_shortb} gives the explicit relationship between velocity and displacement mappings, while Eq.~\eqref{eq:to_short_shorta} has been rewritten as a function of displacement mapping only, {\em i.e.} by using Eq.~\eqref{eq:to_short_shortb} to eliminate $\Upsu$ terms.\\
Let us define $\Xiu^{(3)}_{klm}$ as minus the right hand side of Eq. \eqref{eq:to_short_shorta} to improve readability of following equations:
\begin{equation}
	\Xiu_{klm}^{(3)} = \Gu(\Psiu^{(2)}_{kl},\Psiu^{(1)}_{m}) + \Gu(\Psiu^{(1)}_{k},\Psiu^{(2)}_{lm}) + \Hu(\Psiu^{(1)}_{k},\Psiu^{(1)}_{l},\Psiu^{(1)}_{m}).
\end{equation}
From Eq.~\eqref{eq:to_short_short}, third order nonlinearities can be set  to zero  ($f^{(3)}_{sklm}=0$) whenever the map $\Psiu^{(3)}$ can be computed as:
\begin{align}
\label{eq:to_map}
	\Psiu^{(3)}_{klm} =   -\left[ \left( \lambda_{k}+\lambda_{l}+\lambda_{m} \right)^{2}\Mu  + \Ku \right]^{-1} \Xiu^{(3)}_{klm},\quad \forall\,k,l,m = 1,...,2\Nr. 
\end{align}
Contrary to Eq.~\eqref{eq:so_map}, there are always $(k,l,m)$ combinations such that the resulting system is singular regardless of the associated eigenfrequency, because of the purely imaginary complex conjugate spectrum. This is for instance observed for $(k,l,m)$ combinations of the type $(s,s,mod(s+\Nr,2\Nr))$ with $mod(\cdot,\cdot)$ modulo operation. Resonance conditions that do not depend on the eigenfrequency value are called {\em trivial resonances}. Their presence and consequences for the development of reduced models is detailed in Sec.~\ref{sec:ROM}, see also~\cite{TouzeCISM} for general considerations.\\

In the present work, application to large-scale models are reported by truncating the mappings at second-order, {\em i.e.} $(\Psiu^{(n)},\Upsu^{(n)})=0,\, \forall\; n\geq3$. Dynamical coefficients are truncated at third-order, {\em i.e.}  $\fuu^{(n)}=0,\, \forall\; n\geq4$. These assumptions lead to the second-order DNF method as introduced in~\cite{artDNF2020}, which has already proven to give accurate results. Extensions to higher orders are possible, see {\em e.g.}~\cite{PONSIOEN2018}; however for the present study we restrict ourselves to this second-order DNF in Section~\ref{sec:results} where numerical examples are presented.

The complete normal form of the problem can be written by expressing the third-order coefficients $\fuu^{(3)}$. Due to the complexity of the calculations, it is easier to present them by using the modal basis. Under the hypothesis of no second-order resonances, {\em i.e.} $\fuu^{(2)}=0$, the general third-order dynamics can be written  by projecting Eq.~\eqref{eq:to_short_short}  onto the modal basis, hence yielding:
\begin{subequations}\label{eq:to_proj}
	\begin{align}
		& \forall\,k,l,m = 1,...,2\Nr,\nonumber\\
		& \Phiu^{\mathrm{T}}\sum_{s=1}^{2\Nr} \lambda_{s}\Mu  \Psiu_s^{(1)}f^{(3)}_{sklm} =  -\Phiu^{\mathrm{T}}\Xiu^{(3)}_{klm},  \\
		& \Phiu^{\mathrm{T}}\sum_{s=1}^{2\Nr}\Mu \Psiu_s^{(1)}f^{(3)}_{sklm} = 0,
	\end{align}
\end{subequations}
from which the following relationships are retrieved, thus providing explicit expressions of $f^{(3)}_{sklm}$:
\begin{subequations}
\label{eq:f3eq}
	\begin{align}
		& \forall\,k,l,m = 1,...,2\Nr, \, \forall\,s = 1,...,\Nr, \nonumber \\
		& (\lambda_{s} - \lambda_{s+\Nr}) f^{(3)}_{sklm} =  -\phiu_s^{\mathrm{T}}\Xiu^{(3)}_{klm},  \\
		& f^{(3)}_{sklm} + f^{(3)}_{(s+\Nr)klm} = 0.
	\end{align}
\end{subequations}
One can note in particular that the previous equations can be used to write the complete normal form of the problem in the $2\Nr$-dimensional setting, and thus they involve the full matrix of eigenvectors $\Phiu$. However these operations are never computed in a reduction perspective, as it will be explained in Section~\ref{sec:ROM}. Indeed, reduction will be performed by selecting a few master coordinates so that $\Phiu$ will be restricted to the selected master eigenvectors only. Examples will be provided with direct reduction to single-mode and multi-mode dynamics.\\

Following the procedure developed within the present section, the method of normal transform can be expanded up to arbitrary order, thus allowing algorithmic implementations, such as those reported in~\cite{PONSIOEN2018}. A brief summary of higher order expansions, highlighting the  hierarchical structure of the procedure, is reported in Appendix~\ref{sec:hot} even if, as recalled, in the present work, only the second-order DNF is used as it proves highly efficient for the applications addressed in Section~\ref{sec:results}.\\

Finally, it is worth stressing that the generic structure of the method can be summarised in three steps: for each order $n$ of the expansion, the homological equation of order $n$ is written, thus providing mappings $\Psiu^{(n)}$, $\Upsu^{(n)}$ and dynamics coefficients $\fuu^{(n)}$ that cannot be set to zero. Then, $\fuu^{(n+1)}$ needs to be computed through projection of the homological equation of order $n+1$, since the order $n$ terms will create new monomials at order $n+1$ that need to be tracked properly. The procedure then starts again by identifying mappings and non zero terms at the next order.
Implementation details of the method up to second order are reported in Appendix~\ref{sec:algo}, and the case of second-order internal resonance will be specifically detailed in Section~\ref{subsec:IRmulti}, underlining how the different steps of the method may create higher-order terms in the asymptotics.

\section{Model Order Reduction}
\label{sec:ROM}

In this section, the reduction method is explained in its general settings, then two particular cases are investigated for the sake of clarity. First, reduction to a single master mode is detailed in Section~\ref{subsec:rom1master}, while the case of second-order internal resonance, needing extra calculations and at least two master coordinates, is developed in Section~\ref{subsec:IRmulti}.

\begin{figure*}
\centering
	\begin{tikzpicture}
    	\node[anchor=south west] at (0,0) {\includegraphics[width=0.99\textwidth]{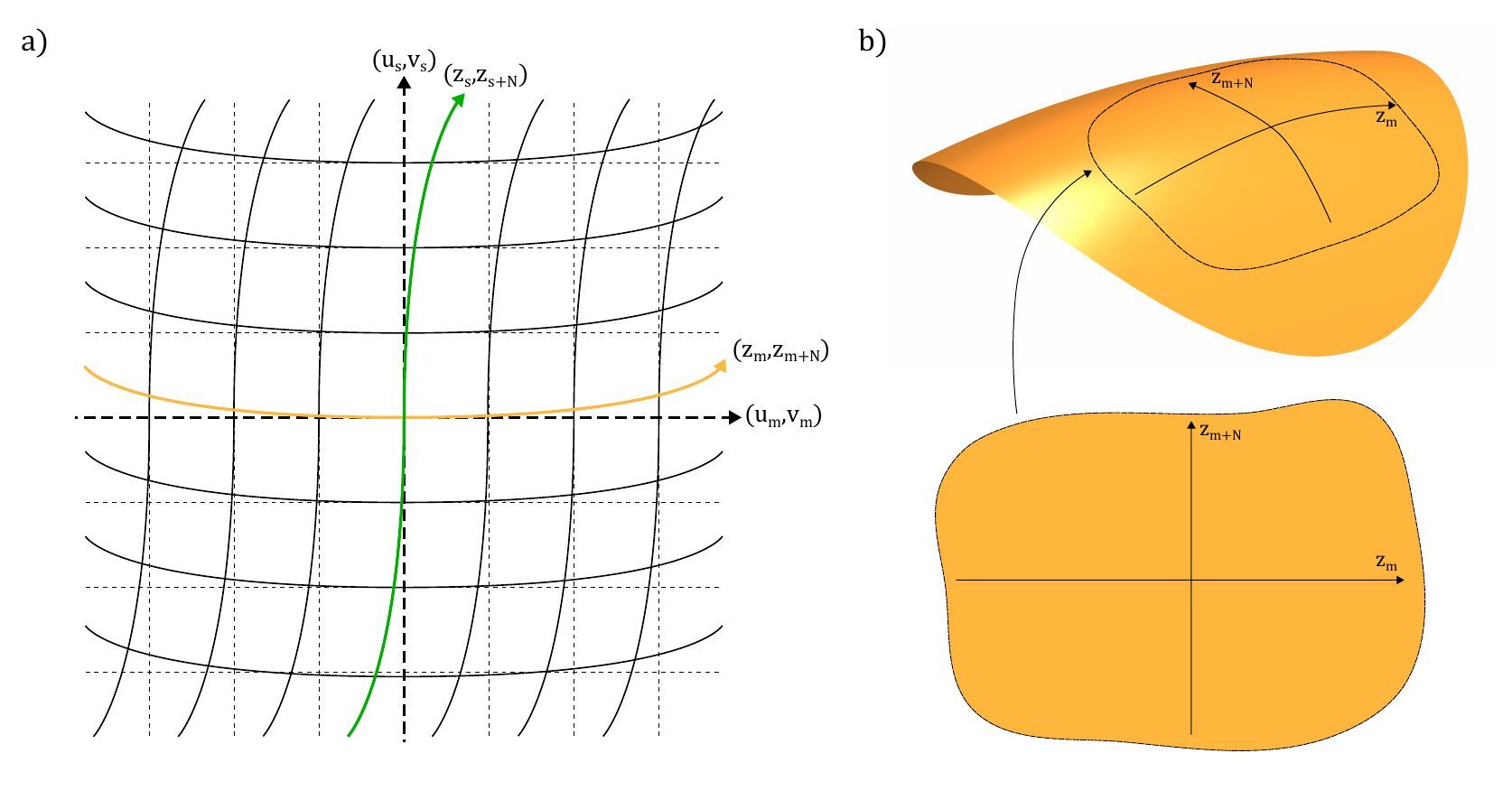}};
    	\begin{scope}
    	\node[text width=2cm] at (11.9,4.85) {$\Psiu(\zu),\Upsu(\zu)$};
    	\end{scope}
	\end{tikzpicture}
	\caption{Schematic representation of the nonlinear mappings in phase space. (a) Effect of the normal transform, expressing the dynamics in an invariant-based span: $(\ur_s,\vr_s)$ and $(\ur_m,\vr_m)$ are the modal coordinates (orthogonormal basis), $(\zr_s,\zr_{s+\Nr})$ and $(\zr_m,\zr_{m+\Nr})$ the associated {\em normal} coordinates. Orange and Green lines represent the invariant manifolds (LSM). (b) Schematic representation of the two-dimensional invariant manifold associated to $\phiu_m$, and effect of the nonlinear mappin and the parametrisation. Picture (a) is a rework of the image reported in \cite{touze03-NNM}.}\label{fig:invariance}
\end{figure*}

\subsection{Reduction and normal form style parametrisation}
\label{subsec:reduceNF}

The normal form procedure as exposed in Section~\ref{sec:theory} is a complete nonlinear change of coordinates where the number of input and output variables is the same. As underlined in~\cite{touze03-NNM}, reduction can then be applied by selecting only a subset of few master normal variables, and cancelling all the others, such that the size of the $\zu$ vector falls from $2\Nr$ to $2m$ with $m \ll \Nr$ the number of master coordinates.
By doing so, one modifies the nonlinear one-to-one  diffeomorphism to a parametrisation procedure of the underlying invariant manifold (or Lyapunov subcenter manifold (LSM) in a conservative context), the existence of which is guaranteed by the Lyapunov center theorem~\cite{Lyapunov1907,Kelley69}. The derived expressions are then identical to the one found by first assuming a parametrisation, then solving the resulting tangent and normal homological problems using the normal form style, as stated in~\cite{Haro}.

In the framework of mechanical systems, the parametrisation of invariant manifolds has been extensively treated in ~\cite{Haller2016,PONSIOEN2018}; the application of the normal form approach to model order reduction is instead adopted in the developments led in~\cite{touze03-NNM,TOUZE:JSV:2006,Wagg2019}.
\\

A geometrical interpretation  is reported in Fig.~\ref{fig:invariance}(a) where a two-dimensional representation of the four-dimensional subspace associated to two modes $\phiu_s$ and $\phiu_m$ is reported. With modal coordinates, the dynamics is expressed within an orthogonal basis formed by the eigenmodes. Due to the presence of invariant-breaking terms~\cite{TouzeCISM}, linear eigensubspaces lose the invariance property. Invariant manifolds (LSM) emerging from modal spectral subspaces are highlighted in red and green respectively, their curvature expressing the non-resonant coupling between modes and thus the amplitude-dependence of modal quantities in a nonlinear framework. The normal form allows expressing the nonlinear dynamics in an invariant-based span of the phase space. Once that  obtained, reduction to an invariant subspace is simply given by cancelling slave normal variables.  Fig. \ref{fig:invariance}(b) illustrates the application of the nonlinear mappings. Thanks to normal form transformation, non-resonant coupling terms are removed, hence the invariant manifold (LSM) tangent at the origin to a given mode lies on the plane identified by the $(\zr_m,\zr_{m+\Nr})$ pair with associated $\Psiu^{(1)}_{m}=\Psiu^{(1)}_{m+\Nr}=\phiu_m$. Therefore, the normal transform allows expressing the dynamics with normal coordinates that are nonlinearly related to the original ones through the nonlinear mappings $\Psiu(\zu)$, $\Upsu(\zu)$.\\
In a more general framework, reduction is achieved by selecting a set $\Phiu_{m}$ of master modes and identifying the set $\cZ$ of indices such that $\Psiu^{(1)}_s\in\Phiu_{m}$. All other coordinates are set to zero:
\begin{equation}
\forall\; s \not\in \mathcal{Z}, s=1,...,2\Nr: \quad \zr_{s} = 0,
\end{equation}
which yields the following mapping expansions in indicial notation:
\begin{subequations}
	\begin{align}
		\Uu = \sum_{s\in\cZ} \Psiu^{(1)}_{s}\zr_s + \sum_{k,l\in\cZ}\Psiu^{(2)}\zr_k\zr_l + \sum_{k,l,m\in\cZ}\Psiu^{(3)}\zr_k\zr_l\zr_m + O(\|\zu\|^{4}),\\
		\Vu = \sum_{s\in\cZ} \Upsu^{(1)}_{s}\zr_s + \sum_{k,l\in\cZ}\Upsu^{(2)}\zr_k\zr_l + \sum_{k,l,m\in\cZ}\Upsu^{(3)}\zr_k\zr_l\zr_m + O(\|\zu\|^{4}),
	\end{align}
\end{subequations}
and reduced dynamics:
\begin{equation}
	\dot{\zr}_s = \lambda_s\zr_s + \sum_{k,l\in\cZ}f^{(2)}_{skl}\zr_k\zr_l + \sum_{k,l,m\in\cZ}f^{(3)}_{sklm}\zr_k\zr_l\zr_m + O(\|\zu\|^{4}), \quad \forall\;s\in\cZ.
\end{equation}
Note in particular that under the assumption of no internal resonance then the quadratic term vanishes, $f^{(2)}_{skl}=0$, and explicit expressions for $f^{(3)}_{sklm}$ are given in Eq.~\eqref{eq:f3eq}. Also, the reduction procedure implies that only quantities associated to the master modes are required during computation, so that the complete eigenfunctions matrix $\Phiu$ is not required.\\
In the remainder of the paper, the case of the second-order DNF is selected. With this choice only second-order internal resonances have to be tracked. On the other hand, the dynamics up to order three contains all the monomials, be they resonant or non-resonant. Applications of the third-order mapping would lead to  cancellation of the non-resonant third-order monomials, but it would create higher-order terms. Keeping the development at second-order has the advantage that all higher-order internal resonance, starting from third-order, do not need a special attention since they are all preserved in the reduced dynamics. The next Section details the case of a single master mode while Sec.~\ref{subsec:IRmulti} is focused on the treatment of second-order internal resonance.

\subsection{Single Mode Reduction}\label{subsec:rom1master}

In this Section, the simplest case of a reduction to a single master coordinate is detailed, which is performed under the  assumption of no internal resonance. Let us denote as $\phiu_m$ the master mode of interest, the master normal coordinates are thus selected as $(\zr_{m},\zr_{m+\Nr})$. The reduction procedure consists in setting to zero all other coordinates:
\begin{equation}
\forall p\neq m, p=1,...,\Nr: \quad \zr_{p} = 0, \quad \zr_{p+\Nr} = 0.
\end{equation}
Let us set $\alpha=m$ and $\beta=m+\Nr$. The associated displacement expansion adopting a second-order mapping reads:
\begin{align}\label{eq:single_dof_map}
	\Uu =  \Psiu^{(1)}_{m}(\zr_{\alpha}+\zr_{\beta}) + \Psiu_{\alpha\alpha}^{(2)}\zr_{\alpha}\zr_{\alpha} +  \Psiu_{\beta\beta}^{(2)}\zr_{\beta}\zr_{\beta} + \left(\Psiu_{\alpha\beta}^{(2)}+\Psiu_{\beta\alpha}^{(2)}\right)\zr_{\alpha}\zr_{\beta},
\end{align}
with $\Psiu^{(1)}_{m}=\phiu_m$. Eq.~\eqref{eq:single_dof_map} represents the parametrisation of the invariant manifold (LSM), as shown in~\cite{Haller2016,PONSIOEN2018}. The tensors are computed from Eq.~\eqref{eq:so_map}, yielding the following relationships:
\begin{subequations}
	\begin{align}
		&\, \Psiu_{\alpha\alpha}^{(2)} = \Psiu_{\beta\beta}^{(2)} =  -\left[ (\lambda_\alpha+\lambda_\alpha)^{2}\Mu+\Ku \right]^{-1}\Gu\left(\Psiu^{(1)}_{m},\Psiu^{(1)}_{m}\right),
 \\
 		&\, \Psiu_{\alpha\beta}^{(2)} = \Psiu_{\beta\alpha}^{(2)} =  -\left[ (\lambda_\alpha+\lambda_{\beta})^{2}\Mu+\Ku \right]^{-1}\Gu\left(\Psiu^{(1)}_{m},\Psiu^{(1)}_{m}\right).
	\end{align}
\end{subequations}
with $\lambda_{m}=\img\omega_{m}$, $\lambda_{m+\Nr}=-\img\omega_{m}$. Third-order reduced dynamics coefficients $\fuu^{(3)}$ are obtained as in Eq.~\eqref{eq:to_proj}. Since only $\zr_\alpha$, $\zr_\beta$ are different from zero, then only the projection on $\phiu_m$ is required:
\begin{subequations}
	\begin{align}
		& \forall\,k,l,p = \{\alpha,\beta\},\nonumber \\
		& \phiu_{m}^{\mathrm{T}}\sum_{s\in\cZ}\left[ \lambda_{s}\Mu \right] \Psiu_{s}^{(1)}f^{(3)}_{sklp} = -\phiu_{m}^{\mathrm{T}}\Xiu^{(3)}_{klp},  \\
		& \phiu_{m}^{\mathrm{T}}\sum_{s\in\cZ}\Mu \Psiu_{s}^{(1)}f^{(3)}_{sklp} = 0,
	\end{align}
\end{subequations}
which yields the following relation to compute the reduced dynamics coefficients:
\begin{subequations}
	\begin{align}
		& \forall\,k,l,p = \{\alpha,\beta\}, \nonumber \\
		& (\lambda_{\alpha} - \lambda_{\beta}) f^{(3)}_{\alpha klp} =  -\phiu_{m}^{\mathrm{T}}\Xiu^{(3)}_{klp},  \\
		& f^{(3)}_{\alpha klp} + f^{(3)}_{\beta klp} = 0.
	\end{align}
\end{subequations}
The associated third-order reduced dynamics reads:
\begin{align}\label{eq:single_dof_red_dyn}
	&\, \left\{\begin{array}{l}
		\dot{\zr}_{\alpha} \\
		\dot{\zr}_{\beta}
	\end{array}\right\} = 
	\left[\begin{array}{cc}
		\img\omega_{m} & 0 \\
		0 & -\img\omega_{m}
	\end{array}\right] 
	\left\{\begin{array}{l}
		{\zr}_{\alpha} \\
		{\zr}_{\beta}
	\end{array}\right\} + \nonumber \\
	&\,
	\left\{\begin{array}{l}
		{f}_{\alpha\alpha\alpha\alpha}^{(3)} \\
		{f}_{\beta\alpha\alpha\alpha}^{(3)}
	\end{array}\right\}{\zr}_{\alpha}^{3} + 
	\left\{\begin{array}{l}
		{f}_{\alpha\alpha\alpha\beta}^{(3)}+{f}_{\alpha\alpha\beta\alpha}^{(3)}+{f}_{\alpha\beta\alpha\alpha}^{(3)} \\
		{f}_{\beta\alpha\alpha\beta}^{(3)} + {f}_{\beta\alpha\beta\alpha}^{(3)} + {f}_{\beta\beta\alpha\alpha}^{(3)}
	\end{array}\right\}{\zr}_{\alpha}^{2}{\zr}_{\beta} +  \nonumber \\
	&\,
	\left\{\begin{array}{l}
		{f}_{\alpha\alpha\beta\beta}^{(3)} + {f}_{\alpha\beta\alpha\beta}^{(3)} + {f}_{\alpha\beta\beta\alpha}^{(3)} \\
		{f}_{\beta\alpha\beta\beta}^{(3)} + {f}_{\beta\beta\alpha\beta}^{(3)} + {f}_{\beta\beta\beta\alpha}^{(3)}
	\end{array}\right\}{\zr}_{\alpha}{\zr}_{\beta}^{2} + 
	\left\{\begin{array}{l}
		{f}_{\alpha\beta\beta\beta}^{(3)} \\
		{f}_{\beta\beta\beta\beta}^{(3)}
	\end{array}\right\}{\zr}_{\beta}^{3}.
\end{align}
Since $\Gu$ is symmetric with respect to permutations of its arguments,  the reduced dynamics can also be written as:
\begin{align}\label{eq:sart_point_rv_sdof}
	&\,
	\left\{\begin{array}{l}
		\dot{\zr}_{\alpha} \\
		\dot{\zr}_{\beta}
	\end{array}\right\} = 
	\left[\begin{array}{cc}
		\img\omega_{m} & 0 \\
		0 & -\img\omega_{m}
	\end{array}\right] 
	\left\{\begin{array}{l}
		{\zr}_{\alpha} \\
		{\zr}_{\beta}
	\end{array}\right\} + \nonumber \\
	&\,
	\left\{\begin{array}{l}
		-\frac{1}{2\img\omega_{m}}\phiu_{m}^{\mathrm{T}}\left[ \Gu(\Psiu^{(2)}_{\alpha\alpha}+\Psiu^{(2)}_{\alpha\beta},\Psiu^{(1)}_{m})  \right]  \\
		+\frac{1}{2\img\omega_{m}}\phiu_{m}^{\mathrm{T}}\left[ \Gu(\Psiu^{(2)}_{\alpha\alpha}+\Psiu^{(2)}_{\alpha\beta},\Psiu^{(1)}_{m})  \right]
	\end{array}\right\}({\zr}_{\alpha}+{\zr}_{\beta})^{3} + \nonumber \\
	&\,
	\left\{\begin{array}{l}
		-\frac{1}{2\img\omega_{m}}\phiu_{m}^{\mathrm{T}}\left[  \Hu(\Psiu^{(1)}_{m},\Psiu^{(1)}_{m},\Psiu^{(1)}_{m}) \right]  \\
		+\frac{1}{2\img\omega_{m}}\phiu_{m}^{\mathrm{T}}\left[  \Hu(\Psiu^{(1)}_{m},\Psiu^{(1)}_{m},\Psiu^{(1)}_{m}) \right]
	\end{array}\right\}({\zr}_{\alpha}+{\zr}_{\beta})^{3} + \nonumber \\
	&\,
	\left\{\begin{array}{l}
		-\frac{1}{2\img\omega_{m}}\phiu_{m}^{\mathrm{T}}\left[ \Gu(\Psiu^{(2)}_{\alpha\alpha}-\Psiu^{(2)}_{\alpha\beta},\Psiu^{(1)}_{m})  \right] \\
		+\frac{1}{2\img\omega_{m}}\phiu_{m}^{\mathrm{T}}\left[ \Gu(\Psiu^{(2)}_{\alpha\alpha}-\Psiu^{(2)}_{\alpha\beta},\Psiu^{(1)}_{m}) \right]
	\end{array}\right\}({\zr}_{\alpha}+{\zr}_{\beta})({\zr}_{\alpha}-{\zr}_{\beta})^{2}.
\end{align}
This last equation represents the nonlinear dynamics, up to the third-order, of the system along the corresponding $m^{\mathrm{th}}$ invariant manifold, with the key feature that all coefficients can be computed directly with operations from the FE dofs in physical space. It also underlines that the derived formulas are explicit enough so that a non-intrusive implementation of all calculations can be targeted. Indeed, evaluation of $\Gu$ and $\Hu$ terms does not require an explicit decomposition of the internal force in their linear, quadratic and cubic parts, but they can be evaluated for instance through non-intrusive methods, retrieving the direct algorithm proposed in~\cite{artDNF2020} with a real formalism.
This approach is particularly appealing for application of the present method with commercial finite element software since there is no special need to implement finite element routines to obtain the ROMs given by the DNF approach. More details about the potential non-intrusive implementation of the method are provided in Sec.~\ref{sec:real}.


\subsection{Second-order Internal Resonance and Multiple Master Modes}
\label{subsec:IRmulti}

In this Section, more detailed explanations on handling the case of a second-order internal resonance are given. This case has not been treated in earlier developments of the normal form approach for reduced-order modelling in~\cite{touze03-NNM,TOUZE:JSV:2006,artDNF2020}. A difficulty relies in the fact that resonant monomials will stay in the second-order normal form, since one cannot simply cancel all quadratic terms by stating $f^{(2)}_{skl}=0,\,\forall\, s,k,l=1,...,2\Nr$ as in Sec.~\ref{subsec:homolog}. In turn, second-order components of the mappings are not solvable anymore, see  Eq.~\eqref{eq:so_map}, since a second-order internal resonance makes the matrix to invert singular. The consequence of these two facts is that new cubic terms will arise in the third-order reduced dynamics and they must be properly computed.\\

As anticipated in Sec. \ref{subsec:homolog}, since the spectrum is composed by purely imaginary conjugate pairs, a second-order internal resonance stems from any combination of three frequencies of the type:
\begin{equation}\label{eq:secondorderIRgen}
	(\pm\omega_k\pm\omega_l)^{2}= (\pm\omega_{r})^{2},\quad \forall\,k,l,r=1,...,\Nr.
\end{equation}
Or by using the entries of $\Lambdau$:
\begin{equation}
	(\lambda_k  +\lambda_l)^2 = -\omega_r^{2}, \quad \forall\;k,l=1,...,2\Nr,\;\mathrm{and}\;\forall\;r=1,...,\Nr.
\end{equation}
If this condition is met, then Eq.~\eqref{eq:so_map} cannot be solved since the operation $[(\lambda_k  +\lambda_l)^2\Mu+\Ku]$ filters $\omega_r^{2}$ from the eigenspectrum, and the resulting matrix is singular. Furthermore, the number of linearly dependent columns of the system is equal to the number of modes with eigenfrequency $\omega_r$. This is better highlighted by projecting to the modal basis via:
\begin{equation}
	\Phiu^{\mathrm{T}}[-\omega^{2}_{r}\Mu+\Ku]\Phiu = -\omega_r^2\Iu + \Omegau^{2},
\end{equation}
showing clearly that any diagonal entry of $\Omegau^{2}$ equal to $\omega_r^{2}$ is cancelled. Hence, all $(k,l)$ such that Eq.~\eqref{eq:secondorderIRgen} is fulfilled must be treated with care.

The treatment of internal resonance with a direct approach operating in physical space is more difficult as compared to normal form computation operating in the modal basis. Indeed, the tracking of monomials associated to the resonance is immediate in the modal basis, and the known remedy consisting in cancelling the term in the mapping, and thus keeping the resonant monomial in the reduced dynamics, can be transparently applied~\cite{touze03-NNM,Haro,PONSIOEN2018}. With a direct approach in physical space,  this is no longer possible since a single $\Psiu^{(2)}_{kl}$ vector embeds terms required to cancel all $kl$ monomials. This implies that if one sets crudely $\Psiu^{(2)}_{kl}=0$, then all $\zu$ coordinates would remain coupled by the resulting $f^{(2)}_{skl}$ monomials, whereas the correct strategy consists in tracking only the terms that are involved in the resonance. A proper approach is to force the solution of the system  to be orthogonal to the kernel of the matrix to invert in Eq.~\eqref{eq:so_map}. This requires identifying all triplets $(r,k,l)$ involved in the resonance relationship \eqref{eq:secondorderIRgen}, and enforcing $\Psiu^{(2)}_{kl}$ and $\Upsu^{(2)}_{kl}$ to be mass-orthogonal to the set of involved modes. Let us define as $\Phiu_{\mathrm{R}}$ the set of eigenmodes such that a given pair of indices $(k,l)$ satisfies a resonance condition, so that the enforced orthogonality condition reads:
\begin{subequations}\label{eq:orth_condition}
\begin{align}
	&\phiu_{r}^{\mathrm{T}}\Mu\Psiu^{(2)}_{kl}=\Psi^{(2)}_{rkl}=0,\quad \forall\;\phiu_r\in\Phiu_{\mathrm{R}},\\
	&\phiu_{r}^{\mathrm{T}}\Mu\Upsu^{(2)}_{kl}=\Upsilon^{(2)}_{rkl}=0,\quad \forall\;\phiu_r\in\Phiu_{\mathrm{R}}.
\end{align}
\end{subequations}
Since the mapping is imposed to be orthogonal to $\Phiu_{\mathrm{R}}$, resonant monomials that remain in the reduced dynamics are obtained by projecting Eqs.~\eqref{eq:so_homologicalPsi} and \eqref{eq:so_homologicalUpsilon} on each $\phiu_r\in\Phiu_{\mathrm{R}}$. In the context of model order reduction, only the subset of modes defined by $\cZ^{(1/2)}$, with $\cZ^{(1/2)}$ subset of indices of $\cZ$ lower or equal to $\Nr$, is considered, hence coefficients are estimated as:
\begin{subequations}
	\begin{align}
		& \forall\; k,l\in\cZ, \,\forall\; r\in\cZ^{(1/2)} :(\lambda_k+\lambda_l)^{2}=-\omega_r^2,\; \mathrm{and} \; \,\phiu_r\in \Phiu_{\mathrm{R}}, \nonumber \\
		& (\lambda_{r}-\lambda_{(r+\Nr)}) f_{rkl}^{(2)} = -\phiu_r^{\mathrm{T}} \Gu(\Psiu_{k}^{(1)},\Psiu_{l}^{(1)}), \\
		& f_{rkl}^{(2)}+f_{(r+\Nr)kl}^{(2)}=0.
	\end{align}
\end{subequations}
Implementation details of the orthogonality condition are reported in Appendix~\ref{sec:algo}.\\
The presence of $\fuu^{(2)}$ in Eq.~\eqref{eq:red_dyn} modifies the third order homological equation given by Eq.~\eqref{eq:to_short_short}, hence one needs to rely on its full form. In the present work, second-order DNF is used, hence $\Psiu^{(3)}=0$ and the third-order homological equation is used only to compute monomials at order $n+1$, {\em{i.e.}} at third-order. Equation~\eqref{eq:to_short_short} in presence of second-order resonances modifies to:
\begin{subequations}\label{eq:to_short_full}
	\begin{align}
	& \forall\,k,l,m \in\cZ, \nonumber\\
	& \sum_{s\in\cZ}\left(\lambda_{s}\Mu	\Psiu^{(1)}_{s}f_{sklm}^{(3)} + \Mu\Upsu^{(2)}_{sm}f^{(2)}_{skl} + \Mu\Upsu^{(2)}_{ks}f^{(2)}_{slm}\right) = - \Xiu^{(3)}_{klm},\label{eq:to_short_fulla}\\
	&\,
	  \sum_{s\in\cZ}\left( \Psiu_{sm}^{(2)}f^{(2)}_{skl} + \Psiu_{ks}^{(2)}f^{(2)}_{slm} + \Psiu^{(1)}_{s}f^{(3)}_{sklm}\right) = 0,\label{eq:to_short_fullb}
	\end{align}
\end{subequations}
where new terms $\Mu\Upsu^{(2)}_{sm}f^{(2)}_{skl} + \Mu\Upsu^{(2)}_{ks}f^{(2)}_{slm}$ in Eq. \eqref{eq:to_short_fulla} and $\Psiu_{sm}^{(2)}f^{(2)}_{skl} + \Psiu_{ks}^{(2)}f^{(2)}_{slm}$ in Eq. \eqref{eq:to_short_fullb} appear. These terms highlight that low order monomials affect the estimation of higher order mapping and reduced dynamics coefficients. The latter are estimated by projecting Eq.~\eqref{eq:to_short_full} onto the modal basis:
\begin{subequations}\label{eq:to_rdyn_coeff_ir}
	\begin{align}	
		& \forall\,k,l,m \in\cZ, \forall s \in\cZ^{(1/2)},	\nonumber \\
		& \lambda_{s}f_{sklm}^{(3)} + \lambda_{(s+\Nr)}f_{(s+\Nr)klm}^{(3)} =  \sum_{p\in\cZ}\left(-\Upsilon_{spm}^{(2)}f^{(2)}_{pkl} -\Upsilon_{skp}^{(2)}f^{(2)}_{plm} \right)  - \phiu_s^{\mathrm{T}}\Xiu^{(3)}_{klm},    \\		
	&	 f_{sklm}^{(3)} + f_{(s+\Nr)klm}^{(3)} = \sum_{p\in\cZ}\left( -\Psi_{spm}^{(2)}f^{(2)}_{pkl} - \Psi_{skp}^{(2)}f^{(2)}_{plm}\right).\label{eq:ff_to_relation} 
	\end{align}
\end{subequations}
Equation \eqref{eq:to_rdyn_coeff_ir} shows that in presence of second-order resonances $f_{sklm}^{(3)}\neq -f_{(s+\Nr)klm}^{(3)}$. This has consequences when mapping the reduced dynamics to real-valued quantities. Further details are reported in Appendix \ref{sec:second_order_ir_coeff}.\\
Overall, the concepts introduced up to this point allow deriving the second-order DNF in the most general case. Second-order mappings $\Psiu^{(2)}$ are obtained from Eq. \eqref{eq:so_map} with potential application of the constraint given by Eq. \eqref{eq:orth_condition} if resonance conditions are met. Non-trivial resonance conditions in a FE framework are identified as detailed in Appendix \ref{sec:r_cond}. Third-order reduced dynamics coefficients are then computed from Eq. \eqref{eq:to_rdyn_coeff_ir}, which reduces to Eq. \eqref{eq:f3eq} if $\fuu^{(2)}=0$.\\
In order to be more specific, let us give more detail on the case of a 1:2 internal resonance. Let us assume that $\phiu_{\Ir}$ and $\phiu_{\Ir\Ir}$ satisfy the relation $\omega_{\Ir\Ir}=2\omega_{\Ir}$, such that a special treatment is needed for these two modes as compared to all others. As a result, combinations $(\lambda_{\Ir}+\lambda_{\Ir})$,  $(\lambda_{\Ir+\Nr}+\lambda_{\Ir+\Nr})$ filter $\omega_{\Ir\Ir}$ from the matrix required to compute $\Psiu^{(2)}$, and combinations $(\lambda_{\Ir}+\lambda_{\Ir\Ir+\Nr})$, $(\lambda_{\Ir\Ir}+\lambda_{\Ir+\Nr})$ filter $\omega_\Ir$, hence $\fuu^{(2)}\neq 0$. Due to the 1:2 resonance, pairs $(\zr_\Ir,\zr_{\Ir+\Nr})$ and $(\zr_{\Ir\Ir},\zr_{\Ir\Ir+\Nr})$ do not individually define invariant subspaces: a strong coupling exists between the two and a reduced-order model cannot involve only one of these two modes, they are intimately related and create a four-dimensional invariant manifold. A ROM with these two master coordinates is simply built by selecting
\begin{equation}
	\forall \;p\neq \{\Ir,\Ir\Ir\}, p=1,...,\Nr :\quad \zr_p = 0,\quad \zr_{p+\Nr} = 0.
\end{equation}

The development of the reduced model is then the application of the same procedure detailed for single master mode reduction in Sec. \ref{subsec:rom1master} but with two modes and with proper treatment of resonance conditions as detailed in the present section. By defining $\alpha=\Ir$, $\beta=\Ir\Ir$, $\gamma=\Ir+\Nr$, and $\delta=\Ir\Ir+\Nr$, the reduced dynamics is then expressed as:
\begin{align}\label{eq:two_dof_red_dyn}
	&\, \left\{\begin{array}{c}
		\dot{\zr}_{\alpha} \\
		\dot{\zr}_{\beta} \\
		\dot{\zr}_{\gamma} \\
		\dot{\zr}_{\delta} \\
	\end{array}\right\} = 
	\left[\begin{array}{cccc}
		\img\omega_{\Ir} & 0 & 0 & 0 \\
		0 &  \img\omega_{\Ir\Ir} & 0 & 0 \\
		0 & 0 & -\img\omega_{\Ir} & 0 \\
		0 & 0 & 0 & -\img\omega_{\Ir\Ir}
	\end{array}\right] 
	\left\{\begin{array}{l}
		{\zr}_{\alpha} \\
		{\zr}_{\beta} \\
		{\zr}_{\gamma} \\
		{\zr}_{\delta} \\
	\end{array}\right\} + \left\{\begin{array}{c}
		0 \\
		{f}_{\beta \alpha\alpha}^{(2)} \\
		0 \\
		{f}_{\delta \alpha\alpha}^{(2)}
	\end{array}\right\} {\zr}_{\alpha}^{2} +\nonumber \\
	&\,
	\left\{\begin{array}{c}
		0 \\
		{f}_{\beta \gamma\gamma}^{(2)} \\
		0 \\
		{f}_{\delta \gamma\gamma}^{(2)}
	\end{array}\right\} {\zr}_{\gamma}^{2}+ 
			\left\{\begin{array}{c}
		{f}_{\alpha \beta\gamma}^{(2)} + {f}_{\alpha \gamma\beta}^{(2)} \\
		0\\
		{f}_{\gamma \beta\gamma}^{(2)}  + {f}_{\gamma \gamma\beta}^{(2)}\\
		0
	\end{array}\right\} {\zr}_{\beta}{\zr}_{\gamma}+
	\left\{\begin{array}{c}
		{f}_{\alpha \alpha\delta}^{(2)} + {f}_{\alpha \delta\alpha}^{(2)} \\
		0\\
		{f}_{\gamma \alpha\delta}^{(2)}  + {f}_{\gamma \delta\alpha}^{(2)}\\
		0
	\end{array}\right\} {\zr}_{\alpha}{\zr}_{\delta}+
	O(\|\zu\|^{3}),
\end{align}
which highlights the presence of second order resonant monomials. This last equation can be mapped to real quantities to provide an oscillator-like equation. This is developed in Appendix \ref{sec:second_order_ir_coeff}.

\section{Real-Valued Mappings and Reduced Dynamics}
\label{sec:real}

In Secs.~\ref{sec:theory}-\ref{sec:ROM} the system is parametrised using a complex-valued mapping in $\zu$ and the resulting reduced dynamics is complex-valued as well, as a result of the initial choices with a starting point at first-order (state-space formulation). A different point of view had been developed in~\cite{touze03-NNM,TOUZE:JSV:2006,artDNF2020}, where the choice of real-valued mappings and real-valued reduced dynamics was enforced all along the calculations in order to fit the more standard oscillator equations for vibratory systems. This section is devoted to clarify the link between the two approaches by giving formulas allowing to pass from one representation to another. Mapping the reduced dynamics to real-valued quantities is efficient since it presents results that are easier to interpret from a physical point of view. Furthermore, solution of a real-valued system is often preferred.\\

Real-valued parametrisation is based on the introduction of two real {\em normal} coordinates, namely the  displacement $\ru$ and velocity $\su$, which are tangent at the origin to modal displacements and velocities. This choice is not unique but provides an immediate link with quantities that are more familiar to engineers, {\em i.e.} modal quantities. Following the direct mapping proposed in~\cite{artDNF2020}, one can write:
\begin{subequations}
	\begin{align}
		\Uu = \hat{\Psiu}(\ru,\su), \\
		\Vu = \hat{\Upsu}(\ru,\su),
	\end{align}
\end{subequations}
with $\hat{\Psiu}(\ru,\su)$, $\hat{\Upsu}(\ru,\su)$ polynomial mappings in $\ru$, $\su$. Truncation  at second-order yields:
\begin{subequations}\label{eq:rv_map}
	\begin{align}
		\Uu =&\, \sum_{k=1}^{\Nr}\phiu_{k} \rr_k + \sum_{k,l=1}^{\Nr}\left( \hat{\au}_{kl}\rr_k\rr_l + \hat{\bu}_{kl}\sr_k\sr_l + \hat{\cu}_{kl}\rr_k\sr_l\right), \\
		\Vu =&\, \sum_{k=1}^{\Nr}\phiu_{k} \sr_k + \sum_{k,l=1}^{\Nr}\left( \hat{\alphau}_{kl}\rr_k\rr_l + \hat{\betau}_{kl}\sr_k\sr_l + \hat{\gammau}_{kl}\rr_k\sr_l\right) ,
	\end{align}
\end{subequations}
with $\hat{\au}_{kl}$, $\hat{\bu}_{kl}$, $\hat{\cu}_{kl}$, $\hat{\alphau}_{kl}$, $\hat{\betau}_{kl}$, $\hat{\gammau}_{kl}$ unknown maps. 
Equation~\eqref{eq:rv_map} is equivalent to Eq.~\eqref{eq:mapping_nl} even if the derivation of direct expressions for mappings is somehow more involved due to the anti-diagonal structure of the linear part of the associated 
dynamics:
\begin{align}\label{eq:rv_red_dyn}
	\left\{\begin{array}{l}
		\dot{\su} \\
		\dot{\ru}
	\end{array}\right\}
	=
	\left[\begin{array}{cc}
 		\zerou &  -\Omegau^{2} \\
		\Iu &  \zerou
	\end{array}\right] 
	\left\{\begin{array}{l}
		\su \\
		\ru
	\end{array}\right\}
	+ O(\|\ru,\su \|^{2}).
\end{align}
A link between $\zu$ and the couple $(\ru$, $\su)$ can be established in analogy with the linear transform that relates $\pu$ to modal displacements $\uu$ and velocities $\vu$, as reported in  Appendix~\ref{sec:diag}. 
By letting $\yu=\{ \su,\ru \}$, one has:
\begin{equation}\label{eq:z_to_rs_map}
	\yu = \Ru\zu, \qquad \mbox{with} \qquad \Ru = \left[\begin{array}{cc}
		\img\Omegau & -\img\Omegau \\
		\Iu & \Iu
	\end{array}\right] ,
\end{equation}
which yields:
\begin{subequations}\label{eq:z_to_rs}
	\begin{align}
		\zr_s =&\, \frac{1}{2}\left( \rr_s - \img \frac{\sr_s}{\omega_{s}}  \right), \qquad \forall s = 1,...,\Nr,
	\\
		\zr_{s+\Nr} =&\, \frac{1}{2}\left( \rr_s + \img \frac{\sr_s}{\omega_{s}}  \right), \qquad \forall s = 1,...,\Nr.
	\end{align}
\end{subequations}
Focusing on displacement mappings, the link between the complex and real-valued expressions at second-order is retrieved by introducing the following split:
\begin{align}\label{eq:so_map_split}
	\Psiu^{(2)}(\zu,\zu) = &\, \sum_{k,l=1}^{\Nr} \left( \Psiu^{(++)}_{kl}\zr_k\zr_l + \Psiu^{(--)}_{kl}\zr_{k+\Nr}\zr_{l+\Nr} +  \Psiu^{(+-)}_{kl}\zr_k\zr_{l+\Nr} + \Psiu^{(-+)}_{kl}\zr_{k+\Nr}\zr_{l}\right).
\end{align}
Let us define vectors $\Psiu^{(\mathrm{P})}_{kl}$ and $\Psiu^{(\Nr)}_{kl}$ such that
\begin{subequations}
	\begin{align}
		& \Psiu^{(\mathrm{P})}_{kl} = \Psiu^{(++)}_{kl} = \Psiu^{(--)}_{kl} = 
		\left[ (+\omega_{k}+\omega_{l})^{2}\Mu-\Ku \right]^{-1}\Gu\left(\phiu_k,\phiu_l\right), \\
		& \Psiu^{(\Nr)}_{kl} = \Psiu^{(+-)}_{kl} = \Psiu^{(-+)}_{kl}  =  \left[ (+\omega_{k}-\omega_{l})^{2}\Mu-\Ku \right]^{-1}\Gu\left(\phiu_k,\phiu_l\right).
	\end{align}
	\label{eq:realmappingscalculation}
\end{subequations}
From Eq.~\eqref{eq:z_to_rs} it is then possible to establish the relationship between complex-valued and real-valued mappings as:
\begin{subequations}
	\begin{align}
		\hat{\au}_{kl} =&\, \frac{1}{2}\left( \Psiu^{(\mathrm{P})}_{kl} +\Psiu^{(\Nr)}_{kl}    \right), \\
		\hat{\bu}_{kl} =&\, \frac{1}{2\omega_{k}\omega_{l}}\left( \Psiu^{(\Nr)}_{kl} - \Psiu^{(\mathrm{P})}_{kl}  \right),\\
		\hat{\cu}_{kl}  =&\, \zerou,\\
		\hat{\alphau}_{kl} = &\,\zerou, \\
		\hat{\betau}_{kl} = &\,\zerou, \\
		\hat{\gammau}_{kl} = &\,  \frac{\omega_{l}+\omega_{k}}{\omega_l}\Psiu^{(\mathrm{P})}_{kl} + \frac{\omega_{l}-\omega_{k}}{\omega_{l}}\Psiu^{(\mathrm{N})}_{kl} ,
	\end{align}\label{eq:realmappings}
\end{subequations}
which corresponds to the results provided in~\cite{artDNF2020}. This procedure can be further extended for any expansion order.

\subsection{Real-Valued Single Mode Reduction}

Further insights into the relationship between complex and real-valued formulations can be derived by focusing on the case of a single mode reduction. In this case the polynomial mapping truncated  at second-order simply reads:
\begin{subequations}
	\begin{align}
		\Uu =&\, \phiu_{m} \rr_m + \hat{\au}_{mm}\rr_m^2 + \hat{\bu}_{mm}\sr_m^2, \\
		\Vu =&\, \phiu_{m} \sr_m + \hat{\gammau}_{mm}\rr_m\sr_m,
	\end{align}
\end{subequations}
where, from Eqs.~\eqref{eq:realmappings}, the relationships between real and complex mappings simplify to:
\begin{subequations}\label{eq:abcrealcomp}
\begin{align}
\hat{\au}_{mm} &= \frac{1}{2}(\Psiu^{(\Nr)}_{mm}+\Psiu^{(\mathrm{P})}_{mm}),\label{eq:abcrealcompa}\\
\hat{\bu}_{mm} &= \frac{1}{2\omega_{m}^{2}}(\Psiu^{(\Nr)}_{mm}-\Psiu^{(\mathrm{P})}_{mm}),\label{eq:abcrealcompb}\\
\hat{\gammau}_{mm} &= 2\,\Psiu^{(\mathrm{P})}_{mm}.\label{eq:abcrealcompc}
\end{align}
\end{subequations}

Recalling Eqs.~\eqref{eq:realmappingscalculation}, the evaluation of $\Psiu^{(\Nr)}_{mm}$ and $\Psiu^{(\mathrm{P})}_{mm}$ is performed through:
\begin{subequations}\label{eq:ZsZd}
	\begin{align}
		& \Psiu^{(\mathrm{P})}_{mm} =
		\left[ (2\,\omega_m)^{2}\Mu-\Ku \right]^{-1}\Gu\left(\phiu_m,\phiu_m\right),
		\label{eq:Zs} \\
		& \Psiu^{(\Nr)}_{mm} = 
		\left[-\Ku \right]^{-1}\Gu\left(\phiu_m,\phiu_m\right).\label{eq:Zd}
	\end{align}
\end{subequations}

The explicit expressions for $\Psiu^{(\Nr)}_{mm}$ and $\Psiu^{(\mathrm{P})}_{mm}$ allows for a non-intrusive implementation of the method in any FE software. Focusing on Eqs.~\eqref{eq:ZsZd}, they both consists in the solution of a linear system with unknown displacement-like vector $\Psiu_{mm}$ and imposed force vector $-\Gu\left(\phiu_m,\phiu_m\right)$. In particular, to solve Eq.~\eqref{eq:Zd}, a linear static analysis has to be performed, whereas to solve Eq.~\eqref{eq:Zs} a linear combination of stiffness and mass matrix has to be created in the software before solving the linear system. Regarding the computation of the force vector
$\Gu\left(\phiu_m,\phiu_m\right)$, it can also be non-intrusively calculated by means of the STEP method \cite{muravyov,Perez2014}. 
For example, by imposing a displacement along the mode $\phiu_m$, firstly with a positive, then with a negative modal amplitude, it is possible to extract the vectors $\Gu\left(\phiu_m,\phiu_m\right)$ and $\Hu(\phiu_{m},\phiu_{m},\phiu_{m})$ respectively.

The reduced dynamics reported in Eq.~\eqref{eq:sart_point_rv_sdof} in the case of single master mode $m$ is obtained by selecting the master coordinates as  $\rr_m={\zr}_{\alpha}+{\zr}_{\beta}$ and $\sr_m=\img\omega_{m}({\zr}_{\alpha}-{\zr}_{\beta})$. Application of Eq.~\eqref{eq:z_to_rs_map} to the reduced dynamics in Eq. \eqref{eq:sart_point_rv_sdof} then yields:
\begin{subequations}
	\begin{align}
		\dot{\rr}_m = &\, \sr_m, \\
		\dot{\sr}_m= &\, -\omega_{m}^{2}\rr_{m}  -\phiu_m^{\mathrm{T}}\left[ \Gu(\Psiu^{(\mathrm{P})}_{mm}+\Psiu^{(\Nr)}_{mm},\phiu_{m}) + \Hu(\phiu_{m},\phiu_{m},\phiu_{m}) \right] \rr_{m}^{3}+
\nonumber
\\
&\,-\phiu_m^{\mathrm{T}}\left[ \Gu(\Psiu^{(\Nr)}_{mm}-\Psiu^{(\mathrm{P})}_{mm},\phiu_{m}) \right]\frac{1}{\omega_{m}^{2}}\rr_{m}\sr_{m}^{2},
	\end{align}
\end{subequations}
where $\Psiu^{(1)}_{m}=\phiu_m$. Finally, using Eqs.~\eqref{eq:abcrealcompa}-\eqref{eq:abcrealcompb} to make $\hat{\au}_{mm}$ and $\hat{\bu}_{mm}$ appear, the third-order single-mode reduced dynamics finally reads:
\begin{equation}\label{eq:single_dof_real_touze}
\ddot{\rr}_m+\omega_{m}^{2}\rr_{m}+
(\hr_{m}+\Ar_{m}) \rr_{m}^{3}+
\Br_{m}\rr_{m}\dot{\rr}_{m}^{2}=0,
\end{equation}
where the introduced coefficients are given by:
\begin{subequations}
\begin{align}\label{eq:coeff_real_rdyn}
& \hr_{m}=\phiu_m^{\mathrm{T}}\Hu(\phiu_m,\phiu_m,\phiu_m),
\\
& \Ar_{m}=2\;\phiu_m^{\mathrm{T}}\Gu(\hat{\au}_{mm},\phiu_m),
\\
& \Br_{m}=2\;\phiu_m^{\mathrm{T}}\Gu(\hat{\bu}_{mm},\phiu_m).
\end{align}
\end{subequations}

The explicit form of these equations makes the computation of the reduced dynamics possible without the need of extrapolating the full tensors $\Gu$ and $\Hu$. In fact, the nonlinear force vectors $\Hu(\phiu_m,\phiu_m,\phiu_m)$, $\Gu(\hat{\au}_{mm},\phiu_m)$, and $\Gu(\hat{\bu}_{mm},\phiu_m)$ can be either computed in an intrusive manner by integration over each element in the FE code, or non-intrusively by means of the STEP method. 
This result is also in full accordance with previous derivations reported in~\cite{artDNF2020}. As detailed in~\cite{artDNF2020}, the STEP method allows extrapolating the required nonlinear force vectors in a non-intrusive manner from any finite element code by imposing a series of displacements on the structure and subsequently extracting the resulting forces.
If mappings and reduced dynamics coefficients are not truncated, then higher order terms appear in the reduced dynamics and the whole process can be pushed further. It is possible to show that also in the case of higher order calculations, the method can be implemented non-intrusively even though the number of calculations required in the STEP method would grow consistently with the order.

\section{Applications}
\label{sec:results}


The DNF method is here validated on three different MEMS structures with complex geometries,  with and without internal resonance. MEMS are known to be operated in near-vacuum packages and are thus subjected to very small damping values. They are also generally  actuated at resonance to fulfil technological requirements such as high sensitivity or large drive displacements. This in turn makes MEMS systems subjected to large displacements and geometric nonlinearities, hence making the DNF approach ideal for developing a predictive ROM strategy for design purposes. MEMS actuation is performed using either electrostatic, magnetic, or piezoelectric actuation. These types of actuation introduce nonlinearities that are not taken into account in the present work, their treatment via the DNF approach being postponed to further studies.

\begin{figure*}[t]
    \centering
    \includegraphics[width = .99\linewidth]{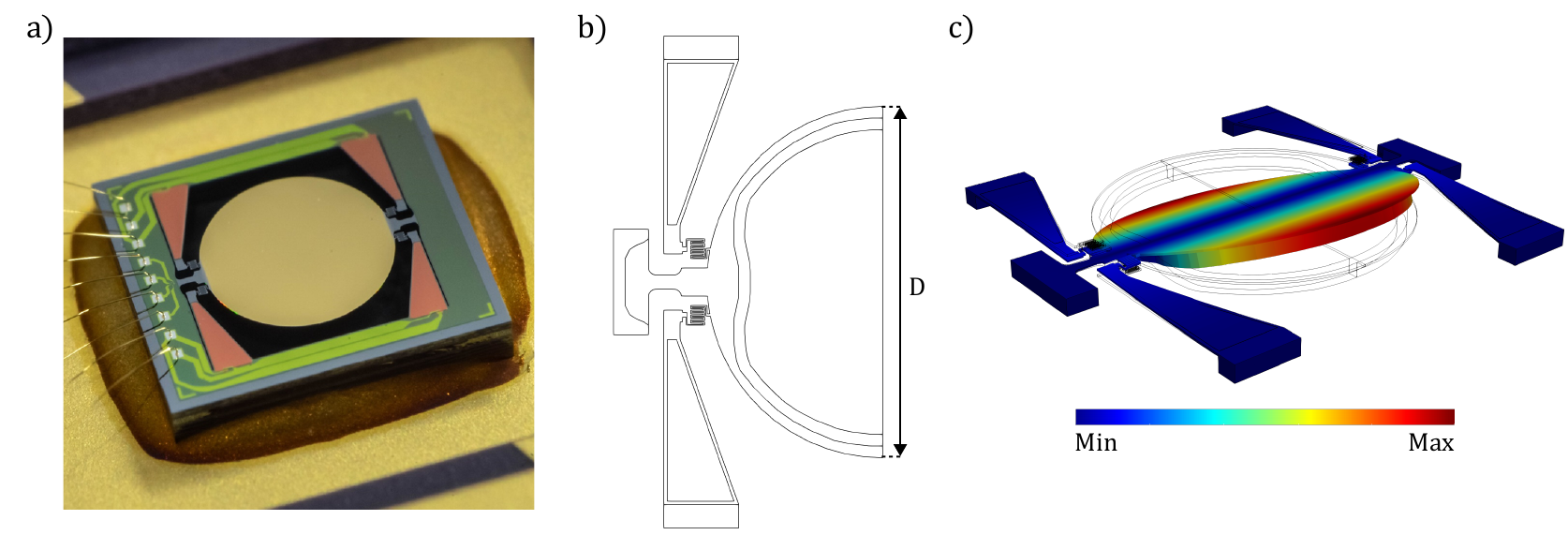}
    \caption{(a) Optical microscope picture of the MEMS micromirror. (b) Schematic representation of half-mirror geometry, D = 3000 $\mu$m. (c) Magnified displacement field associated to the first resonant mode of the structure.}\label{fig:m1_geom}
\end{figure*}

Let us denote as $\phiu_\Ir$ the actuated mode, hereafter referred to as drive mode. We assume that the excitation is provided thanks to a body force proportional to $\phiu_\Ir$ with a driving frequency in the vicinity of $\omega_\Ir$, the eigenfrequency of the drive mode. The losses are modelled with the assumption of mass-proportional damping, with coefficient equal to $\omega_\Ir/Q$, where $Q$ is the quality factor of the drive mode. The resulting model for MEMS structures in the physical space of the dofs of the FE mesh then reads:
\begin{align}\label{eq:fom_gobat}
    \mathbf{M} \ddot{\mathbf{U}}+\frac{\omega_\Ir}{Q}\mathbf{M}\dot{\mathbf{U}}+\mathbf{K U}+ 
    \mathbf{G (U, U)}+\mathbf{H (U, U, U)} = \kappa\mathbf{M}\boldsymbol{\phi}_\Ir \cos{(\omega t)},
\end{align}
with $\kappa$ a scalar load multiplier representing the intensity of the forcing. The damping model is pragmatic and corresponds to an established practice in microsystems. Indeed  
simplified models are often utilised to identify a unique damping parameter, {\em i.e.} a single $Q$ value, for the whole structure. This assumption is of current use in MEMS community, and it has shown to be reliable for a number of test cases \cite{GuerrieriContact,KennySingleQ}.

The ROMS are derived thanks to the second-order DNF method, where the nonlinear mapping is truncated at second-order while the reduced dynamics is truncated at third-order. The reduced model is then obtained by mapping Eq.~\eqref{eq:red_dyn} to real variables, following the guidelines given in Section~\ref{sec:real}. Modal forcing is then simply added at the right-hand side of the reduced dynamics. This assumption has already been used in~\cite{touze03-NNM,TOUZE:JSV:2006,TOUZE:CMAME:2008,artDNF2020}, and has been justified in~\cite{touze03-NNM,TOUZE:JSV:2006} by the fact that the variations of the time-dependent manifolds are of second-order as compared to other effects. A  mathematical justification of this assumption has also been  given in~\cite{VERASZTO}. Since a mass-proportional damping is selected and dissipation values are small, the simplest solution consists in adding directly the modal damping factors to the selected master modes. This simplification can be overcome by adding a more complex formulation of the damping, that takes into account the losses of all the slave modes to better represent the damping on the invariant manifold, following the general formula presented in~\cite{TOUZE:JSV:2006,artDNF2020}, and leading to a nonlinear form of the damping in the reduced dynamics. As shown in~\cite{artDNF2020}, this more complete model is particularly meaningful when using stiffness-proportional damping. In case of mass-proportional Rayleigh damping, the decay rates of the slave modes are rapidly negligible so that the assumption of modal master damping is sufficient. 

The proposed model is applied for the analysis of three MEMS structures showing different types of nonlinearity: a MEMS micromirror that undergoes large rotations, a beam resonator featuring 1:3 internal resonance, and a shallow arch resonator showing 1:2 internal resonance.

\subsection{MEMS Micromirror: Single Mode Reduction}\label{sec:single_dof_red}

The first device addressed is a MEMS micromirror developed by STMicroelectronics\textsuperscript{\textregistered}. The device in shown in Fig. \ref{fig:m1_geom}(a-b). The structure is composed of a circular reflective surface with a diameter of 3000 $\mu$m and is connected to the substrate through two torsional springs. Actuation of the device is performed with two pairs of lead-zirconate titanate PZT patches deposited on top of four trapezoidal beams. Patches are highlighted in orange in Fig. \ref{fig:m1_geom}(a). Silicon is anisotropic and a general discussion of its mechanical properties can be found {\em e.g.} in  \cite{hopcroft2010,opreni2021}.

The device is operated at resonance in the vicinity of its first mode with eigenfrequency $f_0=2266$\,Hz, corresponding to a rotation of the main mirror around the axis that passes through the two torsional springs. The displacement field of the mode is represented in Fig.~\ref{fig:m1_geom}(c). The first six eigenfrequencies of the structure are reported in Table \ref{tab:mirror_eig}, highlighting that no second-order resonance condition occurs between the first  and higher modes. 

\begin{table}[htb]
\centering
\caption{First six eigenfrequencies of the MEMS micromirror.}
\begin{tabular}{|c|c|c|}
\hline
Mode & Frequency {[}MHz{]} & Ratio $f_n/f_1$ \\ \hline
1    & 0.00227             & 1     \\ \hline
2    & 0.00726             & 3.198 \\ \hline
3    & 0.02523             & 11.11 \\ \hline
4    & 0.02527             & 11.13 \\ \hline
5    & 0.05605             & 24.69 \\ \hline
6    & 0.07335             & 32.31 \\ \hline
\end{tabular}
\label{tab:mirror_eig}
\end{table}
The ROM is obtained by imposing $\Psiu^{(1)}_{\alpha}=\Psiu^{(1)}_{\beta}=\phiu_{m}$ with $m=1$ (master mode), and by computing the mappings in Eq.~\eqref{eq:single_dof_map}. Noticeably, only two linear systems with symmetric matrices must be solved to obtain the second-order mapping and third-order reduced dynamics by exploiting the symmetries of the $\Lambdau$ operator. Furthermore, only the eigenvalue and the eigenfunction of the master mode is required. The operator  $\Gu(\Psiu^{(1)}_{\alpha},\Psiu^{(1)}_{\beta})$ needs to be computed only for the master mode and this step has moderate impact on the total running time of the analysis. Overall the technique is suitable for heavy parallelisation and vectorisation on modern processors. In the present contribution, all calculations have been realised thanks to a custom FE code developed by the authors without resorting to any STEP-like computation.\\

Reduced dynamics coefficients are obtained as detailed in Sec.~\ref{sec:ROM}.
The results of the proposed model are compared with full-order simulations of the device, which are performed with the Harmonic Balance Finite Element Method (HBFEM) with pseudo-arc length continuation. The solution of the HBFEM model is obtained with Fourier series expansion up to 7 to ensure model convergence. The geometry is discretised with quadratic (15 nodes) wedge elements and the HBFEM model has 15,341 nodes, corresponding to 690,345 Fourier coefficients. Also the ROM is solved with the Harmonic Balance (HB) method with an order 9 Fourier expansion.\\

\begin{figure}[htb]
    \centering
    \includegraphics[width = .99\linewidth]{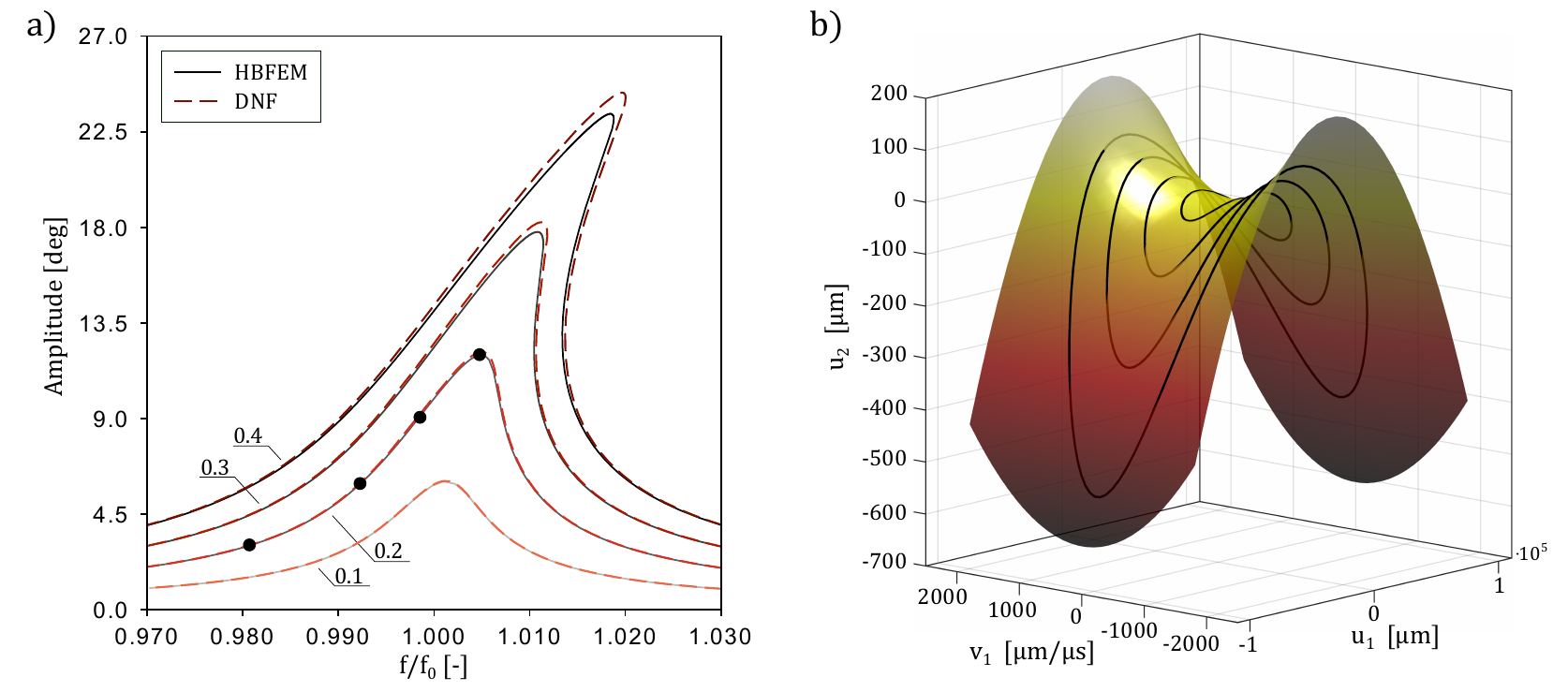}
    \caption{ (a) Maximum opening angle as a function of the frequency, comparison of HBFEM with the DNF method with single-mode reduction. Four increasing forcing amplitudes $\kappa$ [$\mu$m/$\mu$s$^{2}$]
are shown: 0.1, 0.2, 0.3 and 0.4. (b) Phase space representationn of the computed invariant manifold obtained with second-order DNF (colored surface), as compared to trajectories computed with the full-order model (black lines). The data are reported along the three-dimensional subspace $(\ur_1,\vr_1,\ur_2)$, corresponding to the modal displacement $\ur_1$ and velocity $\vr_1$ of the master mode to the modal displacement $\ur_2$ of the second mode of the structure, that is a slave mode. HBFEM trajectories corresponds to the points of the FRF in (a) highlighted by the bullets.}\label{fig:m1_results}
\end{figure}

The comparison between HBFEM and single-mode ROM solutions, obtained for four $\kappa$ values: 0.1, 0.2, 0.3, and 0.4 $\mu$m/$\mu$s$^{2}$, is reported in Fig. \ref{fig:m1_results}(a). The chart presents the maximum rotation angle reached by the mirror for a single steady-state oscillation cycle and for a given frequency value. The results show an excellent agreement between the HBFEM and ROM solution at any rotation angle, up to the curve with the highest $\kappa$ value. Furthermore, the invariant manifold defined by the expansion coincides with the trajectories of the HBFEM solutions, as highlighted by the phase-space representation of the system  in Fig. \ref{fig:m1_results}(b), where the velocity-dependence of the manifold is clearly observable. A slight departure between the ROM and the full-order solution appears at very large rotations, as a a consequence of the second-order DNF method used. Using higher-order developments would likely improve the match even at higher amplitudes.

The reduction to a single master mode offers the main benefit of embedding within the equation of motion of the master mode all contributions of non-resonant slave modes, hence accounting for the curvature of the manifold. This in turn yields impressive computational performance compared to the HBFEM simulations. The time required to integrate the ROM was less than 1 minute using a custom HB solver. The time required to compute the four HBFEM solutions reported in Fig.~\ref{fig:m1_results} was equal to 1 week, hence highlighting how the present method offers outstanding results with little computational resources. Stability is not reported in the present analysis since the HBFEM solution has not this feature implemented. Considerations on the stability are postponed to Section~\ref{subsec:stability}.


\begin{figure}[htb]
    \centering
    \includegraphics[width = .99\linewidth]{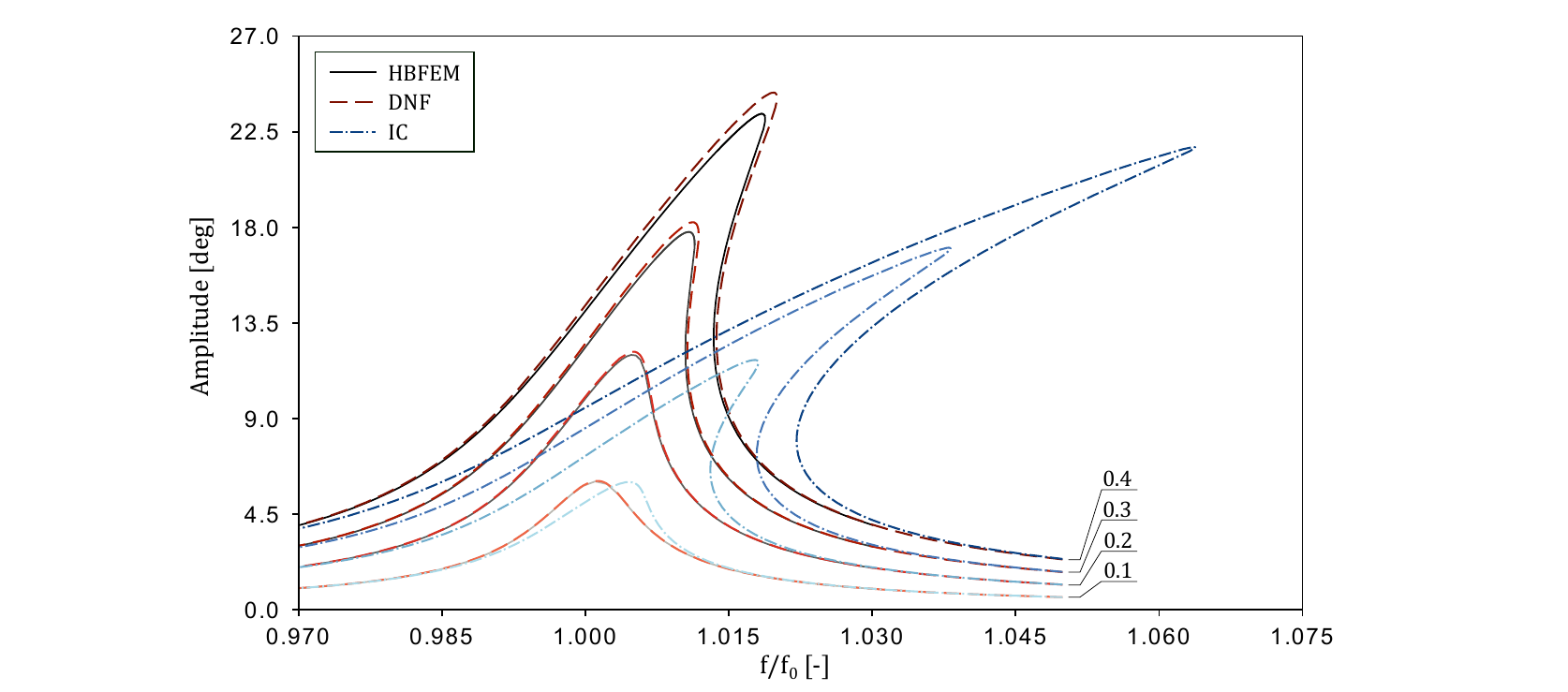}
    \caption{Comparison between HBFEM solution, DNF method, and implicit condensation method (IC). $\kappa$ [$\mu$m/$\mu$s$^{2}$] values associated to each curve are highlighted by the tags in the chart.}\label{fig:m1_results_SC}
\end{figure}

In order to better highlight the quality of the results obtained with the DNF method, we compare with another reduction strategy using implicit condensation (IC) as proposed in~\cite{Hollkamp2008,FRANGI2019}.  As shown for example in~\cite{YichangVib,NicolaidouIceKE} with the case of a cantilever beam, IC method fails at reproducing correctly the nonlinear dynamics of structures when inertia nonlinearities are importantly excited. In particular, the IC method condenses statically the contribution of slave modes on the master mode trajectory, hence neglecting velocity dependent terms~\cite{YichangICE}. The comparison between the DNF and IC methods is reported in Fig.~\ref{fig:m1_results_SC}, highlighting the benefits of adopting the proposed method for modelling systems subjected to large rotations and, more in general, to large transformations. The IC reduction method overestimates the hardening behaviour, and cannot be used in such context of large rotations, whereas the DNF approach is giving uniformly valid solutions without the need of any extra assumption thanks to the invariance property. 

Further insight 
can be gained by looking more closely at the individual values of the coefficients of the reduced dynamics with single mode, Eq.~\eqref{eq:single_dof_real_touze}. The MEMS micromirror under investigation does not show mode interactions even at large rotation amplitudes, yet methods based on static condensation show difficulties in modelling the nonlinear response of the device. For this device, the physical cubic coefficient $\hr_{m}$ in Eq.~\eqref{eq:single_dof_real_touze} is equal to 6.479$\cdot$10$^{-10}$ $\mu$N/$\mu$m$^3$. 
On the other hand, the first correction term computed by the normal transform $\Ar_{m}$ equals -6.479$\cdot$10$^{-10}$ $\mu$N/$\mu$m$^3$, hence their sum cancels out on the leading cubic term $r_m^3$ in Eq.~\eqref{eq:single_dof_real_touze}. This is a remarkable result from a physical standpoint since it can be stated that elastic nonlinearities associated to the excited mode are not the main source of nonlinear response of the structure. On the other hand, $\Br_{m}$ is equal to 4.53$\cdot$10$^{-12}$ $\mu$N$\mu$s$^2$/$\mu$m$^3$. This last term is a velocity dependent nonlinearity which accounts for the change in inertia of the system during motion. Therefore, the nonlinear response of the structure is conveyed by the change in configuration of the structure rather than elastic nonlinearities. Interestingly, similar observations on the particular values of the coefficients have been reported in the case of a cantilever beam.\\

\begin{figure*}[htb]
    \centering
    \includegraphics[width = .99\linewidth]{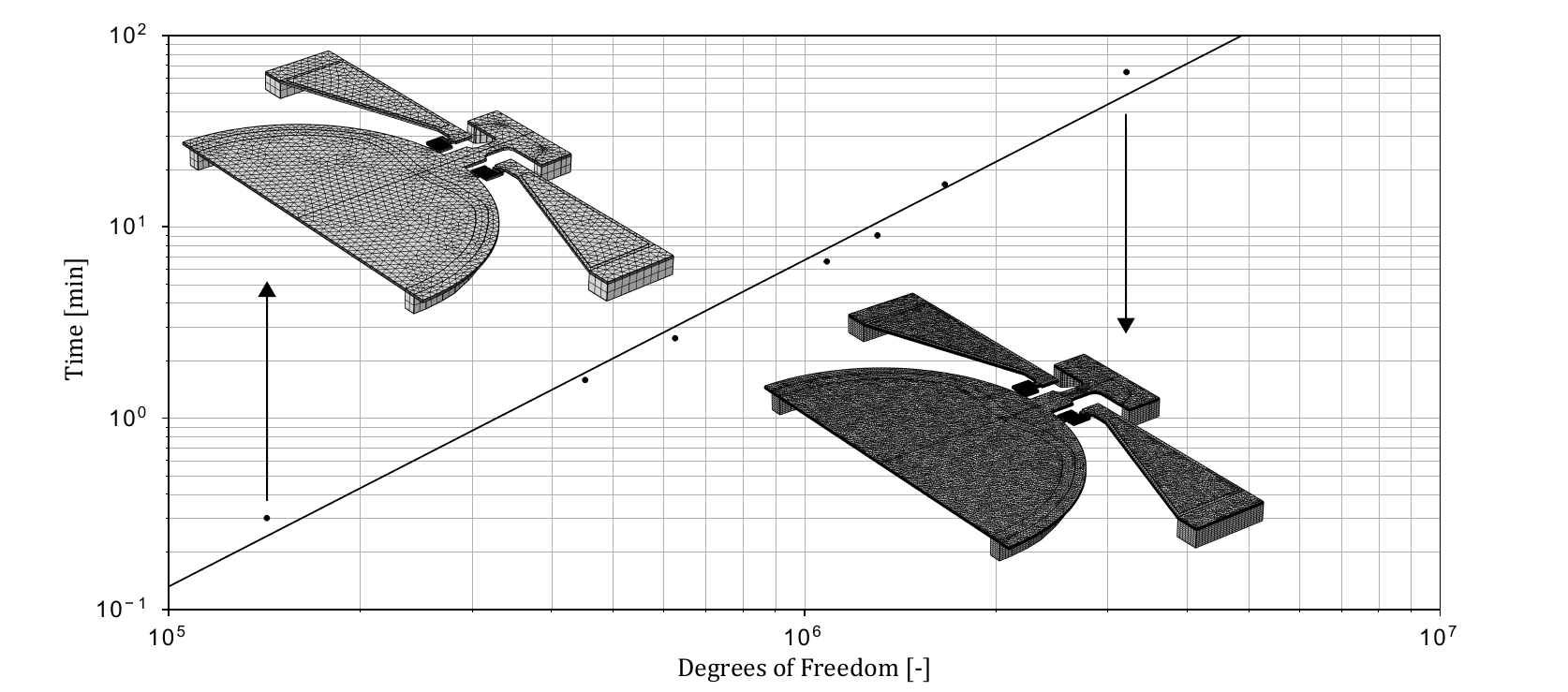}
    \caption{Computational time as a function of the model size for different mesh sizes. Results are obtained for an intrusive version of the method. The computation was run on a desktop workstation with an AMD Ryzen\texttrademark$\,$5 1600 Six-Core Processor 3.20 GHz and 64 GB RAM.}\label{fig:performance}
\end{figure*}

An example of the performance of the proposed method is reported in Fig. \ref{fig:performance}, where the CPU times required to obtain the mappings and reduced dynamics coefficients of the MEMS micromirror for different mesh refinements are reported. The computation is performed on a desktop workstation with an AMD Ryzen\texttrademark$\,$5 1600 Six-Core Processor 3.20 GHz and 64 GB RAM. The computational times highlight how the proposed technique is scalable to large models even with moderate computational resources. The integration of the ROM with the HB is not reported since it is has a negligible cost.

\subsection{Internally resonant MEMS resonators}

In this Section, two different internal resonance scenarios are investigated on MEMS resonators having a beam-like and an arch-like structure, respectively. 
In the first case, a 1:3 internal resonance is investigated, where the second-order DNF method does not need extra assumptions. Indeed the third-order monomials are not affected by the second order mappings and they all remain in the reduced dynamics, making the excitation of the 1:3 resonance possible. 
This advantage of the second-order DNF as compared to higher orders has been already underlined in~\cite{artDNF2020} and will be further commented herein. The second case is a 1:2 resonance which needs extra developments of the DNF approach as discussed in Section~\ref{subsec:IRmulti}, with the reduced dynamics given in Appendix~\ref{sec:second_order_ir_coeff}. In both cases, reduction to the two internally resonant modes is able to catch the complex nonlinear phenomena and retrieve the frequency-response functions of the structures. The dynamics is reduced to a four-dimensional invariant manifold.


\begin{figure*}[htb]
    \centering
    \includegraphics[width = .99\linewidth]{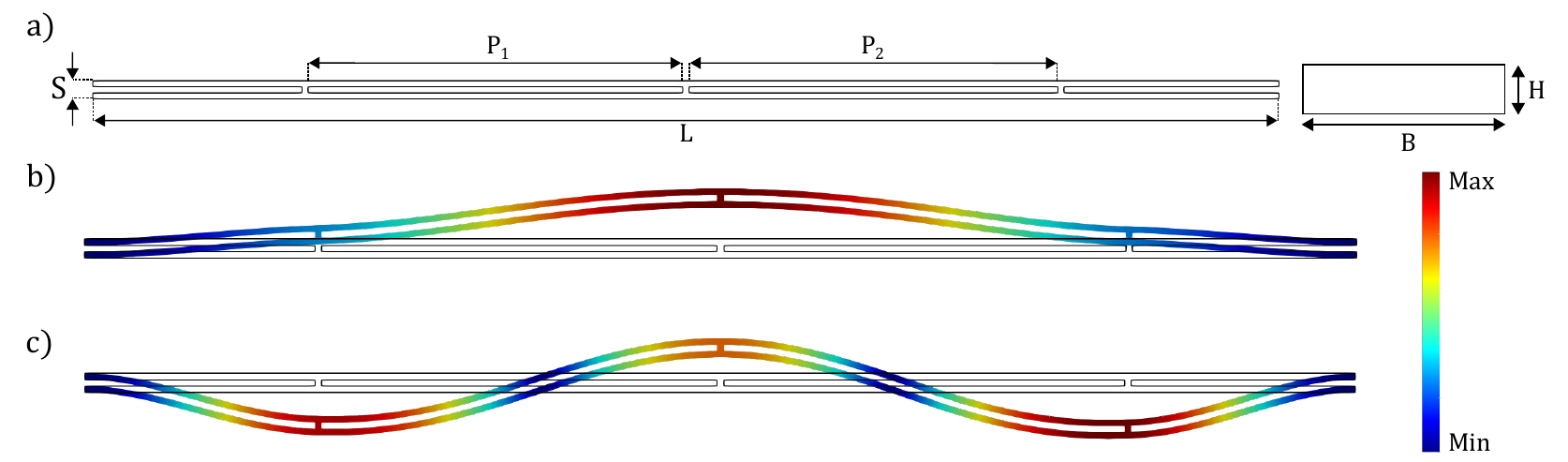}
    \caption{(a) Beam view along the $\eu_1$-$\eu_2$ plane on the left and beam transverse section geometry on the right. The geometrical dimensions of the beam are B = 12 $\mu$m, H = 5 $\mu$m, L = 1000 $\mu$m, P$_1$ = 311 $\mu$m, P$_2$ = 316 $\mu$m, and S = 15 $\mu$m. S size includes the two beams thicknesses as well as the interspace, each of them of 5 $\mu$m. Modes that satisfy the 1:3 internal resonance condition are mode 1 (b) and mode 4 (c).  The colormap reports the magnitude of the displacement field. The deformation of the structure is amplified to improve the visualisation of the field.}\label{fig:beam_geom}
\end{figure*}

\subsubsection{MEMS beam resonator: 1:3 internal resonance}
In this section the analysis of a 1:3 internally resonant beam is reported. The device is shown in Fig. \ref{fig:beam_geom}(a). The structure is made of two beams with a length of 1000 $\mu$m and a thickness of 5 $\mu$m. The two parts are connected in three points. The structure is not symmetric with respect to the center since the external connection points are located at different 
distances (see the values of P$_1$ and P$_2$ reported in the caption of Fig.~\ref{fig:beam_geom}). The material is polycrystalline silicon which is modelled as isotropic with a Young's modulus of 167 GPa and a Poisson ratio of 0.22.

The first six eigenfrequencies of the structure are reported in Table~\ref{tab:beam_eig}, which highlights an almost perfect 1:3 ratio between mode 1 and mode 4, whose displacement field is reported in Fig.~\ref{fig:beam_geom}(b-c).

\begin{table}[htb]
\centering
\caption{First six eigenfrequencies of the 1:3 internally resonant beam.}
\begin{tabular}{|c|c|c|}
\hline
Mode & Frequency {[}MHz{]} & Ratio \\ \hline
1    & 0.08100             & 1     \\ \hline
2    & 0.10411             & 1.285 \\ \hline
3    & 0.19690             & 2.431 \\ \hline
4    & 0.24306             & 3.001 \\ \hline
5    & 0.28683             & 3.541 \\ \hline
6    & 0.41043             & 5.067 \\ \hline
\end{tabular}
\label{tab:beam_eig}
\end{table}

The ROM is built by taking the pairs $(\zr_p,\zr_{p+\Nr})$ and $(\zr_{q},\zr_{q+\Nr})$ such that $\Psiu^{(1)}_{p}=\Psiu^{(1)}_{p+\Nr}=\phiu_{1}$ and $\Psiu^{(1)}_{q}=\Psiu^{(1)}_{q+\Nr}=\phiu_{4}$. All remaining terms of $\zu$ are set to zero. This implies a moderate  increase of the computational burden as compared to the previous case of single-mode reduction, since only mappings and reduced dynamics coefficients that multiply a non-zero entry of $\zu$ need to be computed. 

\begin{figure*}[htb]
    \centering
    \includegraphics[width = .99\linewidth]{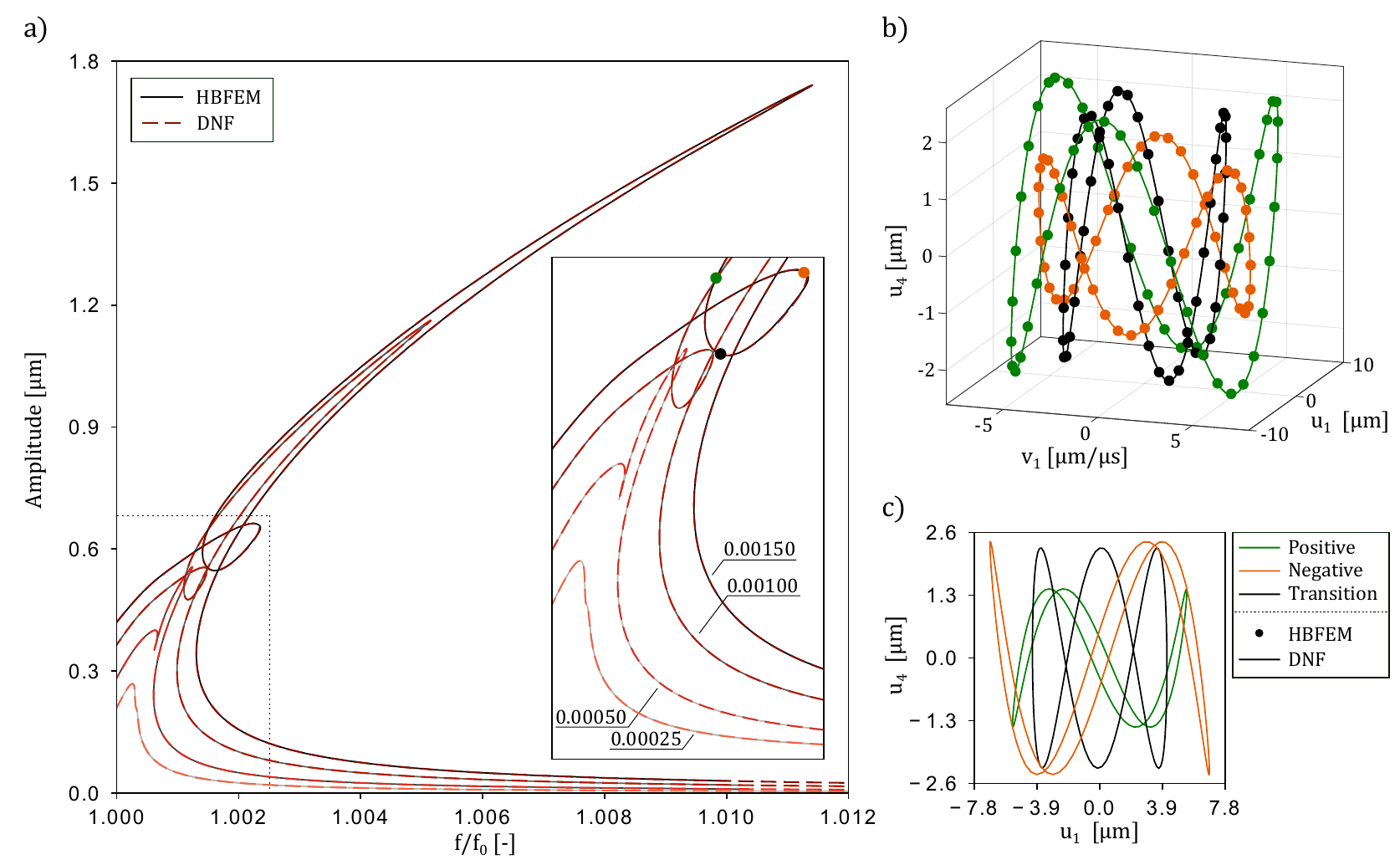}
    \caption{
    (a) Frequency-response function (FRF) of the beam MEMS resonator featuring 1:3 resonance. Comparison between HBFEM and ROM for increasing values of $\kappa$  [$\mu$m/$\mu$s$^{2}$], from 2.5$\cdot$10$^{-4}$ to 1.5$\cdot$10$^{-3}$. (b) Representation of three trajectories in phase space along the $\ur_1$-$\vr_1$-$\ur_4$ coordinates, comparison between the HBFEM solution (bullets) and the ROM solution (continuous line). Trajectories are sampled from the regions highlighted by the bullets in  (a), corresponding to the largest forcing amplitude. (c) Projection of the sampled trajectories along the $\ur_1$-$\ur_4$ plane.
    }\label{fig:beam_results}
\end{figure*}

As for the MEMS micromirror, the ROM is validated with full order HBFEM solutions with a Fourier expansion coefficient of order 7. The geometry is discretised with quadratic (15 nodes) wedge elements  and the resulting HBFEM model is made of 12,906 dofs, corresponding to 193,590 nodal unknowns of the HBFEM problem. The quality factor of the model is set to $3000$. Analyses are performed for $\kappa$ values equal to 0.00025, 0.0005, 0.001, and 0.0015 $\mu$m/$\mu$s$^{2}$. The ROM is solved again with the HB method with a Fourier expansion order equal to 9. Stability is not reported.\\

The comparison between HBFEM and ROM solutions is presented in Fig. \ref{fig:beam_results}(a). The accuracy of the reduced model is remarkable for any $\kappa$ value. The strong nonlinear interaction is put in evidence by the loop appearing in the frequency response function, even at small $\kappa$ values. 
A further comparison is shown in Fig.~\ref{fig:beam_results}(b), where three trajectories are represented in the space ($\ur_1$, $\vr_1$, $\ur_4$). They have been obtained for the largest $\kappa$ value selected, and have been chosen in the vicinity of the 1:3 resonance loop, as marked by the coloured points in Fig.~\ref{fig:beam_results}(a). While the orange point is before the 1:3 resonance, the black point is exactly at resonance, and the green point after. Comparison between full and reduced solutions shows a very good match for the trajectories. Interestingly, the trajectories show a cubic shape with either positive or negative linear term, and they can be related to the two families of periodic orbits (backbone curves) arising in the 1:3 internal resonance. 
As a matter of fact, the solutions keep these shapes all along the FRF, before and after the 1:3 resonance. At exact resonance when the forcing frequency is equal to $\omega_4/3$, a transition form is obtained and reported in black in Fig.~\ref{fig:beam_results}(b-c). The gain in computational time is again impressive, with a factor of 3000 between the ROM and the full-order HBFEM solution since a single FRF of the device requires approximately two days, while the solution of the reduced model requires less than a minute.


\subsubsection{MEMS arch resonator with 1:2 internal resonance}
\label{subsec:archRES12}

\begin{figure*}[htb]
    \centering
    \includegraphics[width = .99\linewidth]{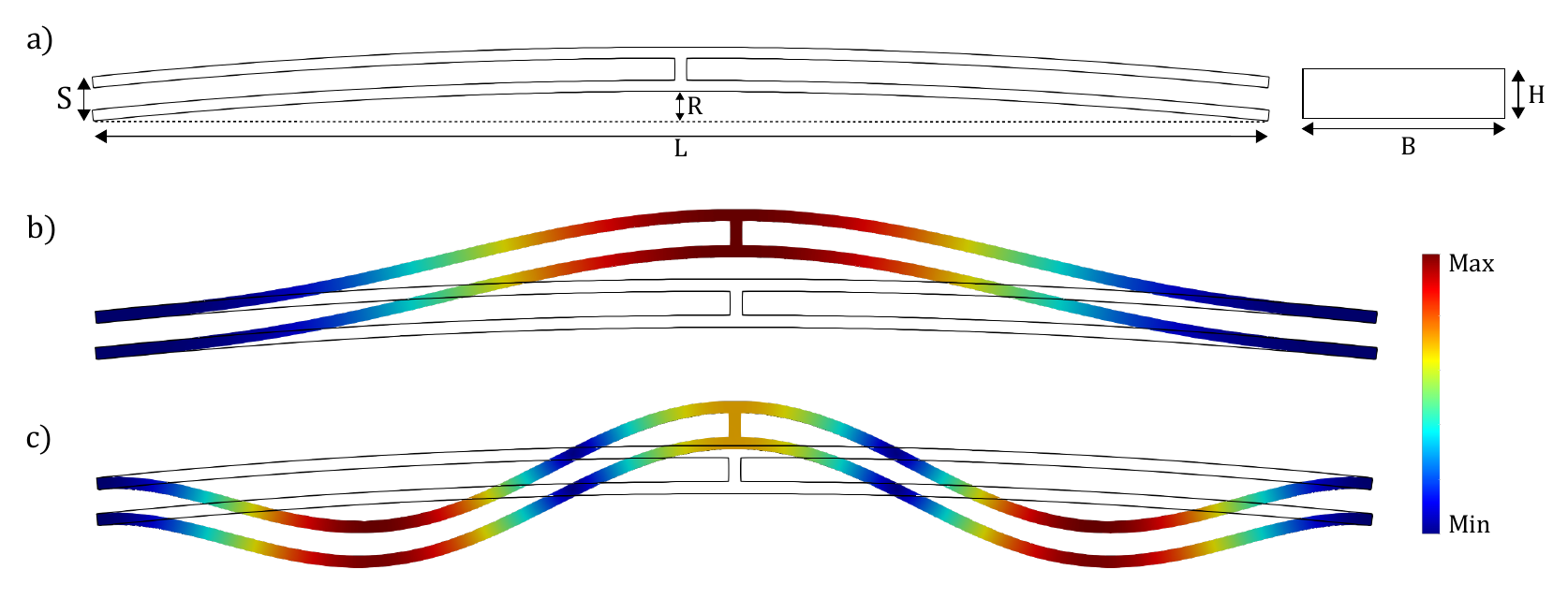}
    \caption{Arch view along the $\eu_1$-$\eu_2$ plane on the left and arch beam section geometry on the right (a). The geometrical dimensions of the arch are B = 20 $\mu$m, H = 5 $\mu$m, L = 530 $\mu$m, R = 13.4 $\mu$m, and S = 20 $\mu$m. S size includes the beams thickness. Modes that satisfy the 1:2 internal resonance condition are mode 1 (b) and mode 6 (c).  The colormap reports the magnitude of the displacement field. The deformation of the structure is amplified to improve the visualisation of the field.}\label{fig:arc_geom}
\end{figure*}

The structure under investigation is the MEMS arch resonator depicted in Fig.~\ref{fig:arc_geom}(a). The arch has a length of 530 $\mu$m and it is made by two arched beams with a thickness of 5 $\mu$m connected at their midpoint. The structure is made in polycrystalline which is modelled as isotropic with a Young's modulus of 167 GPa and a Poisson's ratio of 0.22.

\begin{table}[htb]
\centering
\caption{First six eigenfrequencies of the 1:2 internally resonant arch.}
\begin{tabular}{|c|c|c|}
\hline
Mode & Frequency {[}MHz{]} & Ratio \\ \hline
1    & 0.43416             & 1     \\ \hline
2    & 0.52597             & 1.211 \\ \hline
3    & 0.60391             & 1.391 \\ \hline
4    & 0.66759             & 1.537 \\ \hline
5    & 0.75695             & 1.743 \\ \hline
6    & 0.86367             & 1.989 \\ \hline
\end{tabular}
\label{tab:arc_eig}
\end{table}

As reported in Table \ref{tab:arc_eig}, the structure shows an almost perfect 1:2 internal resonance between  mode 1 and mode 6. The displacement field associated to the two modes is reported in Fig. \ref{fig:arc_geom}(b-c). Due to the 1:2 internal resonance, the reduced model is built by taking the $(\zr_p,\zr_{p+\Nr})$ and $(\zr_{q},\zr_{q+\Nr})$ such that $\Psiu^{(1)}_{p}=\Psiu^{(1)}_{p+\Nr}=\phiu_{1}$ and $\Psiu^{(1)}_{q}=\Psiu^{(1)}_{q+\Nr}=\phiu_{6}$, hence yielding a two oscillators model as the one reported in Sec.~\ref{sec:ROM}, while the reduced dynamical equations are given in Appendix~\ref{sec:second_order_ir_coeff}. Damping and forcing is added following the general guidelines given at the beginning of Section~\ref{sec:results}.

\begin{figure*}[htb]
    \centering
    \includegraphics[width = .99\linewidth]{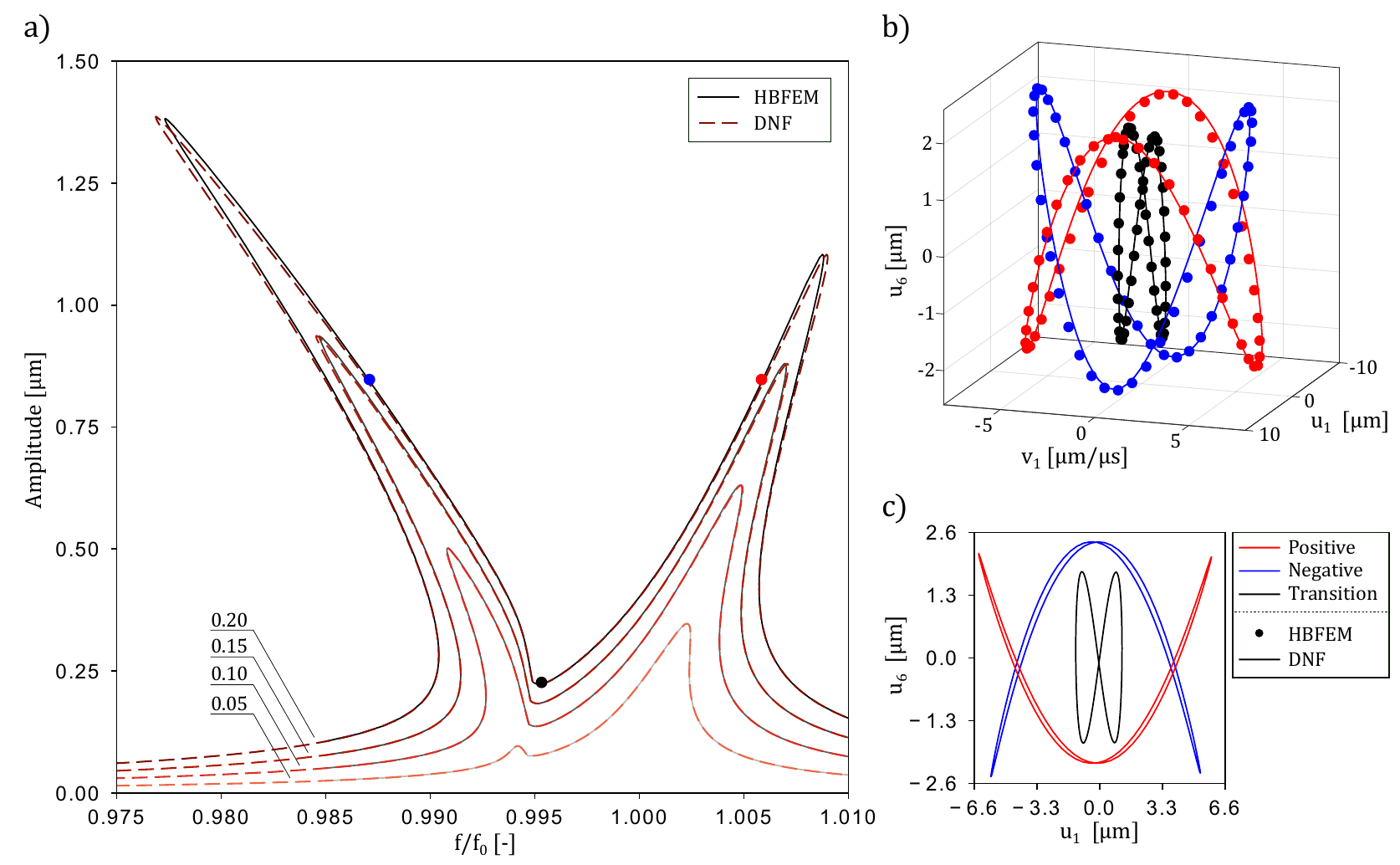}
    \caption{(a) Frequency-response function (FRF) of the arch MEMS resonator featuring 1:2 resonance. Comparison between HBFEM and ROM for increasing values of $\kappa$  [$\mu$m/$\mu$s$^{2}$], from 0.05 to 0.2. (b) Representation of three trajectories in phase space along the $\ur_1$-$\vr_1$-$\ur_6$ coordinates, comparison between the HBFEM solution (bullets) and the ROM solution (continuous line). Trajectories are sampled from the regions highlighted by the bullets in  (a), corresponding to the largest forcing amplitude. (c) Projection of the sampled trajectories along the $\ur_1$-$\ur_6$ plane.
}\label{fig:arch_results}
\end{figure*}

The results provided by the ROM are compared with the full-order HBFEM solution of the MEMS device. HBFEM Fourier expansion order is taken up to order 9 to ensure convergence of the method. The geometry is discretised with quadratic (15 nodes) wedge elements and and the resulting HBFEM model is made of 5,913 dofs, hence the total number of nodal unknowns is equal to 112,347. The ROM is solved with the HB method with order 9. Curves are computed assuming a quality factor $Q$ equal to 500 and analyses are performed for four $\kappa$ values: 0.05, 0.1, 0.15, and 0.2 $\mu$m/$\mu$s$^{2}$.
The comparison between HBFEM and ROM solution is reported in Fig. \ref{fig:arch_results}a. The chart presents the maximum displacement reached by the device during an oscillatory cycle for each frequency value. The data highlight a perfect agreement between HBFEM and ROM solutions for each $\kappa$ value. Only at the highest amplitude a small discrepancy is observed.\\

A set of trajectories of the solution is reported in phase space in Fig. \ref{fig:arch_results}(b-c). The trajectories are selected from the points in the FRF shown by bullets in Fig.~\ref{fig:arch_results}(a),
obtained for the largest forcing amplitude $\kappa$=0.2 $\mu$m/$\mu$s$^{2}$. Comparison between full and reduced-order models on the trajectories shows again an excellent agreement. Interestingly, the shape of the trajectories in phase space can be related to the properties of the underlying backbones of the 1:2 system, as investigated in~\cite{Giorgio12}. Indeed, two families of periodic orbits exist in that case, and are named as parabolic $p^{+}$ and $p^{-}$ modes depending on the sign of the curvature. The left-hand part of the FRF follows the $p^{-}$ backbone such that the shape of the trajectories reproduce the negative parabola in the correct axis, as shown in Fig.~\ref{fig:arch_results}(b-c). The right-hand part corresponds to $p^{+}$ solutions. 
When the forcing frequency equals $\omega_6/2$, a transition form is obtained.


\subsubsection{Stability Analysis}\label{subsec:stability}

\begin{figure*}[htb]
    \centering
    \includegraphics[width = .99\linewidth]{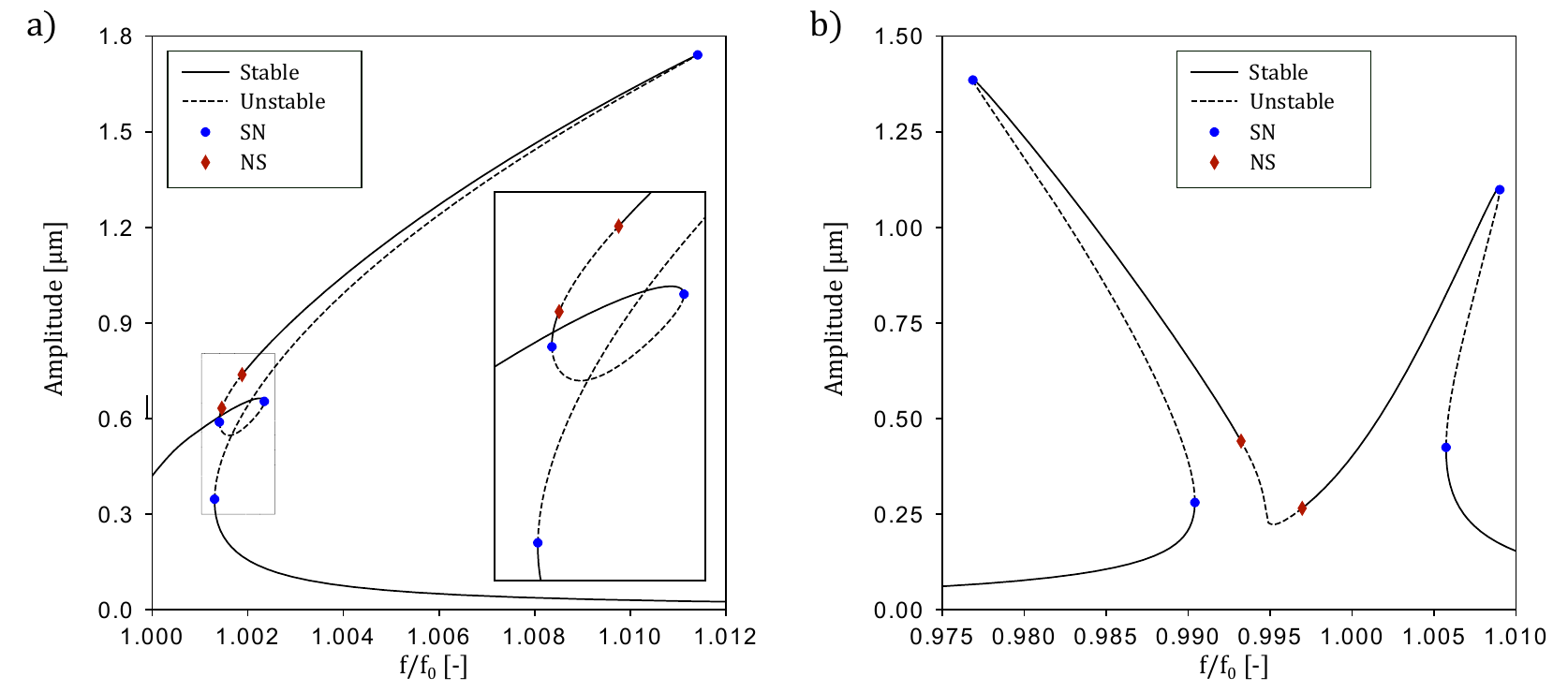}
    \caption{Stability analysis of the internally resonant MEMS resonators. Chart (a) details the stability analysis of the 1:3 internally resonant beam. Data are shown for $\kappa=0.0015\,\mu$m/$\mu$s$^{2}$. Chart (b) reports the same analysis on the 1:2 internally resonant arch resonator. Data are shown for $\kappa=0.2\,\mu$m/$\mu$s$^{2}$. SN: saddle-node bifurcation point, NS: Neimarck-Sacker bifurcation point.}\label{fig:stability}
\end{figure*}

A complete characterisation of the nonlinear dynamic response associated to the structures reported in this Section requires evaluating the stability of the response and the detection of bifurcation points. However, 
the direct implementation of stability in large scale numerical procedures such as the HBFEM is challenging and asks for special attention. In particular, in our implementation of HBFEM, the stability is not computed yet.
On the other hand, a number of available open-source continuation codes can be applied 
to small-scale systems like the ROMs previously discussed. 
Herein, stability analyses are performed using the ManLab continuation package \cite{Guillot2019} on the reduced dynamics, in order to give insight to the results presented in previous sections.\\
The 1:3 internally resonant beam is addressed in Fig. \ref{fig:stability}(a) for a load multiplier $\kappa$ value equal to 0.0015 $\mu$m/$\mu$s$^{2}$. Starting from lower frequency values and moving towards larger values, the data highlight the presence of an unstable region enclosed between two saddle-node (SN) bifurcations. Afterwards, the response becomes stable again, until the system encounters another unstable region enclosed between two Neimark-Sacker (NS) bifurcations, hence suggesting the onset of quasi-periodic response of the system. The system becomes then stable again until a new unstable region enclosed between another set of SN bifurcations is found. Finally, the system becomes stable again and it retains the usual quasi-static response at higher frequencies. Overall, stability and bifurcation analysis of the system highlights important features that are compulsory for correct estimation of the structure response. Furthermore, this result was obtained with a negligible computational cost.

The stability analysis of the 1:2 internally resonant arch is reported in Fig. \ref{fig:stability}(b) for $\kappa$ equal to 0.0015 $\mu$m/$\mu$s$^{2}$. Starting from low frequencies, the chart shows the same sequence of stable and unstable regions observed for the 1:3 internally resonant beam. Indeed, starting from lower frequencies and moving upward, the system starts as stable. Then, an unstable region enclosed between two SN bifurcations is found. Afterwards, the system becomes stable again. An unstable region enclosed between two NS bifurcations is then observed. Stability is then recovered, only to be broken once again by another unstable region enclosed between two other pairs of SN bifurcations. Therefore, both 1:2 and 1:3 internal resonance examples show two pairs of SN bifurcations and one pair of NS bifurcations, all enclosing an unstable region of the system.\\

Overall, the proposed reduction technique does not only allow for an efficient estimation of the frequency response function of the device, but it also enables refined stability analysis that would be computationally expensive for large scale models.

\section{Conclusions}\label{sec:conclusions}

In this paper, a direct normal form approach for model order reduction of the discretised equation of linear momentum for mechanical systems subjected to geometric nonlinearities is proposed. The direct computation has been rewritten with a state-space formulation as starting point, resulting in a symmetric and complex formulation of all the main quantities (mappings and reduced dynamics). The major outcome is the generalisation of the method by a formulation of the homological equations that keeps trace of the mechanical setting, thus facilitating further developments of the method. The complete rewriting also allows to formulating a comparison with the real-valued approach of the DNF originally proposed in~\cite{artDNF2020}, underlining the equivalence of the different methodologies. A special emphasis has also been put on the treatment of second-order internal resonance, a feature that had not been developed in all previous studies using this approach.\\

For illustration purposes, only second-order DNF has been used for numerical examples in this contribution, where the nonlinear mapping is truncated at second-order while third-order dynamics is computed. The method proves efficient for developing ROMs of large-scale finite element models, showing excellent computational performance and accuracy. Furthermore, the method does not rely on the slow-fast assumption between master and slave modes, since the velocity is directly accounted by the parametrisation procedure. This is a major achievement since it offers a uniformly valid and simulation-free method that could be blindly applied to any structure without the need of extra assumptions~\cite{Vizzaccaro:NNMvsMD,YichangICE,YichangVib,YichangNODYCON}, for the same computational cost of other methods. 
Normal form theory makes the distinction between resonant and non-resonant couplings. Non-resonant couplings are automatically embedded in the master mode thanks to the curvature of the invariant manifold, so that there is no need of computing extra vectors (like MDs or dual modes) to achieve convergence of the ROM. Resonant couplings create strong energy exchange between the internally resonant modes and drastically modify the nonlinear dynamics of the system, such that modes in internal resonance must be appended to the ROM dimension.\\

A further advantage of the method is that it does not need to rely on the knowledge of all nonlinearities of the full order model, since it needs to evaluate system nonlinearities only for the mapping vectors associated to the master modes, hence further enhancing the scalability of the method. Furthermore, it has been underlined that the explicit results of the mappings presented in this contribution open the door to a non-intrusive coding of the method,  hence making the DNF easily integrable within finite element commercial software.\\

The efficiency of the method is proved by studying  nonlinear structures with strong reduction to 1 or 2 master modes. Excellent accuracy of the results have been demonstrated on three different examples of increasing dynamical complexity.  Higher order developments may be required for studying more complex structures, for instance structures that show an initial softening behaviour followed by hardening response, a topic that is left for future works, together with the treatment of damping and forcing within the procedure.


\section*{Acknowledgments}
The authors are grateful to Giorgio Gobat for the design of the 1:2 and 1:3 internally resonant structures and to STMicroelectronics\textsuperscript{\textregistered} for the model of the MEMS Micromirror.

\section*{Funding}
The work received no additional funding.

\section*{Conflict of interest} 
The authors declare that they have no conflict of interest.

\section*{Data availability statement}
The codes written to run most of
the simulations presented in this paper can be available upon
simple request to the authors.

\bibliographystyle{unsrt}
\bibliography{biblio}

\appendix

\section{Notation}
\label{sec:notation}

Compact expressions are introduced throughout the paper in order to improve readability.
For polynomial representation of nonlinear terms, functional and indicial formulations are linked with the following relationships for $\Tu_{s},\Tu_{kl},\Tu_{klm}$  vectors of coefficients:
%
\begin{subequations}\label{eq:tensors_general}
	\begin{align}
		\Tu\xu =& \sum_{s}\Tu_{s}\xr_s, \\
		\Tu(\xu,\xu) =& \sum_{k,l}\Tu_{kl}\xr_k\xr_l, \\
		\Tu(\xu,\xu,\xu) =& \sum_{k,l,m}\Tu_{klm}\xr_k\xr_l\xr_m,
	\end{align}
\end{subequations}
%
The distributive property can be expressed as:
\begin{equation}\label{eq:distributive}
	\Tu(\xu,\xu+\yu,\xu) = \Tu(\xu,\xu,\xu) + \Tu(\xu,\yu,\xu).
\end{equation}
Moreover, since $\Tu_{s},\Tu_{kl},\Tu_{klm}$ are constant vectors of coefficients used for polynomial representation, the following definitions hold for time derivatives:
\begin{subequations}\label{eq:time_derivative}
	\begin{align}
		\frac{d}{dt}\left({\Tu}\xu\right) =&\, \Tu\dot{\xu}, \\
		\frac{d}{dt}\left({\Tu}(\xu,\xu)\right) =&\, \Tu(\dot{\xu},\xu) + \Tu({\xu},\dot{\xu}), \\
		\frac{d}{dt}\left({\Tu}(\xu,\xu,\xu)\right) =&\, \Tu(\dot{\xu},\xu,\xu) + \Tu({\xu},\dot{\xu},\xu) + \Tu({\xu},{\xu},\dot{\xu}),
	\end{align}
\end{subequations}
%

\section{Explicit expressions for geometric nonlinear terms}
\label{sec:nl_expr}

We develop analytical expressions for quadratic and cubic nonlinearities in mechanical systems subjected to large transformations (geometric nonlinearities). 
The weak form of the linear momentum conservation equation mapped to reference configuration is expressed as:
\begin{equation}\label{eq:lm_weak_material}
    \int_{\Omega_0} \rho_0\,\ddot{\uuu}\cdot\www\,d\Omega_0 + \int_{\Omega_0} \Su:(\Fu^{\mathrm{T}}\nabla\www)\,d\Omega_0 = 0, \qquad \forall\,\www\in\cC(\zerou),
\end{equation}
with $\Omega_0$ the domain in reference configuration, $\rho_0$ the reference density, $\uuu$ the displacement field, $\www$ the test field defined over the space $\mathcal{C}(\mathbf{0})$ of functions that vanish on the portion of boundary where Dirichlet boundary conditions are prescribed. $\Su$ is the second Piola-Kirchhof stress tensor, $\Fu$ is the deformation gradient, {\em i.e.} $\Iu+\nabla\uuu$ with $\nabla(\cdot)$ denoting material gradient. The first integral in Eq.~\eqref{eq:lm_weak_material} is the kinetic power while the second term represents the power of internal forces. We introduce the Green-Lagrange strain tensor as $\eu=(\Fu^{\mathrm{T}}\Fu-\Iu)/2$. 
Using a Saint-Venant Kirchhoff constitutive model, {\em i.e.} $\Su=\cA:\eu$ with $\cA$ fourth order elasticity tensor, Eq.~\eqref{eq:lm_weak_material} is expressed as:
\begin{equation}\label{eq:lm_weak_cm}
    \int_{\Omega_0} \rho_0\,\ddot{\uuu}\cdot\www\,d\Omega_0 + \int_{\Omega_0} \delta\eu:\cA:\eu\,d\Omega_0 = 0, \qquad \forall\,\www\in\cC(\zerou).
\end{equation}
We notice that since $\mathbf{e}$ is nonlinear with respect to the displacement, Eq.~\eqref{eq:lm_weak_cm} is nonlinear. We express $\mathbf{e}$ and its first variation in terms of the displacement gradient:
\begin{subequations}
    \begin{align}
        \eu(\uuu) = & \frac{1}{2}\left(\nabla\uuu + \nabla^{\mathrm{T}}\uuu + \nabla^{\mathrm{T}}\uuu\,\nabla\uuu\right), \\
        \delta\eu(\uuu,\www) = & \frac{1}{2}\left(\nabla\www + \nabla^{\mathrm{T}}\www + \nabla^{\mathrm{T}}\www\,\nabla\uuu + \nabla^{\mathrm{T}}\uuu\,\nabla\www \right).
    \end{align}
\end{subequations}
Let $\vepsuu(\auu)=(\nabla\auu+\nabla^{\mathrm{T}}\auu)/2$ and the nonlinear operator $\eu^{ns}(\auu,\buu) = (\nabla^{\mathrm{T}}\auu\nabla\buu+\nabla\auu\nabla^{\mathrm{T}}\buu)/2$. 
With these definitions:
\begin{align}
    \delta\eu:\cA:\eu = \vepsuu(\www):\cA:\vepsuu(\uuu) + \frac{1}{2} \vepsuu(\www):\cA:\eu^{ns}(\uuu,\uuu) +   \eu^{ns}(\www,\uuu):\cA:\vepsuu(\www) + \frac{1}{2} \eu^{ns}(\www,\uuu):\cA:\eu^{ns}(\uuu,\uuu).
\end{align}
The first terms is linear, the second and the third are quadratic and the last term is cubic. The power of internal forces term can then be rewritten as:
\begin{equation}
	\int_{\Omega_0} \delta\eu:\cA:\eu \,d\Omega_0 = \kuu(\uuu,\www) + \guu(\uuu,\uuu,\www) + \huu(\uuu,\uuu,\uuu,\www).
\end{equation}
By taking a triplet of displacement fields $\uuu_k$, $\uuu_l$, and $\uuu_m$ the three terms $\kuu$, $\guu$, and $\huu$ are expressed as:
\begin{subequations}\label{eq:decomposition}
\begin{align}
	\kuu(\uuu_k,\www) = & \int_{\Omega_0} \vepsuu(\www) : \mathcal{A} : \vepsuu(\uuu_k) \,d\Omega_0, 		\label{eq:decompositiona}\\
    \guu(\uuu_k,\uuu_l,\www) = & \frac{1}{2}\int_{\Omega_0} \vepsuu(\www):\mathcal{A}:\eu^{ns}(\uuu_k,\uuu_l) + \frac{1}{2}\eu^{ns}(\www,\uuu_k):\mathcal{A}:\vepsuu^{ns}(\uuu_l) +   \frac{1}{2}\eu^{ns}(\www,\uuu_l):\mathcal{A}:\vepsuu^{ns}(\uuu_k) \,d\Omega_0 ,\label{eq:decompositionb}\\
         \huu(\uuu_k,\uuu_l,\uuu_m,\www) = & \frac{1}{3}\int_{\Omega_0} \mathbf{e}^{ns}(\www,\uuu_k):\mathcal{A}:\mathbf{e}^{ns}(\uuu_l,\uuu_m) +   \mathbf{e}^{ns}(\www,\uuu_l):\mathcal{A}:\mathbf{e}^{ns}(\uuu_m,\uuu_k) +  \mathbf{e}^{ns}(\www,\uuu_m):\mathcal{A}:\mathbf{e}^{ns}(\uuu_k,\uuu_l) \,d\Omega_0,\label{eq:decompositionc}
	\end{align}
\end{subequations}
where Eq.~\eqref{eq:decompositiona} is linear with respect to the displacement, Eq.~\eqref{eq:decompositionb} is quadratic, and Eq.~\eqref{eq:decompositionc} is cubic. Upon finite element discretisation of Eqs.~\eqref{eq:decomposition}, one obtains:
\begin{subequations}
    \begin{align}
         \kuu_{h}(\uuu_k,\www) = &\,\UTu^{\mathrm{T}}\mathbf{K}\mathbf{U}_k,\\
        \guu_{h}(\uuu_k,\uuu_l,\www) = &\,\UTu^{\mathrm{T}}\mathbf{G}(\mathbf{U}_k, \mathbf{U}_l), \\
        \huu_{h}(\uuu_k,\uuu_l,\uuu_m,\www) =&\, \UTu^{\mathrm{T}}\mathbf{H}(\mathbf{U}_k, \mathbf{U}_l, \mathbf{U}_m),
    \end{align}
\end{subequations}
where $\UTu$ denotes the vector of nodal values for the test function and subscript $(\cdot)_{h}$ identifies approximated quantities. In the present work, the terms in Eq.~\eqref{eq:decomposition} are computed directly as proposed in~\cite{LazarusThomas2012} and in~\cite{Touze2014a} for MITC shell elements. 
Non-intrusive methods can be used as well, as for instance the STEP~\cite{muravyov,givois2019}.

\section{Physical System diagonalisation}
\label{sec:diag}

In this appendix we recall, for the sake of completeness, some well-known results regarding the
diagonalisation of the state-space formulation used in Eq.~\eqref{eq:damped_dyn_first}.  
Let the matrix $\Phiu$ collect column-wise the real-valued eigenvectors $\phiu_s$ 
of the linear part of Eq.~\eqref{eq:damped_dyn}. 
By setting $\Uu=\Phiu\uu$ and pre-multiplying Eq.~\eqref{eq:damped_dyn} 
by $\Phiu^{\mathrm{T}}$, one gets:
\begin{equation}\label{eq:dynappendix}
	\ddot{\uu} + \Omegau^2 \uu + \gu(\uu,\uu) + \hu(\uu,\uu,\uu) = \zerou,
\end{equation}
where $\Omegau$ is the diagonal matrix that stores the eigenfrequencies of the system, and $\gu(\uu,\uu)$, $\hu(\uu,\uu,\uu)$ are the modal quadratic and cubic nonlinearities:
\begin{subequations}
	\begin{align}
		&\Omegau^2 = \Phiu^{\mathrm{T}}\Ku\Phiu,\\
		&\gu(\uu,\uu) = \Phiu^{\mathrm{T}}\Gu(\Uu,\Uu),\\ 		&\hu(\uu,\uu,\uu) = \Phiu^{\mathrm{T}}\Hu(\Uu,\Uu,\Uu).
	\end{align}
\end{subequations}
Introducing the modal velocity $\vu=\dot{\uu}$, 
Eq.~\eqref{eq:dynappendix} is written as a system of first-order differential equations. 
By letting $\xu=\{ \vu,\uu \}$:
\begin{equation}\label{eq:modal_dyn_first}
	\dot{\xu} = \au{\xu} + \bu(\xu,\xu) + \cu(\xu,\xu,\xu),
\end{equation}
with:
\begin{subequations}
	\begin{align}
		 \au =&\, \left[\begin{array}{cc}
\zerou & -\Omegau^{2} \\
\Iu & \zerou
\end{array}\right], \\
		\bu(\xu,\xu) =&\, 
	\left\{\begin{array}{c}
		-\gu(\uu,\uu) \\
		\zerou
	\end{array}\right\},& \\\
		\cu(\xu,\xu,\xu) =&\, 
	\left\{\begin{array}{c}
		-\hu(\uu,\uu,\uu) \\
		\zerou
	\end{array}\right\}.
	\end{align}
\end{subequations}
Let us now introduce the matrix $\Ru$ defined as:
\begin{equation}
	\Ru = 
	\left[\begin{array}{cc}
		\img\Omegau & -\img\Omegau \\
		\Iu & \Iu
	\end{array}\right] ,
\end{equation}
with $\Iu$ the $\Nr\times \Nr$ identity matrix, and set $\xu=\Ru\pu$, with $\pu$ the $2\Nr$-dimensional vector of generalised coordinates. 
In particular one has:
\begin{subequations}
\label{eq:mod_map_real}
	\begin{align}
		\ur_s = &\, \pr_s + \pr_{s+\Nr}  \qquad  \forall\,s = 1,...,\Nr, \\
	\vr_s = &\, \img\omega_{s}\pr_s -\img\omega_{s}\pr_{s+\Nr} \qquad  \forall\,s = 1,...,\Nr.
	\end{align}
\end{subequations}
The linear operator $\au$ is diagonalised by $\Ru$.
Indeed, pre-multiplying Eq.~\eqref{eq:modal_dyn_first} by $\Ru^{-1}$, the following equation is obtained:
\begin{equation}
	\dot{\pu} = \Lambdau\pu + \Gammau(\pu,\pu) + \Deltau(\pu,\pu,\pu),
\end{equation}
with:
\begin{subequations}
	\begin{align}
		\Lambdau =&\,\Ru^{-1}\au\Ru, \\
		\Gammau(\pu,\pu) =&\, \Ru^{-1}\bu(\xu,\xu), \\ 
		\Deltau(\pu,\pu,\pu) =&\, \Ru^{-1}\cu(\xu,\xu,\xu).
	\end{align}
\end{subequations}
%
In this framework, $\Lambdau$ is a diagonal matrix  which reads
\begin{equation}
\label{eq:lambda}
	\Lambdau = 
	\left[\begin{array}{cc}
		\img\Omegau & \zerou \\
		\zerou & -\img\Omegau
	\end{array}\right] .
\end{equation}

\section{Higher Order Expansions}
\label{sec:hot}

In Secs.~\ref{sec:theory} and \ref{sec:ROM} the normal form approach has been introduced with focus on low order mapping and reduced dynamics. This choice is motivated by the simplicity of low order formulations, which in turn highlights the efficiency of the presented method for the development of reduced-order models. 
However, the method can be implemented in an algorithmic way for arbitrary order expansions. 
Indeed the general structure of homological equations are:
\begin{align}
	\left[\begin{array}{cc}
		\Mu & \zerou \\
		\zerou & -\Iu
	\end{array}\right] 
	\left\{\begin{array}{l}
		\dot{\Upsu}^{(n)} \\
		\dot{\Psiu}^{(n)}
	\end{array}\right\} + 
	\left[\begin{array}{cc}
		\zerou & \Ku \\
		\Iu & \zerou
	\end{array}\right]
	\left\{\begin{array}{l}
		\Upsu^{(n)} \\
		\Psiu^{(n)}
	\end{array}\right\} +  \left\{\begin{array}{l}
		\Gu^{(n)} \\
		\zerou
	\end{array}\right\}	+
	\left\{\begin{array}{l}
		\Hu^{(n)} \\
		\zerou
	\end{array}\right\}	
	= 
	\zerou,
\end{align}
with $\Gu^{(n)}$ defined as the sum of all $\Gu(\Psiu^{(a)},\Psiu^{(b)})$ such that $a+b=n$, and  $\Hu^{(n)}$ equal to the sum of all $\Hu(\Psiu^{(a)},\Psiu^{(b)},\Psiu^{(c)})$ such that $a+b+c=n$. Let $\Xiu^{(n)}=\Gu^{(n)}+\Hu^{(n)}$.
%
For the sake of completeness, we detail hereafter the equations needed up to order five.

\noindent First order:
\begin{subequations}
	\begin{align}
	&\lambda_s \Mu \Upsu^{(1)}_{s} + \Ku \Psiu^{(1)}_{s} = 0,
	\\
	&\lambda_s\Psiu^{(1)}_{s} = \Upsu^{(1)}_{s} .
	\end{align}
\end{subequations}
Second order:
\begin{subequations}
	\begin{gather}
\left[ \left(\lambda_k+\lambda_l \right)\Mu \Upsu_{kl}^{(2)} + \Mu \Upsu^{(1)}_{s}f^{(2)}_{skl} \right] +  \Ku \Psiu^{(2)}_{kl} + \Xiu^{(2)}_{kl} = 0.
	\\
	\left[ \left(\lambda_k+\lambda_l \right)\Psiu_{kl}^{(2)} + \Psiu^{(1)}_{s}f^{(2)}_{skl} \right]  = \Upsu^{(2)}_{kl}.
	\end{gather}
\end{subequations}
Third order:
\begin{subequations}
	\begin{gather}
	\left[ \left(\lambda_k+\lambda_l+\lambda_m \right)\Mu \Upsu_{klm}^{(3)} + \Mu \Upsu^{(1)}_{s}f^{(3)}_{sklm} \right] +  \left[ \Mu \Upsu^{(2)}_{am}f^{(2)}_{akl} + \Mu \Upsu^{(2)}_{ka}f^{(2)}_{alm} \right] +  \Ku \Psiu^{(3)}_{klm} + \Xiu^{(3)}_{klm} = 0.
	\\
	 \left[ \left(\lambda_k+\lambda_l+\lambda_m \right)\Psiu_{klm}^{(3)} + \Psiu^{(1)}_{s}f^{(3)}_{sklm} \right] +  \left[ \Psiu^{(2)}_{am}f^{(2)}_{akl} + \Psiu^{(2)}_{ka}f^{(2)}_{alm} \right]  = \Upsu^{(3)}_{klm} .
	\end{gather}
\end{subequations}
Fourth order:
\begin{subequations}
	\begin{align}
	& \left[ \left(\lambda_k+\lambda_l+\lambda_m+\lambda_n \right)\Mu \Upsu_{klmn}^{(4)} + \Mu \Upsu^{(1)}_{s}f^{(4)}_{sklmn} \right] + \nonumber \\	
	& \left[ \Mu \Upsu^{(2)}_{an}f^{(3)}_{aklm} + \Mu \Upsu^{(2)}_{ka}f^{(3)}_{almn}  \right] + \nonumber \\
	& \left[ \Mu \Upsu^{(3)}_{amn}f^{(2)}_{akl} + \Mu \Upsu^{(3)}_{kan}f^{(2)}_{alm} + \Mu \Upsu^{(3)}_{kla}f^{(2)}_{amn}  \right] + \nonumber \\
	& \Ku \Psiu^{(4)}_{klmn} + \Xiu^{(4)}_{klmn} = 0.
	\\
		& \left[ \left(\lambda_k+\lambda_l+\lambda_m+\lambda_n \right)\Psiu_{klmn}^{(4)} + \Psiu^{(1)}_{s}f^{(4)}_{sklmn} \right] + \nonumber \\	
	& \left[ \Psiu^{(2)}_{an}f^{(3)}_{aklm} + \Psiu^{(2)}_{ka}f^{(3)}_{almn}  \right] + \nonumber \\
	& \left[ \Psiu^{(3)}_{amn}f^{(2)}_{akl} + \Psiu^{(3)}_{kan}f^{(2)}_{alm} + \Psiu^{(3)}_{kla}f^{(2)}_{amn}  \right] \nonumber \\
	& = \Upsu^{(4)}_{klmn} .
	\end{align}
\end{subequations}
Fifth order:
\begin{subequations}
	\begin{align}
	& \left[ \left(\lambda_k+\lambda_l+\lambda_m+\lambda_n+\lambda_p \right)\Mu \Upsu_{klmnp}^{(5)} + \Mu \Upsu^{(1)}_{s}f^{(5)}_{sklmno} \right] + \nonumber \\	
	& \left[ \Mu \Upsu^{(2)}_{ap}f^{(4)}_{aklmn} + \Mu \Upsu^{(2)}_{ka}f^{(4)}_{almnp}  \right] + \nonumber \\
	& \left[ \Mu \Upsu^{(3)}_{anp}f^{(3)}_{aklm} + \Mu \Upsu^{(3)}_{kap}f^{(3)}_{almn} + \Mu \Upsu^{(3)}_{kla}f^{(3)}_{amnp}  \right] + \nonumber \\
	& \left[ \Mu \Upsu^{(4)}_{amnp}f^{(2)}_{akl} + \Mu \Upsu^{(4)}_{kanp}f^{(2)}_{alm} +
	 \Mu \Upsu^{(4)}_{klap}f^{(2)}_{amn} + \Mu \Upsu^{(4)}_{klma}f^{(2)}_{anp}  \right] + \nonumber \\
	& \Ku \Psiu^{(5)}_{klmnp} + \Xiu^{(5)}_{klmnp} = 0.
	\\
		& \left[ \left(\lambda_k+\lambda_l+\lambda_m+\lambda_n+\lambda_p \right)\Psiu_{klmnp}^{(5)} + \Psiu^{(1)}_{s}f^{(5)}_{sklmno} \right] + \nonumber \\	
	& \left[ \Psiu^{(2)}_{ap}f^{(4)}_{aklmn} + \Psiu^{(4)}_{ka}f^{(4)}_{almnp}  \right] + \nonumber \\
	& \left[ \Psiu^{(3)}_{anp}f^{(3)}_{aklm} + \Psiu^{(3)}_{kap}f^{(3)}_{almn} + \Psiu^{(3)}_{kla}f^{(3)}_{amnp}  \right] + \nonumber \\
	& \left[ \Psiu^{(4)}_{amnp}f^{(2)}_{akl} + \Psiu^{(4)}_{kanp}f^{(2)}_{alm}  \right] + \nonumber \\
	& \left[  \Psiu^{(4)}_{klap}f^{(2)}_{amn} + \Psiu^{(4)}_{klma}f^{(2)}_{anp}  \right] \nonumber \\
	& = \Upsu^{(5)}_{klmnp}.
	\end{align}
\end{subequations}
%
This set of equations highlights the general structure of the problems to be solved and gives the framework for computation of higher-order automated solutions; the implementation of the method for a generic order is out of the scope of this work and will be the subject of future contributions by the authors.

\section{Implementation Algorithm}
\label{sec:algo}

We provide herein some details of the implementation procedure, as schematically reported in Algorithm \ref{alg:implementation}.\\
The input parameters of the model are a discretisation of the system domain (MESH), together with boundary conditions (BCS), and material properties (PROP). 
Furthermore, the user should provide the list of master modes to be selected in the reduction, collected in the subset $\cZ^{(1/2)}$ following the notation used in the main text.

The algorithm starts by computing the mass $\Mu$ and stiffness $\Ku$ matrices (LinInt) as well as the eigenmodes $\Phiu_m$ and the eigenfrequencies $\Omegau_m$ (GenEig) of the selected master modes only. 
Since our basic assumption is that of a conservative vibratory system, all maps $\Psiu^{(2)}$ are real-valued and the reduced dynamics coefficients $\fuu^{(2)}$, $\fuu^{(3)}$ are purely imaginary. This implies that the formulation allows working on real double precision arithmetic even for the coefficients of the reduced dynamics by storing only the imaginary parts.\\
First-order reduced dynamics coefficients $\fuu^{(1)}$ are obtained as $\pm\img\Omegau$. First-order displacement mappings are the linear eigenmodes.\\
Second order reduced dynamics coefficients $\fuu^{(2)}$ and mappings $\Psiu^{(2)}$ are then computed by iterating over all $kl$ permutations. For $\nr$ master modes, then each index spans from $1$ to $2\nr$. For each permutation, first we check that $(\lambda_{k}+\lambda_{l})$ does not yield a resonance condition (CheckResonance2). Then, the matrix $\Du$ required to compute $\Psiu^{(2)}$ is assembled (AssemblyMapMatrix):
\begin{equation}
	\Du = \left[\begin{array}{cc}
 		\left[(\lambda_{k}+\lambda_{l})^2\Mu + \Ku \right] &  \left[\Mu\Phiu_{\mathrm{R}}\right]   \\
		\left[\Mu\Phiu_{\mathrm{R}}\right]^{\mathrm{T}} &  \zerou 
	\end{array}\right],
\end{equation}
where the two blocks $\Mu\Phiu_{\mathrm{R}}$ are used to impose mass-orthogonality between the mapping and the eigenmodes of the resonant monomials. The right hand side of Eq.~\eqref{eq:so_map}, {\em i.e.} $\Gu(\Psiu_k^{(1)},\Psiu_l^{(1)})$, is then assembled in a vector $\Lu$ of dimensions compatible with $\Du$ (IntegrateFNL2). By solving the resulting linear system:
\begin{align}\label{eq:full_matrix} 
	 \left\{\begin{array}{c}
		\Psiu^{(2)}_{kl} \\
		\left( \Lambdau_{\mathrm{R}}- \Lambdau_{\mathrm{R}+\Nr} \right) {\fuu}^{(2)}_{Rkl} 
	\end{array}\right\} = \Du^{-1}\Lu = 
	\left[\begin{array}{cc}
 		\left[(\lambda_{k}+\lambda_{l})^2\Mu + \Ku \right] &  \left[\Mu\Phiu_{\mathrm{R}}\right]   \\
		\left[\Mu\Phiu_{\mathrm{R}}\right]^{\mathrm{T}} &  \zerou 
	\end{array}\right]^{-1} 
	\left\{\begin{array}{c}
		-\Gu(\Psiu^{(1)}_{k},\Psiu^{(1)}_{l}) \\
		\zerou 
	\end{array}\right\},
\end{align}
with $\Lambdau_{\mathrm{R}}$ and $\Lambdau_{\mathrm{R}+\Nr}$ being the diagonal matrices with entries $\lambda_r$ and $\lambda_{r+\Nr}$ such that $\Psiu_{r}^{(1)}=\Psiu_{r+\Nr}^{(1)}=\phiu_{r}\in \Phiu_{\mathrm{R}}$, and $\fuu^{(2)}_{Rkl}$ is the vector of monomials $f^{(2)}_{rkl}$ that cannot be cancelled in the reduced dynamics. From Eq.~\ref{eq:full_matrix} both mappings and reduced dynamics coefficients are obtained.\\
For second order DNF, a second iteration is then performed to compute $\fuu^{(3)}$. Therefore all $klm$ permutations are spanned. For each permutation, $\Lu$ is assembled as (IntegrateFNL3):
\begin{align}
	\Lu = - \Mu\Upsu^{(2)}_{sm}f^{(2)}_{skl} - \Mu\Upsu^{(2)}_{ks}f^{(2)}_{slm} - \Xiu^{(3)}_{klm},
\end{align}
Finally, $f^{(3)}_{sklm}$ terms are computed by projecting $\Lu$ onto the master modes and by exploiting Eq.~\eqref{eq:ff_to_relation} (step "Project" in the algorithm).

\begin{algorithm}
\SetAlgoLined
\KwIn{MESH,BCS,PROP,$\cZ^{(1/2)}$}
\KwResult{$\Psiu^{(1)},\Psiu^{(2)},\fuu^{(1)},\fuu^{(2)},\fuu^{(3)}$}
 $\Ku$, $\Mu$ $\leftarrow$ LinInt(MESH,BCS,PROP) \;
 $\Phiu_m$, $\Omegau_m$ $\leftarrow$ GenEig($\Ku$, $\Mu$, $\cZ^{(1/2)}$) \;
 $\fuu^{(1)}$ $\leftarrow$ $\Omega_m$ \;
 $\Psiu^{(1)}$ $\leftarrow$ $\Phiu_m$ \;
 \For{$k\in\cZ$}{
    \For{$l\in\cZ$}{
       $\Phiu_{\mathrm{R}}$ $\leftarrow$ CheckResonance2($\Phiu_m$,$\Omegau_m$,k,l) \;
       $\Lu$ $\leftarrow$  IntegrateFNL2($\Psiu^{(1)}$, k,l) \;
       $\Du$ $\leftarrow$  AssemblyMapMatrix($\Ku$,$\Mu$,$\Phiu_{\mathrm{R}}$,$\Omegau_m$, k,l) \;
       $\Psiu^{(2)},{\fuu}^{(2)}$ $\leftarrow$ $\Du^{-1}\Lu$ \;
    }
 }
 \For{$k\in\cZ$}{
    \For{$l\in\cZ$}{
        \For{$m\in\cZ$}{
    	    $\Lu$ $\leftarrow$  IntegrateFNL3($\Psiu^{(1)}$,$\Psiu^{(2)}$,$\fuu^{(2)}$, k,l,m) \;
    	    $\fuu^{(3)}$ $\leftarrow$ Project($\phiu_m,\Lu$)\;
    	}
    }
 }
 \Return
 \caption{Implementation algorithm.}\label{alg:implementation}
\end{algorithm}



\section{Third Order Monomials for 1:2 Internally Resonant System}
\label{sec:second_order_ir_coeff}
The linear mapping relating complex to real-valued quantities reported in Eq.~\eqref{eq:z_to_rs_map} can be applied to reduced-order models with an arbitrary number of master modes. The case of a 1:2 internally resonant system is here considered as an illustration, 
since this case is explicitly investigated in the example of Section~\ref{subsec:archRES12}. 
Let us assume that the two modes satisfying the 1:2 relation are mode 1 and mode 2 for the sake
of simplicity. Therefore, model order reduction with the present method requires selecting as master coordinates $(\zr_1,\zr_{1+\Nr})$ and $(\zr_2,\zr_{2+\Nr})$. The reduced dynamics given in Eq.~\eqref{eq:red_dyn} when mapped to real-valued quantities using Eq.~\eqref{eq:z_to_rs_map} yields the following equations for $\dot{\ru}$:
\begin{subequations}\label{eq:rom_w_nasty_on_dotr}
	\begin{align}
		& 
		\dot{\rr}_{1} + \gr_{112} (
			4 \hat{\br}_{111}  \rr_1 \rr_2 \sr_1 +
			2 \hat{\br}_{122}  \rr_1^2 \sr_2
		)
		=  \sr_1, \\
		& 
		\dot{\rr}_{2} + \gr_{211} (
			4 \hat{\br}_{212} \rr_1 \rr_2 \sr_2 + 
			2 \hat{\br}_{212} \rr_1^2 \sr_1 + 
			2 \hat{\br}_{222} \rr_1^2 \sr_2
		)
		=  \sr_2.
	\end{align}
\end{subequations}
The same operation yields the following equations for $\dot{\su}$:
\begin{subequations}\label{eq:rom_w_nasty_on_dots}
	\begin{align}
		& 
		\dot{\sr}_{1} + \omega_1^{2}\rr_1 + (\gr_{112}+\gr_{121})\rr_1\rr_2 + \nonumber \\ 
		& 
		(\Ar_{1111}+\hr_{1111})\rr_1^{3} + 3(\Ar_{1122}+\hr_{1122})\rr_1\rr_2^{2} +  3(\Ar_{1112}+\hr_{1112})\rr_1^{2}\rr_2 + (\Ar_{1222}+\hr_{1222})\rr_2^{3} + \nonumber \\ 
		& 
		\Br_{1111}\rr_1{\sr}_1^{2} + \Br_{1211}\rr_2{\sr}_1^{2}  +   \Br_{1122}\rr_1 {\sr}_2^{2} + 2\Br_{1112}\rr_1{\sr}_1{\sr}_2 +  2\Br_{1212}\rr_2{\sr}_1{\sr}_2 + \Br_{1222}\rr_2{\sr}_2^{2}   + \nonumber \\
		& 
		-\gr_{112}(
		4\hat{\ar}_{111} - 4\hat{\br}_{111}\omega_1^2 +
		2\hat{\ar}_{122} - 2\hat{\br}_{122}\omega_2^2
		)\rr_1^{2}\rr_2 
		=  0, \\
		& 
		\dot{\sr}_{2} + \omega_2^{2}\rr_2 + \gr_{211}\rr_1^2 + \nonumber \\ 
		& 
		(\Ar_{2111}+\hr_{2111})\rr_1^{3} + 3(\Ar_{2122}+\hr_{2122})\rr_1\rr_2^{2} +  3(\Ar_{2112}+\hr_{2112})\rr_1^{2}\rr_2 + (\Ar_{2222}+\hr_{2222})\rr_2^{3} + \nonumber \\ 
		& 
		\Br_{2111}\rr_1{\sr}_1^{2} + \Br_{2211}\rr_2{\sr}_1^{2}  +  \Br_{2122}\rr_1 {\sr}_2^{2} + 2\Br_{2112}\rr_1{\sr}_1{\sr}_2 +  2\Br_{2212}\rr_2{\sr}_1{\sr}_2 + \Br_{2222}\rr_2{\sr}_2^{2}   + \nonumber \\ 
		& 
		-\gr_{211} (
			(2 \hat{\ar}_{212} - 2\hat{\br}_{212}\omega_1^2) \rr_1^3 + 
			(2 \hat{\ar}_{222} - 2\hat{\br}_{222}\omega_2^2) \rr_1^2 \rr_2 + 
			(4 \hat{\ar}_{212} - 4\hat{\br}_{212}\omega_2^2) \rr_1 \rr_2^2
		) =  0.
	\end{align}
\end{subequations}
An important remark compared to previous developments is that the normal velocity $\su$ is not equal to the time derivative of the normal displacement $\dot{\ru}$. This is a noticeable result that cannot be observed for real-valued reduced dynamics truncated at third order if no second order resonances between master modes are observed. By taking the derivatives of Eq.~\eqref{eq:rom_w_nasty_on_dotr} with substitution truncated at third-order:
\begin{subequations}\label{eq:rom_w_nasty_on_dottr}
	\begin{align}
		& 
		\ddot{\rr}_{1} + \gr_{112} (
			4 \hat{\br}_{111}  \rr_2 \sr_1^2
			+ (4 \hat{\br}_{111} + 4 \hat{\br}_{122}) \rr_1 \sr_1 \sr_2 
			- (4\hat{\br}_{111}\omega_1^2 + 2\hat{\br}_{122}\omega_2^2) \rr_1^2 \rr_2 
		)
		=  \dot{\sr}_1, \\
		& 
		\ddot{\rr}_{2} + \gr_{211} (			
			- 2 \hat{\br}_{212}\omega_1^2  	\rr_1^3
			- 2 \hat{\br}_{222}\omega_2^2 	\rr_1^2 \rr_2
			- 4 \hat{\br}_{212}\omega_2^2  	\rr_1 \rr_2^2 
			+ 4 \hat{\br}_{222} 				\rr_1 \sr_1 \sr_2 
			+ 4 \hat{\br}_{212} 				\rr_2 \sr_1 \sr_2
			+ 4 \hat{\br}_{212} 				\rr_1 \sr_2^2			
			+ 4 \hat{\br}_{212} 				\rr_1 \sr_1^2 			
		)
		=  \dot{\sr}_2.
	\end{align}
\end{subequations}
Summing Eqs.~\eqref{eq:rom_w_nasty_on_dottr} and Eqs.~\eqref{eq:rom_w_nasty_on_dots} and substituting $\su$ with $\dot{\ru}$: 
\begin{subequations}\label{eq:rom_w_nasty}
	\begin{align}
		& 
		\ddot{\rr}_{1} + \omega_1^{2}\rr_1 + (\gr_{112}+\gr_{121})\rr_1\rr_2 + \nonumber \\ 
		& 
		(\Ar_{1111}+\hr_{1111})\rr_1^{3} + (3\Ar_{1122}+3\hr_{1122})\rr_1\rr_2^{2} +  (3\Ar_{1112}+3\hr_{1112}+\mathrm{P}_{1112})\rr_1^{2}\rr_2 + (\Ar_{1222}+\hr_{1222})\rr_2^{3} + \nonumber \\ 
		& 
		\Br_{1111}\rr_1\dot{\rr}_1^{2} + (\Br_{1211}+\Qr_{1211})\rr_2\dot{\rr}_1^{2}  +   \Br_{1122}\rr_1 \dot{\rr}_2^{2} + (2\Br_{1112}+\Qr_{1112})\rr_1\dot{\rr}_1\dot{\rr}_2 +  2\Br_{1212}\rr_2\dot{\rr}_1\dot{\rr}_2 + \Br_{1222}\rr_2\dot{\rr}_2^{2}   =  0, \\
		& 
		\ddot{\rr}_{2} + \omega_2^{2}\rr_2 + \gr_{211}\rr_1^2 + \nonumber \\ 
		& 
		(\Ar_{2111}+\hr_{2111}+\mathrm{P}_{2111})\rr_1^{3} + (3\Ar_{2122}+3\hr_{2122}+\mathrm{P}_{2122})\rr_1\rr_2^{2} +  (3\Ar_{2112}+3\hr_{2112}+\mathrm{P}_{2112})\rr_1^{2}\rr_2 + (\Ar_{2222}+\hr_{2222})\rr_2^{3} + \nonumber \\ 
		& 
		(\Br_{2111}+\Qr_{2111})\rr_1\dot{\rr}_1^{2} + \Br_{2211}\rr_2\dot{\rr}_1^{2}  +  (\Br_{2122}+\Qr_{2122})\rr_1 \dot{\rr}_2^{2} + \nonumber \\
		& (2\Br_{2112}+\Qr_{2112})\rr_1\dot{\rr}_1\dot{\rr}_2 +  (2\Br_{2212}+\Qr_{2212})\rr_2\dot{\rr}_1\dot{\rr}_2 + \Br_{2222}\rr_2\dot{\rr}_2^{2}   =  0.
	\end{align}
\end{subequations}
This last set of equations represent the exact normal form of the system in real-valued formalism if truncated to second-order, with the presence of the only two resonant monomials.  Third order monomials are all kept in the reduced dynamics in accordance with the assumption of truncating the nonlinear mapping at second-order and the reduced dynamics at third. Writing the exact normal form up to the third-order would need application of third-order mappings that would cancel all non-resonnt third-order monomials in Eqs.~\eqref{eq:rom_w_nasty}. Second-order physical coefficients are obtained as usual with:
	\begin{align}
		\gr_{pkl} = \phiu_p^{\mathrm{T}}\Gu(\phiu_k,\phiu_l),
	\end{align}
which also applies to the cubic physical modal coupling coefficients $\hr_{pklm}$, reading
	\begin{align}
		\hr_{pklm}=\phiu_p^{\mathrm{T}}\Hu(\phiu_k,\phiu_l,\phiu_m).
	\end{align}
From $\Gu$ and mappings $\hat{\au}$, $\hat{\bu}$, one obtains the values for $\Ar_{pklm}$ and $\Br_{pklm}$ as
\begin{subequations}
	\begin{align}
		\Ar_{pklm}=&\,2\;\phiu_p^{\mathrm{T}}\Gu(\phiu_k,\hat{\au}_{lm}),\\
		\Br_{pklm}=&\,2\;\phiu_p^{\mathrm{T}}\Gu(\phiu_k,\hat{\bu}_{lm}).
	\end{align}
\end{subequations}
These equations follows the general computation guidelines for all the coefficients used in the case where no internal resonance is present, {\em i.e.} by direct application of the formula given in~\cite{artDNF2020,touze03-NNM}. The only difference being that the coefficients $\hat{a}_{klm}$ and $\hat{b}_{klm}$ corresponding to the 1:2 resonant monomials, are now equal to zero, which does not impact the general formula. Finally, new terms that arise due to the presence of non-zero second order terms in a 1:2 internally resonant system read:
\begin{subequations}
	\begin{align}
		&\mathrm{P}_{1112} = \gr_{112}\left( -4\hat{\ar}_{111}-2\hat{\ar}_{122}+8\hat{\br}_{111}\omega_1^{2}+4\hat{\br}_{122}\omega_2^{2} \right),\\
		&\mathrm{P}_{2111} = \gr_{211}\left( -2\hat{\ar}_{212}+4\hat{\br}_{212}\omega_1^{2} \right),\\
		&\mathrm{P}_{2112} = \gr_{211}\left( -2\hat{\ar}_{222}+4\hat{\br}_{222}\omega_2^{2} \right),\\
		&\mathrm{P}_{2122} = \gr_{211}\left( -4\hat{\ar}_{212}+8\hat{\br}_{212}\omega_2^{2} \right),\\
		& \Qr_{1211} = -4\gr_{112} \hat{\br}_{111},\\
		&\Qr_{1112} = -4\gr_{112}\left( \hat{\br}_{111}+\hat{\br}_{112} \right),\\
		&\Qr_{2111}=\Qr_{2122}= \Qr_{2112} = -4\gr_{211}\hat{\br}_{212},\\
		&\Qr_{2212} = -4\gr_{211}\hat{\br}_{222},
	\end{align}
\end{subequations}
with $\hat{\ar}_{m kl}$, $\hat{\br}_{m kl}$ obtained by the pre-multiplying the real-valued mappings $\hat{\au}_{kl}$, $\hat{\bu}_{kl}$ by $\Mu$ and subsequently projecting on $\phiu_m$. These new terms are of particular importance and are completely related to the existence of a 1:2 resonance. They are the consequence of the remaining second-order monomials in the normal form, that in turn creates new cubic terms that have been here computed and made explicit.

\section{On the identification of Non Trivial Resonance Conditions for Discrete Models}
\label{sec:r_cond}


The identification of non trivial internal resonance conditions in finite element systems needs to be decided by giving a tolerance since perfect integer ratios are not possible due to round-off errors. As detailed in Sec.~\ref{sec:theory}, the generic mapping $\Psiu^{(n)}$ can be estimated through solution of a linear system of the type:
\begin{equation}
	\Eu\Psiu^{(n)}=\Lu,
\end{equation}
with $\Eu$ a symmetric matrix which is singular in presence of resonance conditions. 
For trivial resonances, then $det(\Eu)$ is always zero. For non trivial resonances in floating point arithmetic the determinant of $\Eu$ is not zero, yet $\Eu$ is ill-conditioned. 
For each mapping term $\Psiu^{(n)}$, one has to invert a matrix of the form:
\begin{equation}
	\Eu = \left[ \sigmauu^{2}\Mu+\Ku \right],
\end{equation}
with $\sigmauu$ summation of $\Lambdau$ entries. Resonance condition implies that $\sigmauu^{2}$ is equal to minus any of the entries of $\Omegau^2$. Therefore, resonance conditions of a given order can be detected by taking the norm of the difference between $\sigmauu^{2}$ and the square of any eigenfrequency of the master modes $\omega_r$. If the value is below a given tolerance $\varepsilon$:
\begin{equation}
	|\sigmauu^{2} + \omega_r^{2}| < \varepsilon,
\end{equation}
then a resonance condition is assumed.

\end{document}